\documentclass[12pt,reqno,a4paper]{article}
\usepackage{geometry}
\usepackage[blocks]{authblk}
\usepackage{amsmath,amsthm,amstext,amscd,amssymb,euscript,url,mathrsfs,euscript,dsfont,calrsfs,mathtools,comment,relsize,xcolor}
\usepackage{hyperref}
\usepackage{amsmath,amsthm,amstext,amscd,amssymb,euscript,url}
\usepackage{mathrsfs}
\usepackage{calrsfs}
\usepackage{epsfig}
\usepackage[inline]{enumitem}
\usepackage{mathtools}
\usepackage[absolute,overlay]{textpos}
\usepackage{graphicx}
\usepackage{subfig}
\usepackage{tikz}
\usetikzlibrary{matrix}
\usetikzlibrary{arrows,positioning}
\usepackage{comment}

\usepackage{color}

\renewcommand{\P}{\mathbb{P}}

\newcommand{\R}{\mathbb R}
\newcommand{\N}{{\mathbb N}}

\newcommand{\E}{\mathbb E}

\renewcommand{\Pr}{\mathbb P}

\def\1{{\mathchoice {\rm 1\mskip-4mu l} {\rm 1\mskip-4mu l}
{\rm 1\mskip-4.5mu l} {\rm 1\mskip-5mu l}}}

\newtheorem{theorem}{{\small T}{\scriptsize HEOREM}}[section]
\newtheorem{corollary}{{\bf{\small C}{\scriptsize OROLLARY}}}[section]
\newtheorem{proposition}{{\bf{\small P}{\scriptsize ROPOSITION}}}[section]
\newtheorem{lemma}{{\bf{\small L}{\scriptsize EMMA}}}[section]
\newtheorem{remark}{{\bf{\small R}{\scriptsize EMARK}}}[section]
\newtheorem{definition}{{\bf{\small D}{\scriptsize EFINITION}}}[section]
\newtheorem{induction}{{\bf{\small I}{\scriptsize NDUCTIVE HYPOTHESIS}}}[section]
\newtheorem{example}{{\bf{\small E}{\scriptsize XAMPLE}}}[section]

\renewenvironment{proof}[1]
{\noindent{{\bf{\small{ P}{\scriptsize ROOF}}}.}\hspace{0.1cm} #1} {$\;\qed$\newline}

\newcommand{\beq}{\begin{eqnarray}}
\newcommand{\eeq}{\end{eqnarray}}

\newcommand{\ba}{\begin{align*}}
\newcommand{\ea}{\end{align*}}

\newcommand{\be}{\begin{equation}}
\newcommand{\ee}{\end{equation}}

\newcommand{\bl}{\begin{lemma}}
\newcommand{\el}{\end{lemma}}

\newcommand{\br}{\begin{remark}}
\newcommand{\er}{\end{remark}}

\newcommand{\bt}{\begin{theorem}}
\newcommand{\et}{\end{theorem}}

\newcommand{\bd}{\begin{definition}}
\newcommand{\ed}{\end{definition}}

\newcommand{\bind}{\begin{induction}}
\newcommand{\eind}{\end{induction}}

\newcommand{\bp}{\begin{proposition}}
\newcommand{\ep}{\end{proposition}}

\newcommand{\bc}{\begin{corollary}}
\newcommand{\ec}{\end{corollary}}

\newcommand{\bpr}{\begin{proof}}
\newcommand{\epr}{\end{proof}}

\newcommand{\bi}{\begin{itemize}}
\newcommand{\ei}{\end{itemize}}

\newcommand{\ben}{\begin{enumerate}}
\newcommand{\een}{\end{enumerate}}

\newcommand{\caG}{{\mathcal G}}

\newcommand{\caS}{{\mathcal S}}

\newmuskip\pFqmuskip

\newcommand\pFq[6][8]{%
	\begingroup 
	\pFqmuskip=#1mu\relax
	\mathcode`\,=\string"8000
	\begingroup\lccode`\~=`\,
	\lowercase{\endgroup\let~}\pFqcomma
	{}_{#2}F_{#3}{\left[\genfrac..{0pt}{}{#4}{#5};#6\right]}%
	\endgroup
}
\newcommand{\pFqcomma}{\mskip\pFqmuskip}

\newtheorem{asu}{Assumption}
\newcommand{\basu}{\begin{asu}}
\newcommand{\easu}{\end{asu}}

\DeclareMathOperator\supp{supp}

\begin{document}

\title{Sitewise Generator Decomposition and Non-Equilibrium Steady States for Generalized Zero-Range Processes}
\author{Berend van Tol}

\affil{Institute of Applied Mathematics, Delft University of Technology, Delft,
The Netherlands}

\date{\today}

\maketitle

\begingroup
\renewcommand{\thefootnote}{}
\footnotetext{Email: \texttt{B.T.vanTol@tudelft.nl}}
\addtocounter{footnote}{-1}
\endgroup
\enlargethispage{3\baselineskip}

\begin{abstract}
In this paper we study the generalized zero-range processes considered in \cite{greenblatt2006product, evans2006factorized, fajfrova2016invariant} with boundary driving. In these processes mass is moved between neighboring sites with rates depending only on the departure site. For a wide class of such models, the equilibrium stationary states are known to be of product form. We show that for such models, the boundary-driven generator can be decomposed into a sum of generators acting on single sites when integrating against (possibly inhomogeneous) products of the equilibrium marginals. We exploit this structure to derive new results for the harmonic models introduced in \cite{frassek2020non}. In particular, we prove that the propagation of mixture measures in \cite{giardina2025intertwining} is still valid on graphs with asymmetric edges. On the chain geometry with asymmetric edges, we identify a class of models for which the NESS is given precisely by the NESS of the symmetric harmonic model on the chain. In addition to these results for harmonic models, we also present a new class of asymmetric models on the chain for which the NESS can be found explicitly and exhibits shock behavior. The large deviations for the density profile corresponding to these NESS are stated and proved.
\end{abstract}

\section{Introduction}

Boundary-driven interacting particle systems play a crucial role in non-equilibrium statistical physics, where they provide microscopic models of mass and energy
transport far from equilibrium. Their stationary measures, usually referred to as non-equilibrium steady states (NESS), encode characteristic features of non-equilibrium systems, such as long-range correlations and stationary currents. In contrast with equilibrium Gibbs measures, such stationary measures are rarely available in explicit form. Notable exceptions are zero-range processes \cite{levine2005zero}, harmonic models \cite{frassek2020non, frassek2022exact}, and exclusion processes \cite{derrida1998exactly,derrida2002large,derrida2003exact,enaud2004large}.\\

In recent years, it has become clear that the NESS of these three classes of models share deep structural connections. It was shown that the NESS of the harmonic models can be interpreted as stochastic mixtures of inhomogeneous product measures \cite{carinci2025large}, with mixing laws given by ordered Dirichlet distributions. In other words, these NESS have the form

\be
\mu(\eta) = \E[\otimes_{i\in V} \nu_{\Theta_i}(\eta_i)],
\ee

where the expectation is taken with respect to the random parameters $\Theta_i$ and the one-site measures $\nu_{\theta}$ are the equilibrium marginals with mean $\theta$.  For the harmonic models on the chain $1,2, \ldots,N$, the random parameters then have an ordered Dirichlet law.

Soon after, it was shown in \cite{giardina2025intertwining} that the stochastic-mixture structure of these NESS arises through a hidden process that describes the evolution of the harmonic models by evolving the inhomogeneous random parameters $\Theta_i$ instead of the original interacting particle system. In the same paper, the NESS of the harmonic models were recovered by determining the stationary distribution of these random parameters. Remarkably, this can be done by finding the stationary distribution of a system containing only a single site coupled to two reservoirs, because the generator of the hidden process naturally decomposes into generators acting on the individual sites of the chain. Using the spatial Markov structure of the Dirichlet distribution, the stationary distribution of the one-site hidden process can then be extended to a stationary distribution on the whole chain.

In a subsequent paper \cite{gabrielli2026markovbridges}, stochastic-mixture measures whose random parameters are ordered and have a spatial Markov structure were identified precisely as measures of matrix-product ansatz form. The matrix product ansatz was, in turn, used by Derrida and collaborators to determine the NESS of TASEP and other exclusion-type processes. For the harmonic NESS, a matrix-product representation, along with several other representations, was provided in \cite{frassek2026representations}.\\

In this paper, we show that the single-site decomposition for the generator of the hidden process associated with harmonic models is a manifestation of a much more general decomposition principle. More precisely, we show that, for a wide class of processes that includes the harmonic and zero-range models, the generator naturally decomposes into single-site generators after integration with respect to an inhomogeneous product measure with equilibrium marginals. This decomposition provides a general strategy for extending the NESS of a one-site system to the NESS on a larger graph, such as the chain $1,2, \ldots,N$, using spatial Markov structure. This is analogous to the way in which the stationary state of the harmonic models was recovered in \cite{giardina2025intertwining}. The distinction is that, in this more general setting, the existence of an associated hidden process is not required.

The boundary-driven processes to which the decomposition theorem applies form a subclass of the generalized zero-range processes. They are precisely of the form considered in \cite{greenblatt2006product,evans2006factorized,fajfrova2016invariant}, for which the stationary states of the closed system are homogeneous product measures. In this paper, a generalized zero-range process is a misanthrope process in which particles move independently of the configuration at the target site, as in a zero-range process, but multiple particles may be displaced at once. Related reversible and gradient mass-exchange models have been classified in \cite{capanna2024mass}, partly with a view toward computing stationary large deviations through macroscopic fluctuation theory.\\

We exploit this decomposition to obtain new results for harmonic processes and to construct novel models whose NESS have a stochastic-mixture structure with an underlying spatial Markov property. Our purpose is to illustrate how these stochastic-mixture measures with an underlying Markovian structure, which correspond to the matrix product ansatz, arise naturally as stationary measures of boundary-driven generalized zero-range processes. All the new examples presented here are asymmetric. We generalize the intertwining and propagation of mixtures in \cite{giardina2025intertwining} to harmonic models with asymmetric oriented edge weights. By using the decomposition principle, we bypass the direct use of duality, allowing us to make the oriented edge weights asymmetric. We also identify models on the chain for which the random parameters take only two values, namely the left and right reservoir parameters. In these models, mass moves to the left and right at different rates. This asymmetry is stronger than a difference between the leftward and rightward oriented weights. The rates of leftward and rightward mass transport have genuinely different forms. Because the random parameters take only two values, the models exhibit a shock phenomenon. We provide a more detailed overview of the main new results and their relation to the existing literature in Section \ref{subsec:summary_main_results}.

\subsection{Summary of the main new results} \label{subsec:summary_main_results}

The main new contributions of this paper and their relation to the literature are the following.

\begin{itemize}
    \item[\textbf{1.}] \textbf{Decomposition principle:} We prove the single-site decomposition in Theorem
    \ref{thm:decomposition}.  For generalized zero-range processes admitting
    the equilibrium marginals $(\nu_q)_q$, integration against an
    inhomogeneous product of these marginals turns the full generator into a
    sum of generators acting on individual sites, with neighboring parameters
    playing the role of reservoir parameters.  We then show in Theorem 
    \ref{thm:markov-mixture-stationary} that compatible one-site stationary
    mixing laws yield a stationary measure for the full process through a
    Markov extension.  Product stationary measures for generalized
    zero-range and mass migration processes were characterized in
    \cite{greenblatt2006product,evans2006factorized,fajfrova2016invariant},
    while spatially Markov mixing measures and their relation with Markov
    bridges and matrix products are studied in
    \cite{redig2025non,gabrielli2026markovbridges}.  The new point here is the
    generator-level reduction to one-site boundary problems and the resulting
    stationarity criterion.

    \item[\textbf{2.}] \textbf{Propagation of mixtures for harmonic models with asymmetric oriented edge weights:} We extend the propagation-of-mixtures result for harmonic models to
    asymmetric oriented edge weights in Theorem
    \ref{thm:asymmetric-harmonic-intertwining}.  The symmetric result of
    \cite{giardina2025intertwining} was obtained using duality.  Our proof
    instead uses the single-site decomposition and applies when, at every
    vertex, the sum of the incoming oriented edge weights does not exceed the sum of the outgoing
    oriented edge weights.

    \item[\textbf{3.}] \textbf{Explicit NESS for harmonic models with asymmetric oriented edge weights:} For a chain whose oriented edge weights depend only on the destination
    site, we identify an asymmetric family with exactly the same NESS as the
    symmetric harmonic chain. See Theorem
    \ref{thm:destination-weighted-harmonic-ness}. The symmetric NESS and its
    mixture representations are known from
    \cite{frassek2022exact,carinci2024solvable,carinci2025large}.  Explicit
    stationary states are also known for the different, $q$-deformed
    multiparticle asymmetric diffusion model
    \cite{frassek2024madminess}. 

    \item[\textbf{4.}] \textbf{New models with discrete mixing measure:} We construct asymmetric boundary-driven generalized zero-range processes whose NESS is an explicit mixture of product measures with random parameters concentrating on the values of the left and right reservoir parameters. See Proposition
    \ref{prop:asymmetric-shock-mixture}.  Theorem
    \ref{thm:asymmetric-shock-tail-factorization} characterizes the admissible
    rates, and
    Corollary \ref{cor:generalized-factorial-asymmetric-shock-family} gives a
    broad class expressed through Fonten\'e--Ward generalized factorials. The discrete nature of the random parameters implies that the NESS abruptly changes at a single, random, microscopic point (the ``shock''). Conditioned on the position of this microscopic point, the NESS has product structure, with different parameters for the one-site marginals on either side.
    Earlier zero-range shock results concern product measures whose shock
    position is also a single microscopic point and evolves dynamically
    \cite{balazs2010shock,balazs2019qzerorange}. Contrary to what is done in this paper, these results concern infinite systems and do not allow multiple particles to jump at once.

    \item[\textbf{5.}] \textbf{LDP for models with discrete mixing structure:} In Theorem \ref{thm:asymmetric-ness-profile-ldp} we prove an LDP for the density profile of the newly introduced models with discrete random parameters. More precisely, we find a joint large-deviation principle for the macroscopic shock position and the
    empirical density profile, and obtain the density-profile rate function by
    contraction.  Stationary empirical-density large deviations have
    previously been obtained for an ordinary boundary-driven zero-range
    process with long jumps in \cite{bernardin2022long}, and for the open
    harmonic process in \cite{carinci2025large}.  We remark that our result also fits into the broader Markov-bridge and matrix-product large-deviation framework of \cite{gabrielli2026matrixldp}. The
    contribution here is an explicit rate function for the new shock-mixture NESS.
\end{itemize}

\subsection{Organization of this paper}

This paper consists of three parts. In the first part, Section \ref{sec:gzrp_decomposition}, we introduce the models of interest (Section \ref{subsec:def_generalized_zero_range}), state and prove the decomposition principle (Section \ref{subsec:single-site_decomposition}), and show how mixture measures with a Markovian mixing structure arise naturally as NESS (Section \ref{subsec:stat_meas_Markov_mixture_structure}). In the second part, Section \ref{sec:Harmonic_models}, we apply the results of Section \ref{sec:gzrp_decomposition} to the harmonic models. We prove the propagation of mixtures for asymmetric oriented edge weights (Section \ref{subsec:prop_mixtures}) and identify models with asymmetric oriented edge weights that share the NESS of the symmetric harmonic models (Section \ref{subsec:asymmetric-harmonic-explicit-ness}). In the third part, we treat the new models with discrete mixing structure. The models are constructed in Section \ref{subsec:asym_models}, and the large deviations for the density profile are obtained in Section \ref{subsec:ldp}.

\section{Generalized zero-range processes and single-site decomposition} \label{sec:gzrp_decomposition}

\subsection{Definition of the generalized zero-range processes} \label{subsec:def_generalized_zero_range}

The processes we consider live on a finite simple graph $G=(V,E)$, where particles may be exchanged between two sites $i,j \in V$ if the edge $\{i,j\} \in E$. The state space is $\N_0^V$, where we write $\N_0 = \{0,1,2, \ldots \}$ for the natural numbers including zero. The notation $\N:=\{1,2,\ldots\}$ is used for the natural numbers not including zero. We denote a single configuration by $\eta$, i.e.
\begin{equation}
\eta=(\eta_i)_{i\in V}\in \N_0^V,
\label{eq:configuration-space}
\end{equation}
where $\eta_i$ is the number of particles at site $i$. For $i\in V$, we use $\mathbf e_i$ to denote the configuration containing one particle at site $i$ and no particles elsewhere.

\subsubsection{Rates in the bulk}
The interactions between two neighboring sites $i,j \in V$ are described by the following conservative generator,

\begin{align}
(L_{\{i,j\}}f)(\eta)
={}&
\sum_{k=1}^{\eta_i}
M_{i,j}(k)
\frac{\phi(\eta_i)}{\phi(\eta_i-k)}
\bigl[f(\eta-k\mathbf e_i+k\mathbf e_j)-f(\eta)\bigr]
\nonumber\\
&+
\sum_{k=1}^{\eta_j}
M_{j,i}(k)
\frac{\phi(\eta_j)}{\phi(\eta_j-k)}
\bigl[f(\eta-k\mathbf e_j+k\mathbf e_i)-f(\eta)\bigr].
\label{eq:bulk-generator}
\end{align}

Here,

\begin{itemize}
    \item [1.] $M_{i,j}$ and $M_{j,i}$ are non-negative functions from $\N = \{1,2,\ldots\} $ to $ [0,\infty)$. We refer to $M_{i,j}$ as the \emph{mass function} from $i$ to $j$, since it depends only on the transferred mass (number of particles) $k$.
    \item[2.] $\phi$ is a function from $\N_0 = \{0,1,2, \ldots\}$ to $(0, \infty)$. The function $\phi$ will be referred to as the \emph{potential function}: if $U(n):=-\log\phi(n)$, then 
    
    \begin{equation}
    \frac{\phi(n)}{\phi(n-k)}
    =
    \exp\bigl\{U(n-k)-U(n)\bigr\}.
    \label{eq:potential-difference}
    \end{equation}
    That is, the rate for moving $k$ particles away from a site is proportional to the exponential of the difference in $U$. Hence we interpret $U$ as a potential and refer, with slight abuse of terminology, to $\phi$ as the potential function.
\end{itemize}

The first sum in \eqref{eq:bulk-generator} corresponds to batches of $k$ particles moving from $i$ to $j$, and the second describes the reverse transition. In both directions, the rate depends on the occupation of the departure site but not on the occupation of the arrival site. Suppose that the mass function is symmetric, i.e. for all $\{i,j\} \in E$ we have $M_{i,j} = M_{j,i}$. Without boundary driving, this process has one parameter family of stationary homogeneous product measures

\be
(\otimes_{i \in V} \nu_{q})(\eta) = \prod_{i \in V} \nu_q(\eta_i), 
\ee

with single-site marginals $\nu_q$ as defined in Definition \ref{def:equilibrium-single-site-measures} below. The stationarity is the content of Proposition \ref{prop:detailed_balance_equi_marginal}, and is a consequence of detailed balance.

\begin{definition}[Equilibrium single-site measures]\label{def:equilibrium-single-site-measures}
For $q\in\R_0^+$ such that the normalizing constant below is finite, define
\begin{equation}
Z(q):=\sum_{m\in\N_0}\frac{q^m}{\phi(m)},
\qquad
\nu_q(n):=\frac{1}{Z(q)}\frac{q^n}{\phi(n)},
\quad n\in\N_0.
\label{eq:equilibrium-marginals}
\end{equation}
\end{definition}
The family $(\nu_q)_q$ depends only on the potential function $\phi$. 

\bp \label{prop:detailed_balance_equi_marginal}
Fix $q$ such that $Z(q) < \infty$ and assume that for all $\{i,j\} \in E$ we have $M_{i,j} = M_{j,i}$. Then, the homogeneous product measure $(\otimes_{i \in V} \nu_{q})(\eta)$ is stationary for the process with generator

\be \label{eq:closed_generator}
\sum_{\{i,j\} \in E} L_{\{i,j\}}.
\ee

\ep

\bpr
We verify the stationarity of $\otimes_{i \in V} \nu_q (\eta)$ using detailed balance. For fixed sites $i,j$ and all $k \geq 1$, $\eta_i, \eta_j \in \N$,

\begin{align} \label{eq:detailed_balance_closed}
&M_{i,j}(k) \frac{\phi(\eta_i)}{\phi(\eta_i - k)} \left(\otimes_{l \in V} \nu_q \right) (\eta) \\
&\qquad= M_{i,j}(k) \frac{\phi(\eta_i)}{\phi(\eta_i - k)} \left( \frac{\phi(\eta_i-k)}{\phi(\eta_i)} q^{k} \cdot \frac{\phi(\eta_j+k)} {\phi(\eta_j)} q^{-k} \left(\otimes_{l \in V} \nu_q \right) (\eta -k\mathbf e_i+k\mathbf e_j) \right)\nonumber\\
&\qquad= M_{j,i}(k)\frac{\phi(\eta_j+k)} {\phi(\eta_j)} \left(\otimes_{l \in V} \nu_q\right) (\eta -k\mathbf e_i+k\mathbf e_j). \nonumber
\end{align}

For the first equality we use the structure of the homogeneous product measure, and for the second equality we use the symmetry assumption on the mass functions. Since \eqref{eq:detailed_balance_closed} is satisfied for arbitrary $\{i,j\} \in E$, $k \in \N$, and $\eta_i, \eta_j \in \N_0$, detailed balance holds and $\otimes_{l \in V} \nu_q$ is stationary for the generator in \eqref{eq:closed_generator}.
\epr

\subsubsection{Boundary reservoirs}

In order to couple the bulk dynamics to external reservoirs, we associate with every $i\in V$ a set $V_g(i)$ of \emph{ghost sites} connected to $i$, and write
\begin{equation}
V_g:=\bigcup_{i\in V}V_g(i),
\qquad
\overline V:=V\cup V_g.
\label{eq:extended-vertex-set}
\end{equation}
With every ghost site $j\in V_g$, we associate a fixed boundary parameter $p_j^{\mathrm g}\in\R_0^+$, collected in
\begin{equation}
p^{\mathrm g}:=(p_j^{\mathrm g})_{j\in V_g}\in(\R_0^+)^{V_g}.
\label{eq:fixed-ghost-profile}
\end{equation}

Each ghost site $j$ represents a reservoir and should be interpreted as a site which is unaffected by exchanges with the system on the graph $G$. That is, its distribution is kept at $\nu_{p^g_j}$ throughout the time evolution of the system on the graph. The precise interpretation of the distribution being kept at $\nu_{p^g_j}$ is explained in Remark \ref{rmk:thermalized-boundary}.\\

Suppose that a ghost site $j \in V_g(i)$ associated to $i \in V$ had parameter $p^g_j = q$. Then the insertion and absorption of particles associated to the ghost site $j$ is modeled by the generator

\begin{align}
(L_i^{q}f)(\eta)
={}&
\sum_{k=1}^{\infty}
M_{j,i}(k)q^k
\bigl[f(\eta+k\mathbf e_i)-f(\eta)\bigr]
\nonumber\\
&+
\sum_{k=1}^{\eta_i}
M_{i,j}(k)
\frac{\phi(\eta_i)}{\phi(\eta_i-k)}
\bigl[f(\eta-k\mathbf e_i)-f(\eta)\bigr].
\label{eq:one-boundary-generator}
\end{align}

Here, the functions $M_{i,j}, M_{j,i}$ play the same role as the mass functions in the bulk dynamics in \eqref{eq:bulk-generator}. The first term corresponds to a batch of $k$ particles being inserted at site $i$, and the second corresponds to a batch of $k$ particles being removed from $i$.

\begin{remark}[Thermalized boundary conditions]\label{rmk:thermalized-boundary}
The insertion term in \eqref{eq:one-boundary-generator} is obtained by averaging a bulk transfer from a ghost site over $\nu_{p_j^{\mathrm g}}$. That is, if ghost site $j$ has occupation $\eta_j$ distributed according to $\nu_{p_j^{\mathrm g}}$, then a batch of $k$ particles is transferred from $j$ to $i$ at rate
\begin{equation}
M_{j,i}(k)\frac{\phi(\eta_j)}{\phi(\eta_j-k)} \1(\eta_j \geq k).
\label{eq:conditional-ghost-rate}
\end{equation}
Averaging over $\eta_j$ gives
\begin{align}
\E_{\nu_{p_j^{\mathrm g}}}\left[M_{j,i}(k)\frac{\phi(\eta_j)}{\phi(\eta_j-k)} \1(\eta_j \geq k)\right] &= M_{j,i}(k)
\sum_{n=k}^{\infty}
\frac{\phi(n)}{\phi(n-k)}
\nu_{p_j^{\mathrm g}}(n)
\nonumber\\
&=
\frac{M_{j,i}(k)}{Z(p_j^{\mathrm g})}
\sum_{n=k}^{\infty}
\frac{(p_j^{\mathrm g})^n}{\phi(n-k)}\nonumber\\
&=
M_{j,i}(k)(p_j^{\mathrm g})^k.
\label{eq:thermalized-insertion-rate}
\end{align}
The final equality follows by substituting $m=\eta_j-k$. Hence \eqref{eq:one-boundary-generator} describes thermalized boundary conditions.
\end{remark}

For later use we assign to each site $i \in V$ a parameter $p_i \in \R_0^+$ and denote

\begin{equation}
p:=(p_j)_{j\in V}\in (\R_0^+)^V.
\label{eq:bulk-parameter-profile}
\end{equation}

We will refer to $p$ as the \emph{parameter profile}. Furthermore, we introduce the following notation and conventions:

\begin{itemize}
    \item[1.] For $i\in V$ and $j\in\overline V$, we write $j\sim i$ whenever either $\{i,j\}\in E$ or $j\in V_g(i)$.
    \item[2.] For ghost site $j \in V_g$ we make the convention that $p$ with subscript $j$ is equal to $p^g_j$. That is $p_j := p^g_j$. Thus, $p_j$ denotes the parameter associated with any $j\in\overline V$, while the vector $p$ itself contains only the parameters associated with physical sites.
    \item[3.] For $S\subset V$, define its exterior neighborhood and the corresponding parameter family by
    \begin{equation}
    \partial S
    :=
    \bigl\{
    j\in\overline V\setminus S:
    j\sim i\ \text{for some }i\in S
    \bigr\},
    \label{eq:exterior-neighborhood}
    \end{equation}
    and
    \begin{equation}
    \partial p(S):=(p_j)_{j\in\partial S},
    \qquad
    \partial p(i):=\partial p(\{i\}).
    \label{eq:exterior-parameter-profile}
    \end{equation}
    Here $\partial p(S)$ contains variable parameters for physical neighbors outside $S$ and fixed parameters for ghost neighbors.
\end{itemize}

\begin{definition}[Admissible parameter domain]
\label{def:admissible-parameter-domain}
Define
\begin{equation}
\mathsf D
:=
\left\{
q\in\R_0^+:
Z(q)<\infty
\ \text{and}\
\sum_{k=1}^{\infty}
\bigl(M_{i,j}(k)+M_{j,i}(k)\bigr)q^k<\infty
\ \text{for every }i\in V\text{ and }j\sim i
\right\}.
\label{eq:admissible-parameter-domain}
\end{equation}
Here, \(\R_0^+:=[0,\infty)\). The set \(\mathsf D\) is an interval containing zero.
We call a bulk profile admissible if \(p\in\mathsf D^V\) and assume throughout that the fixed ghost profile satisfies \(p^{\mathrm g}\in\mathsf D^{V_g}\).
\end{definition}

In Definition \ref{def:boundary-driven-generator} we state the definition of the generator for the whole boundary-driven system on graph $G$ as well as the definition of the generator for a boundary-driven subsystem on $S \subset V$. The boundary parameters of the latter are given by $p$.

\begin{definition}[Boundary-driven generalized zero-range process]\label{def:boundary-driven-generator}
For $S\subset V$ and an admissible parameter profile $p\in\mathsf D^V$, the generator on $S$, with every neighbor outside $S$ replaced by its corresponding reservoir, is
\begin{equation}
L_S^{\partial p(S)}
:=
\sum_{\substack{\{i,j\}\in E\\i,j\in S}}
L_{\{i,j\}}
+
\sum_{i\in S}
\sum_{\substack{j\sim i\\j\notin S}}
L_i^{p_j}.
\label{eq:subsystem-generator}
\end{equation}
In particular, the generator of the whole boundary-driven particle system is
\begin{equation}
L_V^{\partial p(V)}
=
\sum_{\{i,j\}\in E}L_{\{i,j\}}
+
\sum_{i\in V}\sum_{j\in V_g(i)}L_i^{p_j}.
\label{eq:full-generator}
\end{equation}

Here, the generators $L_{\{i,j\}}$ and $L_{i}^{p_j}$ are as defined in \eqref{eq:bulk-generator} and \eqref{eq:one-boundary-generator}, respectively.
\end{definition}

The first sum in \eqref{eq:subsystem-generator} contains the conservative bulk dynamics along edges internal to $S$, while the second corresponds to replacing every edge leaving $S$ by a one-site reservoir interaction. In the case where the edge leads to another physical site $j\in V/S$ outside of $S$, the boundary parameter is $p_j$. For $S=V$, only the ghost sites lie outside the system. For $S=\{i\}$, every neighbor of $i$ is replaced by a reservoir and
\begin{equation}
L_{\{i\}}^{\partial p(i)}
=
\sum_{j\sim i}L_i^{p_j}
\label{eq:single-site-generator}
\end{equation}
acts only on site $i$.

\subsection{Single-site decomposition} \label{subsec:single-site_decomposition}

In this section we state and prove the result that underlies the remainder of the paper. This result, namely Theorem \ref{thm:decomposition}, allows us to decompose the whole boundary-driven process on graph $G$ into boundary-driven processes acting only on single sites.\\

Theorem \ref{thm:decomposition} involves integration against the inhomogeneous product measure on $\N_0^V$ with one-site marginals $\nu_{p_i}$, $i \in V$, which we denote

\begin{equation}
\bigotimes_{i\in V}\nu_{p_i}.
\label{eq:inhomogeneous-product-measure}
\end{equation}

Accordingly, for a function $f:\N_0^V\to\R$, we define
\begin{equation}
(\caG f)(p)
:=
\sum_{\eta\in\N_0^V}
\left(\bigotimes_{i\in V}\nu_{p_i}\right)(\eta)f(\eta).
\label{eq:G-transform}
\end{equation}

Thus, $(\caG f)(p)$ is the expectation of $f$ under the generally inhomogeneous product measure determined by $p\in\mathsf D^V$. We stress that the fixed ghost site parameters $p^{\mathrm g}$ are not included in the argument of $\caG f$.\\

The content of Theorem \ref{thm:decomposition} can be explained in words as follows. Although the original process contains interactions along the edges of $G$, after applying $\caG$ its generator can be replaced by a sum of one-site generators. In each one-site generator, a physical neighbor $j$ enters only through its parameter $p_j$, exactly as if it were a reservoir.

\begin{theorem}[Single-site decomposition]\label{thm:decomposition}
For every bulk parameter profile $p\in\mathsf D^V$ for which the expressions below are finite, and every function $f:\N_0^V\to\R$ in the domain of the generators,
\begin{equation}
\left(
\caG\left[L_V^{\partial p(V)}f\right]
\right)(p)
=
\left(
\caG\left[
\sum_{i\in V}L_{\{i\}}^{\partial p(i)}f
\right]
\right)(p) = 
\sum_{i\in V} \left(\caG \left[L_{\{i\}}^{\partial p(i)}  f \right] \right) (p).
\label{eq:decomposition}
\end{equation}
\end{theorem}

\begin{figure}[ht]
\centering
\begin{tikzpicture}[
    site/.style={circle,draw,thick,minimum size=8mm,inner sep=0pt},
    reservoir/.style={circle,draw,dashed,minimum size=8mm,inner sep=0pt},
    every node/.style={font=\small},
    >=stealth
]
\node at (6,0.9) {$L_V^{\partial p(V)}$};
\node[reservoir] (gl) at (0.6,0) {$p_{\ell}^{\mathrm g}$};
\node[site] (x1) at (3,0) {$1$};
\node[site] (x2) at (6,0) {$2$};
\node[site] (x3) at (9,0) {$3$};
\node[reservoir] (gr) at (11.4,0) {$p_r^{\mathrm g}$};
\draw[<->,thick] (gl) -- (x1);
\draw[<->,thick] (x1) -- (x2);
\draw[<->,thick] (x2) -- (x3);
\draw[<->,thick] (x3) -- (gr);

\node at (6,-1.25) {$\Downarrow\quad\text{after applying }\caG$};

\draw[gray,rounded corners] (0.1,-3.75) rectangle (3.5,-1.75);
\node at (1.8,-2.25) {$L_{\{1\}}^{\partial p(1)}$};
\node[reservoir] (a0) at (0.7,-3.2) {$p_{\ell}^{\mathrm g}$};
\node[site] (a1) at (1.8,-3.2) {$1$};
\node[reservoir] (a2) at (2.9,-3.2) {$p_2$};
\draw[<->] (a0) -- (a1);
\draw[<->] (a1) -- (a2);

\draw[gray,rounded corners] (4.3,-3.75) rectangle (7.7,-1.75);
\node at (6,-2.25) {$L_{\{2\}}^{\partial p(2)}$};
\node[reservoir] (b1) at (4.9,-3.2) {$p_1$};
\node[site] (b2) at (6,-3.2) {$2$};
\node[reservoir] (b3) at (7.1,-3.2) {$p_3$};
\draw[<->] (b1) -- (b2);
\draw[<->] (b2) -- (b3);

\draw[gray,rounded corners] (8.5,-3.75) rectangle (11.9,-1.75);
\node at (10.2,-2.25) {$L_{\{3\}}^{\partial p(3)}$};
\node[reservoir] (c2) at (9.1,-3.2) {$p_2$};
\node[site] (c3) at (10.2,-3.2) {$3$};
\node[reservoir] (cr) at (11.3,-3.2) {$p_r^{\mathrm g}$};
\draw[<->] (c2) -- (c3);
\draw[<->] (c3) -- (cr);

\node at (6,-4.5) {$\displaystyle\sum_{i=1}^{3}L_{\{i\}}^{\partial p(i)}$};
\end{tikzpicture}
\caption{The single-site decomposition for a three-site chain. Solid circles are physical sites and dashed circles are reservoirs. The ghost parameters $p_{\ell}^{\mathrm g}$ and $p_r^{\mathrm g}$ are fixed, whereas $p_1,p_2,p_3$ are coming from the argument $p$ of $(\caG L^{\partial p (V)}_V f )(p)$. The separated boxes in the lower panel emphasize that the three generators act on individual sites.}
\label{fig:single-site-decomposition}
\end{figure}
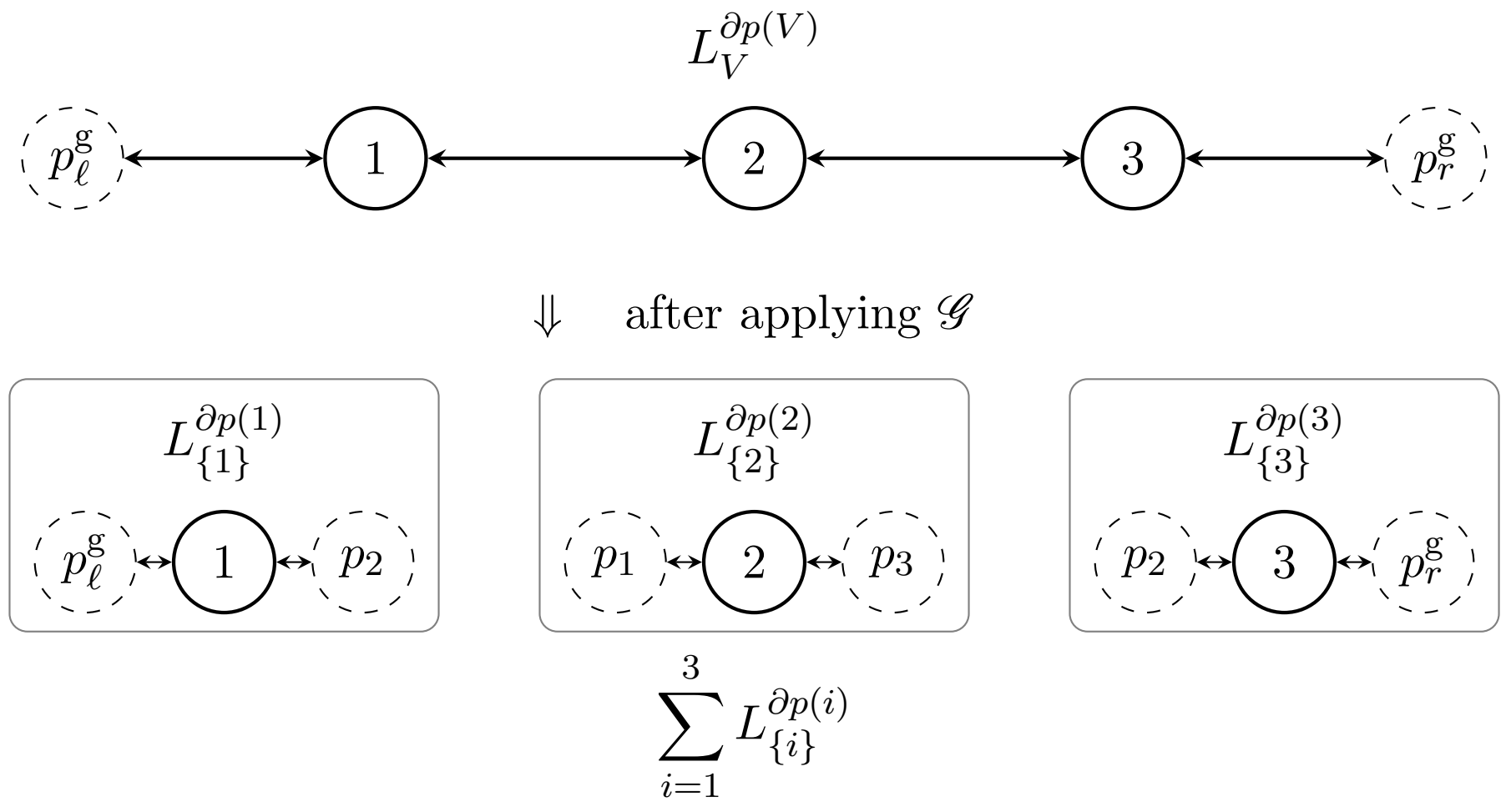

Figure~\ref{fig:single-site-decomposition} illustrates Theorem~\ref{thm:decomposition} on a three-site chain. In the original system, particles move between neighboring physical sites and enter or leave through two fixed ghost reservoirs. In the decomposed representation, each physical site evolves separately, and its former physical neighbors appear as reservoirs carrying the corresponding variable parameters.

\begin{proof}[Theorem \ref{thm:decomposition}]
Fix an edge $\{i,j\}\in E$ and $\gamma_i,\gamma_j\in\N_0$, and define
\begin{equation}
\1_{\gamma_i,\gamma_j}(\eta)
:=
\1_{\{\eta_i=\gamma_i,\,\eta_j=\gamma_j\}}.
\label{eq:two-site-indicator}
\end{equation}
When the two occupations are denoted by $n$ and $m$, respectively, we write
$\1_{\gamma_i,\gamma_j}(n,m)=\1_{\{n=\gamma_i,\,m=\gamma_j\}}$.
We shall use the identity
\begin{equation}
\nu_q(n+k)\frac{\phi(n+k)}{\phi(n)}
=
q^k\nu_q(n),
\qquad n,k\in\N_0,
\label{eq:equilibrium-shift-identity}
\end{equation}
which follows immediately from the definition of $\nu_q$ in
\eqref{eq:equilibrium-marginals}.\\

Our strategy is to compute both $\caG L_{\{i,j\}}\1_{\gamma_i,\gamma_j}$ and $\caG \left(L_i^{p_j}+L_j^{p_i}\right)\1_{\gamma_i,\gamma_j}$ and show that these are equal for arbitrary $\gamma_i, \gamma_j \ \in \N_0$. We then conclude that $\caG L_{\{i,j\}}f$ and $\caG \left(L_i^{p_j}+L_j^{p_i}\right)f$ must be equal for any function $f: \N_0^V \rightarrow \R$, using the fact that $\gamma_i,\gamma_j$ are arbitrary. We first write out the action of $\caG L_{\{i,j\}}$ on $\1_{\gamma_i, \gamma_j}$.

\begin{align}
&\bigl(\caG L_{\{i,j\}}\1_{\gamma_i,\gamma_j}\bigr)(p)
\nonumber\\
&\quad =\sum_{n,m\geq0}\nu_{p_i}(n)\nu_{p_j}(m)\Bigg\{
\sum_{k=1}^{n}M_{i,j}(k)\frac{\phi(n)}{\phi(n-k)}
\bigl[\1_{\gamma_i,\gamma_j}(n-k,m+k)-\1_{\gamma_i,\gamma_j}(n,m)\bigr]
\nonumber\\
&\hspace{40.9mm}+
\sum_{k=1}^{m}M_{j,i}(k)\frac{\phi(m)}{\phi(m-k)}
\bigl[\1_{\gamma_i,\gamma_j}(n+k,m-k)-\1_{\gamma_i,\gamma_j}(n,m)\bigr]\Bigg\}
\nonumber\\
&\quad=\sum_{k\geq1}M_{i,j}(k)
\sum_{\substack{n\geq k\\m\geq0}}\nu_{p_i}(n)\nu_{p_j}(m)
\frac{\phi(n)}{\phi(n-k)}
\bigl[\1_{\gamma_i,\gamma_j}(n-k,m+k)-\1_{\gamma_i,\gamma_j}(n,m)\bigr]
\nonumber\\
&\qquad+\sum_{k\geq1}M_{j,i}(k)
\sum_{\substack{n\geq0\\m\geq k}}\nu_{p_i}(n)\nu_{p_j}(m)
\frac{\phi(m)}{\phi(m-k)}
\bigl[\1_{\gamma_i,\gamma_j}(n+k,m-k)-\1_{\gamma_i,\gamma_j}(n,m)\bigr].\label{eq:edge-generator-on-indicator_1}
\end{align}

For the second equality we interchange the sum over $k$ with the sums over $m,n$. We now evaluate the indicators, which gives

\begin{align}
&\bigl(\caG L_{\{i,j\}}\1_{\gamma_i,\gamma_j}\bigr)(p)
\nonumber\\
&\quad=\sum_{k=1}^{\gamma_j}M_{i,j}(k)
\nu_{p_i}(\gamma_i+k)\nu_{p_j}(\gamma_j-k)
\frac{\phi(\gamma_i+k)}{\phi(\gamma_i)}
\nonumber -\nu_{p_i}(\gamma_i)\nu_{p_j}(\gamma_j)
\sum_{k=1}^{\gamma_i}M_{i,j}(k)
\frac{\phi(\gamma_i)}{\phi(\gamma_i-k)}
\nonumber\\
&\qquad+\sum_{k=1}^{\gamma_i}M_{j,i}(k)
\nu_{p_i}(\gamma_i-k)\nu_{p_j}(\gamma_j+k)
\frac{\phi(\gamma_j+k)}{\phi(\gamma_j)} -\nu_{p_i}(\gamma_i)\nu_{p_j}(\gamma_j)
\sum_{k=1}^{\gamma_j}M_{j,i}(k)
\frac{\phi(\gamma_j)}{\phi(\gamma_j-k)}
\nonumber\\
&\quad=\nu_{p_i}(\gamma_i)
\sum_{k=1}^{\gamma_j}M_{i,j}(k)p_i^k\nu_{p_j}(\gamma_j-k)
+\nu_{p_j}(\gamma_j)
\sum_{k=1}^{\gamma_i}M_{j,i}(k)p_j^k\nu_{p_i}(\gamma_i-k)
\nonumber\\
&\qquad-\nu_{p_i}(\gamma_i)\nu_{p_j}(\gamma_j)\left[
\sum_{k=1}^{\gamma_i}M_{i,j}(k)\frac{\phi(\gamma_i)}{\phi(\gamma_i-k)}
+\sum_{k=1}^{\gamma_j}M_{j,i}(k)\frac{\phi(\gamma_j)}{\phi(\gamma_j-k)}
\right].\label{eq:edge-generator-on-indicator_2}
\end{align}

For the last equality we use \eqref{eq:equilibrium-shift-identity}.

We repeat the computation above for a one-site generator. Writing out the action of $\caG L_i^{p_j}$ on $\1_{\gamma_i,\gamma_j}$ gives
\begin{align}
&\bigl(\caG L_i^{p_j}\1_{\gamma_i,\gamma_j}\bigr)(p)
\nonumber\\
&\quad =\sum_{n,m\geq0}\nu_{p_i}(n)\nu_{p_j}(m)\Bigg\{
\sum_{k\geq1}M_{j,i}(k)p_j^k
\bigl[\1_{\gamma_i,\gamma_j}(n+k,m)-\1_{\gamma_i,\gamma_j}(n,m)\bigr]
\nonumber\\
&\quad \hspace{34mm}+\sum_{k=1}^{n}M_{i,j}(k)\frac{\phi(n)}{\phi(n-k)}
\bigl[\1_{\gamma_i,\gamma_j}(n-k,m)-\1_{\gamma_i,\gamma_j}(n,m)\bigr]
\Bigg\}
\nonumber\\
&\quad=\sum_{k\geq1}M_{j,i}(k)p_j^k
\sum_{n,m\geq0}\nu_{p_i}(n)\nu_{p_j}(m)
\bigl[\1_{\gamma_i,\gamma_j}(n+k,m)-\1_{\gamma_i,\gamma_j}(n,m)\bigr]
\nonumber\\
&\qquad +\sum_{k\geq1}M_{i,j}(k)
\sum_{\substack{n\geq k\\m\geq0}}\nu_{p_i}(n)\nu_{p_j}(m)
\frac{\phi(n)}{\phi(n-k)}
\bigl[\1_{\gamma_i,\gamma_j}(n-k,m)-\1_{\gamma_i,\gamma_j}(n,m)\bigr]
\nonumber
\label{eq:one-site-generator-on-indicator}
\end{align}

For the second equality, we interchanged the sums over $k$ and $n,m$. We proceed by evaluating the indicators and using \eqref{eq:equilibrium-shift-identity} to obtain

\begin{align}
&\bigl(\caG L_i^{p_j}\1_{\gamma_i,\gamma_j}\bigr)(p)
\nonumber\\
&\quad =\nu_{p_j}(\gamma_j)
\sum_{k=1}^{\gamma_i}M_{j,i}(k)p_j^k\nu_{p_i}(\gamma_i-k)
-\nu_{p_i}(\gamma_i)\nu_{p_j}(\gamma_j)
\sum_{k\geq1}M_{j,i}(k)p_j^k
\nonumber\\
&\qquad +\nu_{p_j}(\gamma_j)\sum_{k\geq1}M_{i,j}(k)
\nu_{p_i}(\gamma_i+k)\frac{\phi(\gamma_i+k)}{\phi(\gamma_i)}
-\nu_{p_i}(\gamma_i)\nu_{p_j}(\gamma_j)
\sum_{k=1}^{\gamma_i}M_{i,j}(k)
\frac{\phi(\gamma_i)}{\phi(\gamma_i-k)}
\nonumber\\
&\quad= \nu_{p_j}(\gamma_j)
\sum_{k=1}^{\gamma_i}M_{j,i}(k)p_j^k\nu_{p_i}(\gamma_i-k)
-\nu_{p_i}(\gamma_i)\nu_{p_j}(\gamma_j)
\sum_{k\geq1}M_{j,i}(k)p_j^k
\nonumber\\
&\qquad +\nu_{p_i}(\gamma_i)\nu_{p_j}(\gamma_j)
\sum_{k\geq1}M_{i,j}(k)p_i^k
-\nu_{p_i}(\gamma_i)\nu_{p_j}(\gamma_j)
\sum_{k=1}^{\gamma_i}M_{i,j}(k)
\frac{\phi(\gamma_i)}{\phi(\gamma_i-k)}.
\end{align}

 Interchanging $i$ and
$j$ gives the analogous formula for $L_j^{p_i}$.  When the two formulas are
added, the sums over $k$ from $1$ to $\infty$ cancel pairwise, and therefore
\begin{align}
&\bigl(\caG(L_i^{p_j}+L_j^{p_i})
\1_{\gamma_i,\gamma_j}\bigr)(p)
\nonumber\\
&\quad=\nu_{p_i}(\gamma_i)
\sum_{k=1}^{\gamma_j}M_{i,j}(k)p_i^k\nu_{p_j}(\gamma_j-k)
+\nu_{p_j}(\gamma_j)
\sum_{k=1}^{\gamma_i}M_{j,i}(k)p_j^k\nu_{p_i}(\gamma_i-k)
\nonumber\\
&\qquad-\nu_{p_i}(\gamma_i)\nu_{p_j}(\gamma_j)\left[
\sum_{k=1}^{\gamma_i}M_{i,j}(k)\frac{\phi(\gamma_i)}{\phi(\gamma_i-k)}
+\sum_{k=1}^{\gamma_j}M_{j,i}(k)\frac{\phi(\gamma_j)}{\phi(\gamma_j-k)}
\right]
\nonumber\\
&\quad=\bigl(\caG L_{\{i,j\}}\1_{\gamma_i,\gamma_j}\bigr)(p).
\label{eq:single-edge-decomposition}
\end{align}

Since $\gamma_i$ and $\gamma_j$ are arbitrary, linearity first extends
\eqref{eq:single-edge-decomposition} to every finitely supported function of
the two occupations. For a general $f$ in the stated domain, apply this
identity with all occupations outside $\{i,j\}$ fixed and then average over
them. The assumed finiteness permits the corresponding limiting argument.
Consequently,
\begin{equation}
\bigl(\caG L_{\{i,j\}}f\bigr)(p)
=
\bigl(\caG (L_i^{p_j}+L_j^{p_i})f\bigr)(p).
\label{eq:single-edge-decomposition-general}
\end{equation}

Finally, by \eqref{eq:single-site-generator},
\begin{equation}
\sum_{i\in V}L_{\{i\}}^{\partial p(i)}
=
\sum_{\{i,j\}\in E}(L_i^{p_j}+L_j^{p_i})
+
\sum_{i\in V}\sum_{j\in V_g(i)}L_i^{p_j}.
\label{eq:sum-of-single-site-generators}
\end{equation}
Using \eqref{eq:full-generator} and the linearity of $\caG$, we obtain
\begin{align}
&\bigl(\caG L_V^{\partial p(V)}f\bigr)(p)
\nonumber\\
&\quad=\sum_{\{i,j\}\in E}
\bigl(\caG L_{\{i,j\}}f\bigr)(p)
+\sum_{i\in V}\sum_{j\in V_g(i)}
\bigl(\caG L_i^{p_j}f\bigr)(p)
\nonumber\\
&\quad=\sum_{\{i,j\}\in E}
\bigl(\caG(L_i^{p_j}+L_j^{p_i})f\bigr)(p)
+\sum_{i\in V}\sum_{j\in V_g(i)}
\bigl(\caG L_i^{p_j}f\bigr)(p)
\nonumber\\
&\quad=\left(\caG\left[
\sum_{\{i,j\}\in E}(L_i^{p_j}+L_j^{p_i})
+\sum_{i\in V}\sum_{j\in V_g(i)}L_i^{p_j}
\right]f\right)(p)
\nonumber\\
&\quad=\left(\caG\left[
\sum_{i\in V}L_{\{i\}}^{\partial p(i)}
\right]f\right)(p).
\label{eq:decomposition-final-step}
\end{align}
For the second equality we use
\eqref{eq:single-edge-decomposition-general} on each physical edge and for the last equality
\eqref{eq:sum-of-single-site-generators}.  We recovered
\eqref{eq:decomposition} and hence proved Theorem \ref{thm:decomposition}.
\end{proof}

\subsection{Stationary measures and Markovian mixture structure} \label{subsec:stat_meas_Markov_mixture_structure}

In this paper we aim to find the stationary states of boundary-driven processes that correspond to a generator of the form $L^{\partial p (V)}_V$ by integrating \eqref{eq:decomposition} against a probability measure $\Lambda$ on $\mathsf D^V$, i.e.

\begin{equation}
\int \left(
\caG\left[L_V^{\partial p(V)}f\right]
\right)(p) \, \Lambda(dp)
=
\int \left(
\caG\left[
\sum_{i\in V}L_{\{i\}}^{\partial p(i)}f
\right]
\right)(p) \, \Lambda(dp).
\label{eq:integrated_single_site_decomposition_fst}
\end{equation}

In particular, if the left-hand side of \eqref{eq:integrated_single_site_decomposition_fst} is zero, then

\be 
\mu_\Lambda (\eta) = \int_{\mathsf D^V} \left(\otimes_{i \in V} \nu_{p_i}\right) (\eta) \Lambda(dp), \qquad \eta \in \N_0^V,
\ee

is stationary for the boundary-driven process $L^{\partial p (V)}_V$. The decomposition into single-site generators on the right-hand side of \eqref{eq:integrated_single_site_decomposition_fst} provides a connection between the stationary states of the single-site generators and the stationary state of $L^{\partial p (V)}_V$. Here we demonstrate how this decomposition can be used, in particular cases, to obtain the NESS of $L^{\partial p (V)}_V$ effectively from the NESS of the single-site generators $L^{\partial p (i)}_{\{i\}}$, $i \in V$.\\

\subsubsection{Stationary Markov extensions} \label{subsubsec:Stationary_Markov_extensions}

Let $\Lambda$ be a Borel probability measure on $\mathsf D^V$. On the canonical probability space
\begin{equation}
\bigl(\mathsf D^V,\mathcal B(\mathsf D^V),\Lambda\bigr),
\label{eq:canonical-parameter-space}
\end{equation}

define the random vector

\begin{equation}
P=(P_i)_{i\in V},
\label{eq:canonical-parameter-vector}
\end{equation}

Thus, $\Lambda$ is the law of $P$, while $p=(p_i)_{i\in V}$ denotes a deterministic realization of this vector. The ghost parameters $p_j^{\mathrm g}$, $j\in V_g$, remain fixed and are not coordinates of $P$.

\begin{definition}[Mixture of product measures]
\label{def:product-mixture}
The probability measure on $\N_0^V$ obtained by mixing with respect to $\Lambda$ is
\begin{equation}
\mu_\Lambda(\eta)
:=
\int_{\mathsf D^V}
\left(\bigotimes_{i\in V}\nu_{p_i}\right)(\eta)
\,\Lambda(dp),
\qquad
\eta\in\N_0^V.
\label{eq:mixture_measure_on_V}
\end{equation}

We refer to $\mu_\Lambda$ as a mixture of product measures on $V$ associated with the mixing measure $\Lambda$.
\end{definition}

For every test function $f$, Theorem~\ref{thm:decomposition} gives
\begin{align}
\sum_{\eta\in\N_0^V}\mu_\Lambda(\eta)
\bigl(L_V^{\partial p(V)}f\bigr)(\eta)
&=
\int_{\mathsf D^V}
\bigl(\caG L_V^{\partial p(V)}f\bigr)(p)\,\Lambda(dp)
\nonumber\\
&=
\sum_{i\in V}
\int_{\mathsf D^V}
\bigl(\caG L_{\{i\}}^{\partial p(i)}f\bigr)(p)\,\Lambda(dp).
\label{eq:integrated_single_site_decomposition}
\end{align}
To prove the stationarity of $\mu_\Lambda$, it suffices to choose $\Lambda$ so that every integral in the sum on the second line of \eqref{eq:integrated_single_site_decomposition} vanishes. This happens when the distribution of $P_i$, conditioned on $P_j$ for all $j \neq i$, depends only on $\partial P (i)$ and corresponds to the stationary measure of the single-site generator $L^{\partial p (i)}_{\{i\}}$. In other words, $\Lambda$ satisfies a spatial Markov property with the correct marginals. To make this precise, we need to condition $\Lambda$ on all components of $P$ except the $i$-th. Write

\begin{equation}
P_{\neq i}:=(P_j)_{j\in V\setminus\{i\}},
\qquad
p_{\neq i}:=(p_j)_{j\in V\setminus\{i\}},
\label{eq:parameter-profile-without-i}
\end{equation}
and let $\Lambda_{\neq i}$ be the law of $P_{\neq i}$. Since $\mathsf D^V$ is a standard Borel space, there exists a regular conditional probability kernel
\begin{equation}
K_i(p_{\neq i},A)
:=
\Lambda\bigl(P_i\in A\mid P_{\neq i}=p_{\neq i}\bigr),
\qquad A\in\mathcal B(\mathsf D),
\label{eq:regular-conditional-kernel}
\end{equation}
defined for $\Lambda_{\neq i}$-almost every $p_{\neq i}$. It is characterized by
\begin{equation}
\Lambda\bigl(P_i\in A,\ P_{\neq i}\in B\bigr)
=
\int_B K_i(p_{\neq i},A)\,\Lambda_{\neq i}(dp_{\neq i})
\label{eq:conditional-kernel-characterization}
\end{equation}
for all Borel sets $A\subset\mathsf D$ and $B\subset\mathsf D^{V\setminus\{i\}}$. In differential notation, this disintegration is written as
\begin{equation}
\Lambda(dp)
=
\Lambda_{\neq i}(dp_{\neq i})K_i(p_{\neq i},dp_i).
\label{eq:lambda-disintegration}
\end{equation}

\begin{definition}[Markov extension]
\label{def:markov-extension}
Suppose that we are given a family of Borel probability measures $\lambda^{\partial p (i)}_i$, $p \in \mathsf D^V$ and $i \in V$, on $\mathsf D$. A probability measure $\Lambda$ on $\mathsf D^V$ is a Markov extension of this family of single-site probability measures if, for every $i\in V$,
\begin{equation}
K_i(p_{\neq i},dq)
=
\lambda_i^{\partial p(i)}(dq)
\quad\text{for }\Lambda_{\neq i}\text{-almost every }p_{\neq i}.
\label{eq:spatial_markov}
\end{equation}
\end{definition}

\br
We remark that the Markov extension does not necessarily exist for a given family $\lambda^{\partial p (i)}_i$, $p \in \mathsf D^V$. In Section \ref{subsubsec:Markov_extensions_on_the_chain} we provide examples of families for which the Markov extension exists.
\er

Because the graph is simple, $i\notin\partial\{i\}$, and therefore $\partial p(i)$ is determined by $p_{\neq i}$ together with the fixed ghost parameters. Definition~\ref{def:markov-extension} requires that the conditional law of $P_i$ given the entire complementary random parameter profile depends on that profile only through the random parameters at the neighbors of $i$.\\

To extend the stationary measures of the single-site operators $L^{\partial p (i)}_{\{i\}}$ to the stationary measure of $L^{\partial p(V)}_V$, one first needs a mixture structure for the single-site stationary measures as in Definition \ref{def:single-site-mixing-kernels}. If $\Lambda$ is then the Markov extension of the single-site mixing measures, the stationarity of $\mu_\Lambda$ in \eqref{eq:mixture_measure_on_V} follows. This is the content of Theorem \ref{thm:markov-mixture-stationary}.

For $i\in V$, let
\begin{equation}
\mathsf P_i
:=
\bigl\{
r\in\mathsf D^{\partial\{i\}}:
r_j=p_j^{\mathrm g}\ \text{for every }j\in\partial\{i\}\cap V_g
\bigr\}.
\label{eq:local-parameter-space}
\end{equation}
denote the space of parameter families on the neighbors of $i$, with the coordinates belonging to ghost sites fixed at their prescribed values.

\begin{definition}[Single-site stationary mixing kernels]
\label{def:single-site-mixing-kernels}
For every $i\in V$ and every neighboring profile $r\in\mathsf P_i$, let $\lambda_i^r$ be a probability measure on $\mathsf D$, and define
\begin{equation}
\mu_i^r(n)
:=
\int_{\mathsf D}\nu_q(n)\,\lambda_i^r(dq),
\qquad n\in\N_0.
\label{eq:stat_mix_one_site}
\end{equation}
We call $(\lambda_i^r)_{r\in\mathsf P_i}$ a family of single-site stationary mixing kernels if $\mu_i^r$ is stationary for $L_{\{i\}}^r$ for every admissible $r$.
\end{definition}

\begin{theorem}[Stationarity of a Markovian mixture]
\label{thm:markov-mixture-stationary}
Suppose that $(\lambda_i^r)_{r\in\mathsf P_i}$ is a family of single-site stationary mixing kernels for every $i\in V$, and that $\Lambda$ is a Markov extension of these kernels. Then the mixture of product measures $\mu_\Lambda$ in \eqref{eq:mixture_measure_on_V} is stationary for $L_V^{\partial p(V)}$.
\end{theorem}

\begin{proof}
For $\eta\in\N_0^V$, write $\eta_{\neq i}:=(\eta_j)_{j\in V\setminus\{i\}}$ and use $(\eta_{\neq i},n)$ for the configuration obtained by assigning occupation $n$ to site $i$. Fix $i\in V$ and disintegrate $\Lambda$ with respect to $P_{\neq i}$. Since $L_{\{i\}}^{\partial p(i)}$ changes only the occupation at $i$ and its rates do not depend on $p_i$, Definition~\ref{def:markov-extension} yields
\begin{align}
&\int_{\mathsf D^V}
\bigl(\caG L_{\{i\}}^{\partial p(i)}f\bigr)(p)\,\Lambda(dp)
\nonumber\\
&\quad=
\int\Lambda_{\neq i}(dp_{\neq i})
\sum_{\eta_{\neq i}\in\N_0^{V\setminus\{i\}}}
\left(\prod_{j\ne i}\nu_{p_j}(\eta_j)\right)
\sum_{n\in\N_0}
\left[\int_{\mathsf D}\nu_q(n)\lambda_i^{\partial p(i)}(dq)\right]
\bigl(L_{\{i\}}^{\partial p(i)}f\bigr)(\eta_{\neq i},n)
\nonumber\\
&\quad=0.
\label{eq:single-site-mixture-vanishes}
\end{align}
The final equality follows because the expression in square brackets is $\mu_i^{\partial p(i)}(n)$, which is stationary for the corresponding one-site generator. Summing over $i$ in \eqref{eq:integrated_single_site_decomposition} proves the claim.
\end{proof}

\br \label{rmk:incomplete_single_site_family}
We remark that for specific boundary parameters $\partial p (V)$, one does not necessarily need the whole family of single-site stationary mixing kernels to identify $\Lambda$ as a Markov extension. The reason is that Definition \ref{def:markov-extension} only requires \eqref{eq:spatial_markov} for $\Lambda_{\neq i}$-almost every $p_{\neq i}$. In all the situations we will encounter, $\Lambda_{\neq i}$ does not have full support. Proposition \ref{prop:homogeneous-stationary-state} contains a first example of such a $\Lambda$, whose support is restricted to constant parameter profiles.
\er

An elementary special case is obtained when the mixing law is concentrated on a
constant profile. Consider the measures 

\be
\lambda^{\partial p (i)}_i (p_i) = \delta \left(\forall j \sim i: \left(\partial p (i)\right)_j = p_i = q \right), \qquad i \in V.
\ee

Then,

\be
\Lambda_q (p) := \delta(\forall i \in V: p_i = q )
\ee

satisfies \eqref{eq:spatial_markov}. The mixing measures $\Lambda_q$ and $\lambda^{\partial p (i)}_i$, $i \in V$, correspond to the stationary measures under equilibrium boundary driving for models where, for each $k \in \N$ and $i \in V$, the cumulative mass function directed into site $i$ is equal to the cumulative mass function directed out of site $i$. This is the content of Proposition \ref{prop:homogeneous-stationary-state}. The aforementioned restriction on the mass functions should be understood as in \eqref{eq:cummulative_symmetry_mass_functions}.

\begin{proposition}[Stationary state under equilibrium boundary driving]
\label{prop:homogeneous-stationary-state}
Fix $q\in\mathsf D$. Suppose that
\begin{equation}
p_j^{\mathrm g}=q,
\qquad j\in V_g,
\label{eq:homogeneous-ghost-parameters}
\end{equation}
and that, for every $i\in V$ and $k\in\N$,
\begin{equation}
\sum_{j\sim i}M_{i,j}(k)
=
\sum_{j\sim i}M_{j,i}(k).
\label{eq:cummulative_symmetry_mass_functions}
\end{equation}
Then
\begin{equation}
\mu_q:=\bigotimes_{i\in V}\nu_q
\label{eq:homogeneous-product-stationary-measure}
\end{equation}
is stationary for $L_V^{\partial p(V)}$.
\end{proposition}

\begin{proof}
Let $\mathbf q:=(q)_{i\in V}$ and $\Lambda_q:=\delta_{\mathbf q}$. Its conditional kernels are $K_i(p_{\neq i},dq')=\delta_q(dq')$ for $\Lambda_{q,\neq i}$-almost every $p_{\neq i}$. For fixed $i$, the cumulative-rate assumption \eqref{eq:cummulative_symmetry_mass_functions} and \eqref{eq:equilibrium-marginals} imply, for every $n\in\N_0$ and $k\geq1$,
\begin{align}
\nu_q(n)q^k\sum_{j\sim i}M_{j,i}(k)
&=
\nu_q(n)q^k\sum_{j\sim i}M_{i,j}(k)
\nonumber\\
&=
\nu_q(n+k)\frac{\phi(n+k)}{\phi(n)}
\sum_{j\sim i}M_{i,j}(k).
\label{eq:single-site-detailed-balance}
\end{align}
Thus $\nu_q$ satisfies detailed balance for $L_{\{i\}}^{\partial p(i)}$ at $p=\mathbf q$. Hence $\Lambda_q$ is a Markov extension of the stationary one-site mixing laws $\delta_q$, and Theorem ~\ref{thm:markov-mixture-stationary} gives the result.
\end{proof}

\subsubsection{Markov extensions on the chain} \label{subsubsec:Markov_extensions_on_the_chain}

Theorem ~\ref{thm:markov-mixture-stationary} shows that the stationary measure of the full boundary-driven process can be obtained from stationary mixtures for the single-site generators, provided that their mixing kernels admit a common Markov extension. Finding the stationary state of the whole system from single-site mixing kernels is therefore a compatibility problem. Given the families of single-site stationary mixing kernels, one must construct a probability measure $\Lambda$ on $\mathsf D^V$ whose regular conditional kernels agree with them in the sense of \eqref{eq:spatial_markov}. Such a measure need not exist for an arbitrary collection of single-site kernels, and even when it exists its construction is not immediate.

We identify two classes of compatible kernels for which the Markov extension can be written explicitly. Both constructions concern a chain and use monotone parameter profiles: the first has a two-point parameter space, while the second is its continuous analogue. Let
\begin{equation}
V=\{1,\ldots,N\},
\qquad
E=\bigl\{\{i,i+1\}:1\leq i<N\bigr\},
\label{eq:chain-geometry}
\end{equation}
The left and right ghost sites are labeled $0$ and $N+1$, respectively. To simplify the notation, we denote their fixed parameters by $z_\ell$ and $z_r$, i.e.

\be
p_0^{\mathrm g}:=z_\ell, \qquad
p_{N+1}^{\mathrm g}:=z_r.
\ee

Thus, $z_\ell$ and $z_r$ are precisely the fixed parameters of the left and right boundary reservoirs. In some expressions it is convenient to write $p_0$ and $p_{N+1}$ instead. Hence we identify, as we did before,

\be
p_0 := p_0^{\mathrm g}=z_\ell, \qquad p_{N+1} := p_{N+1}^{\mathrm g}=z_r
\ee

\begin{figure}[ht]
\centering
\begin{tikzpicture}[
    site/.style={circle,draw,thick,minimum size=10mm,inner sep=0pt},
    reservoir/.style={circle,draw,dashed,minimum size=10mm,inner sep=0pt},
    every node/.style={font=\small},
    >=stealth
]
\node[reservoir] (leftghost) at (0,0) {$0$};
\node[site] (siteone) at (2.2,0) {$1$};
\node[site] (sitetwo) at (4.4,0) {$2$};
\node[site] (siten) at (8,0) {$N$};
\node[reservoir] (rightghost) at (10.2,0) {$N+1$};

\draw[<->,thick] (leftghost) -- (siteone);
\draw[<->,thick] (siteone) -- (sitetwo);
\draw[<->,thick] (sitetwo) -- (siten);
\draw[<->,thick] (siten) -- (rightghost);
\node[fill=white,inner sep=3pt] at (6.2,0) {$\cdots$};

\node[font=\scriptsize,align=center] at (0,0.85) {left ghost\\reservoir};
\node[font=\scriptsize,align=center] at (10.2,0.85) {right ghost\\reservoir};
\node at (0,-0.75) {$p_0^{\mathrm g}=z_\ell$};
\node at (2.2,-0.75) {$p_1$};
\node at (4.4,-0.75) {$p_2$};
\node at (8,-0.75) {$p_N$};
\node at (10.2,-0.75) {$p_{N+1}^{\mathrm g}=z_r$};
\node at (5.1,-1.55)
{$z_\ell=p_0\leq p_1\leq p_2\leq\cdots\leq p_N\leq p_{N+1}=z_r$};
\end{tikzpicture}
\caption{Chain geometry and boundary reservoirs used in Propositions~\ref{prop:markov_extensions_chain_discrete} and~\ref{prop:markov_extensions_chain_continuous}. Solid circles represent the physical sites and dashed circles the fixed ghost reservoirs.}
\label{fig:markov-extension-chain}
\end{figure}

In both constructions, the global mixing measure is supported on nondecreasing profiles. Consequently, for every profile in its support,
\begin{equation}
p_{i}\leq p_{i+1},
\qquad i=0,1,\ldots,N.
\label{eq:ordered-neighbor-parameters}
\end{equation}

We first consider a two-point parameter space. That is, we fix $x,y\in\mathsf D$ with $x<y$ and positive weights $g_x,g_y>0$. For $a,b\in\{x,y\}$ with $a\leq b$, define the probability measure
\begin{equation}
\lambda_{\mathrm d}^{a,b}
:=
\begin{cases}
\delta_x, & a=b=x,\\[1mm]
\dfrac{g_x\delta_x+g_y\delta_y}{g_x+g_y},
& a=x,\ b=y,\\[3mm]
\delta_y, & a=b=y.
\end{cases}
\label{eq:discrete-single-site-kernels}
\end{equation}
Thus, $\lambda_{\mathrm d}^{a,b}$ is supported on the values in $\{x,y\}$ lying between the neighboring parameters $a$ and $b$, with relative weights $g_x$ and $g_y$ whenever both values are admissible. The subscript $d$ in $\lambda_{\mathrm d}^{a,b}$ indicates that these measures have support in a discrete set (i.e., $\{x,y\}$).

\begin{proposition}[Discrete Markov extension]
\label{prop:markov_extensions_chain_discrete}
For $z_\ell,z_r\in\{x,y\}$ with $z_\ell\leq z_r$, define a probability measure $\Lambda_{z_\ell,z_r}^{(N)}$ on $\{x,y\}^N$ by

\be
\Lambda_{z_\ell,z_r}^{(N)}(\{p\}):=
\frac{1}{Z_{z_\ell,z_r}^{(N)}}
\left(\prod_{r=1}^N g_{p_r}\right)
\left(\prod_{r=0}^N\1\{p_r\leq p_{r+1}\}\right)
\label{eq:discrete-chain-extension}
\ee

where

\begin{equation}
Z_{z_\ell,z_r}^{(N)}
:=
\sum_{p\in\{x,y\}^N}
\left(\prod_{r=1}^N g_{p_r}\right)
\1\{z_\ell\leq p_1\leq\cdots\leq p_N\leq z_r\}.
\label{eq:discrete-chain-normalization}
\end{equation}
Let $P=(P_1,\ldots,P_N)$ have law $\Lambda_{z_\ell,z_r}^{(N)}$, and denote the law of $P_{\neq i}$ by $\Lambda_{z_\ell,z_r,\neq i}^{(N)}$. Then, for every $i\in V$, the probability kernel associated with $P_i$ conditional on $P_{\neq i}$ is
\begin{equation}
K_i(p_{\neq i},dq)
=
\lambda_{\mathrm d}^{p_{i-1},p_{i+1}}(dq)
\quad
\text{for }\Lambda_{z_\ell,z_r,\neq i}^{(N)}\text{-almost every }p_{\neq i}.
\label{eq:discrete-chain-conditional-kernel}
\end{equation}
In particular, $\Lambda_{z_\ell,z_r}^{(N)}$ is a Markov extension, in the sense of Definition~\ref{def:markov-extension}, of the single-site kernels
\begin{equation}
\lambda_i^{\partial p(i)}
:=
\lambda_{\mathrm d}^{p_{i-1},p_{i+1}},
\qquad i=1,\ldots,N.
\label{eq:discrete-local-kernel-identification}
\end{equation}
\end{proposition}

\begin{proof}
Fix $i\in V$ and a value $p_{\neq i}$ having positive $\Lambda_{z_\ell,z_r,\neq i}^{(N)}$-probability. That is,

\be
z_\ell \leq p_1 \leq \ldots \leq p_{i-1}  \leq p_{i+1} \leq \ldots \leq p_{N} \leq z_r.
\ee

From \eqref{eq:discrete-chain-extension}, we can see that the dependence on the remaining coordinate $p_i$ is entirely contained in $g_{p_i}\1\{p_{i-1}\leq p_i\leq p_{i+1}\}$. Normalizing this factor over the two possible values of $p_i$ gives, for $q\in\{x,y\}$,
\begin{equation}
\Lambda_{z_\ell,z_r}^{(N)}
\bigl(P_i=q\mid P_{\neq i}=p_{\neq i}\bigr)
=
\frac{
g_q\1\{p_{i-1}\leq q\leq p_{i+1}\}
}{
\displaystyle
\sum_{u\in\{x,y\}}
g_u\1\{p_{i-1}\leq u\leq p_{i+1}\}
}.
\label{eq:discrete-conditional-ratio}
\end{equation}

There are three possible pairs of neighboring values. If $(p_{i-1},p_{i+1})=(x,x)$, only $q=x$ is admissible, so the conditional law is $\delta_x$. If $(p_{i-1},p_{i+1})=(y,y)$, only $q=y$ is admissible, giving $\delta_y$. Finally, if $(p_{i-1},p_{i+1})=(x,y)$, both values are admissible and \eqref{eq:discrete-conditional-ratio} assigns them probabilities $g_x/(g_x+g_y)$ and $g_y/(g_x+g_y)$, respectively. These three conditional laws are precisely those in \eqref{eq:discrete-single-site-kernels}. This proves \eqref{eq:discrete-chain-extension}, and therefore the Markov-extension property.
\end{proof}

We next give the continuous analogue. Let $g:\mathsf D^2\to[0,\infty)$ be continuous and assume that all normalizing constants below are strictly positive. For $a,b\in\mathsf D$ with $a<b$, define
\begin{equation}
\lambda_{\mathrm c}^{a,b}(dq)
:=
\frac{
g(a,q)g(q,b)\1\{a\leq q\leq b\}
}{
\displaystyle\int_a^b g(a,t)g(t,b)\,dt
}\,dq,
\qquad
\lambda_{\mathrm c}^{a,a}:=\delta_a.
\label{eq:continuous-single-site-kernels}
\end{equation}

Here, the subscript $\mathrm c$ indicates that these measures have support in the continuous parameter domain $\mathsf D$.

\begin{proposition}[Continuous Markov extension]
\label{prop:markov_extensions_chain_continuous}
For $z_\ell,z_r\in\mathsf D$ with $z_\ell<z_r$, let $\Lambda_{z_\ell,z_r}^{(N)}$ be the probability measure on $\mathsf D^N$ with density
\begin{equation}
\frac{d\Lambda_{z_\ell,z_r}^{(N)}}{dp}(p)
=
\frac{1}{Z_{z_\ell,z_r}^{(N)}}
\left(\prod_{r=0}^{N}g(p_r,p_{r+1})\right)
\1\{z_\ell\leq p_1\leq\cdots\leq p_N\leq z_r\},
\label{eq:continuous-chain-extension}
\end{equation}
where $p_0=z_\ell$, $p_{N+1}=z_r$, and
\begin{equation}
Z_{z_\ell,z_r}^{(N)}
:=
\int_{\mathsf D^N}
\left(\prod_{r=0}^{N}g(p_r,p_{r+1})\right)
\1\{z_\ell\leq p_1\leq\cdots\leq p_N\leq z_r\}\,dp.
\label{eq:continuous-chain-normalization}
\end{equation}
If $P$ has law $\Lambda_{z_\ell,z_r}^{(N)}$ and $\Lambda_{z_\ell,z_r,\neq i}^{(N)}$ denotes the law of $P_{\neq i}$, then, for every $i\in V$,
\begin{equation}
K_i(p_{\neq i},dq)
=
\lambda_{\mathrm c}^{p_{i-1},p_{i+1}}(dq)
\quad
\text{for }\Lambda_{z_\ell,z_r,\neq i}^{(N)}\text{-almost every }p_{\neq i}.
\label{eq:continuous-chain-conditional-kernel}
\end{equation}
Consequently, $\Lambda_{z_\ell,z_r}^{(N)}$ is a Markov extension of the kernels $\lambda_i^{\partial p(i)}:=\lambda_{\mathrm c}^{p_{i-1},p_{i+1}}$. When $z_\ell=z_r$, the corresponding degenerate extension is the point mass on the constant profile $(z_\ell,\ldots,z_\ell)$.
\end{proposition}

\begin{proof}
For fixed $p_{\neq i}$, the only factors in \eqref{eq:continuous-chain-extension} that depend on $p_i$ are
\begin{equation}
g(p_{i-1},p_i)g(p_i,p_{i+1})
\1\{p_{i-1}\leq p_i\leq p_{i+1}\}.
\label{eq:continuous-conditional-factor}
\end{equation}
Normalizing this expression as a function of $p_i$ gives \eqref{eq:continuous-single-site-kernels}. This is also the construction in \cite[Theorem~4.5]{redig2025non}.
\end{proof}

\subsubsection{Single-site stationarity}

In Sections \ref{sec:Harmonic_models} and \ref{sec:Asymmetric models with discrete mixing structure}, we encounter several examples of processes for which the single-site stationary measure has a mixture structure with mixing measures of the form in \eqref{eq:discrete-single-site-kernels} or \eqref{eq:continuous-single-site-kernels}. To obtain these single-site mixtures, we solve the stationarity equations for the process on one site coupled to boundary reservoirs. These stationarity equations are stated in Proposition \ref{prop:single-site_stationarity_equations}. 

\begin{proposition}[Single-site stationarity equations]
\label{prop:single-site_stationarity_equations}
Fix $i\in V$ and an admissible neighboring parameter profile $\partial p(i)$. Let $\lambda_i^{\partial p(i)}$ be a probability measure on $\mathsf D$ and set
\begin{equation}
\mu_i^{\partial p(i)}(n)
:=
\int_{\mathsf D}\nu_q(n)\,\lambda_i^{\partial p(i)}(dq),
\qquad n\in\N_0.
\label{eq:single-site-mixture-for-stationarity}
\end{equation}
Then $\mu_i^{\partial p(i)}$ is stationary for $L_{\{i\}}^{\partial p(i)}$ if and only if, for every $\gamma\in\N_0$,
\begin{equation}
\begin{aligned}
0
={}&
\int_{\mathsf D}\frac{1}{Z(q)}
\Bigg\{
\frac{1}{\phi(\gamma)}
\sum_{k=1}^{\infty}
\left[
\left(\sum_{j\sim i}M_{i,j}(k)\right)q^{\gamma+k}
-
\left(\sum_{j\sim i}M_{j,i}(k)p_j^k\right)q^\gamma
\right]
\\
&\qquad\qquad+
\sum_{k=1}^{\gamma}\frac{1}{\phi(\gamma-k)}
\left[
\left(\sum_{j\sim i}M_{j,i}(k)p_j^k\right)q^{\gamma-k}
-
\left(\sum_{j\sim i}M_{i,j}(k)\right)q^\gamma
\right]
\Bigg\}
\,\lambda_i^{\partial p(i)}(dq).
\end{aligned}
\label{eq:single-site-stationarity-equations}
\end{equation}
\end{proposition}

\begin{proof}
For $\gamma\in\N_0$, let $\1_\gamma(n):=\1\{n=\gamma\}$. We first apply
the definition of the single-site generator:
\begin{equation}
\begin{aligned}
&\sum_{n\in\N_0}\nu_q(n)
\bigl(L_{\{i\}}^{\partial p(i)}\1_\gamma\bigr)(n)
\\
&\quad=
\sum_{n=0}^{\infty}\nu_q(n)
\Bigg\{
\sum_{j\sim i}\sum_{k=1}^{\infty}
M_{j,i}(k)p_j^k
\bigl[\1_\gamma(n+k)-\1_\gamma(n)\bigr]
\\
&\hspace{36mm}+
\sum_{j\sim i}\sum_{k=1}^{n}
M_{i,j}(k)\frac{\phi(n)}{\phi(n-k)}
\bigl[\1_\gamma(n-k)-\1_\gamma(n)\bigr]
\Bigg\}\nonumber\\
&\quad=
\sum_{k=1}^{\infty}
\left(\sum_{j\sim i}M_{j,i}(k)p_j^k\right)
\sum_{n=0}^{\infty}\nu_q(n)
\bigl[\1_\gamma(n+k)-\1_\gamma(n)\bigr]
\\
&\qquad+
\sum_{k=1}^{\infty}
\left(\sum_{j\sim i}M_{i,j}(k)\right)
\sum_{n=k}^{\infty}\nu_q(n)
\frac{\phi(n)}{\phi(n-k)}
\bigl[\1_\gamma(n-k)-\1_\gamma(n)\bigr].
\end{aligned}
\end{equation}
For the second equation we interchanged the sums over $n$ and $k$. We now evaluate the indicators. 

\begin{equation}
\begin{aligned}
&\sum_{n\in\N_0}\nu_q(n)
\bigl(L_{\{i\}}^{\partial p(i)}\1_\gamma\bigr)(n)
\\
&\quad=
\sum_{k=1}^{\gamma}
\left(\sum_{j\sim i}M_{j,i}(k)p_j^k\right)
\nu_q(\gamma-k)
-
\nu_q(\gamma)
\sum_{k=1}^{\infty}
\left(\sum_{j\sim i}M_{j,i}(k)p_j^k\right)
\\
&\qquad+
\sum_{k=1}^{\infty}
\left(\sum_{j\sim i}M_{i,j}(k)\right)
\nu_q(\gamma+k)
\frac{\phi(\gamma+k)}{\phi(\gamma)}
-
\nu_q(\gamma)
\sum_{k=1}^{\gamma}
\left(\sum_{j\sim i}M_{i,j}(k)\right)
\frac{\phi(\gamma)}{\phi(\gamma-k)}.
\end{aligned}
\label{eq:single-site-indicator-calculation}
\end{equation}

Recall that

\begin{equation}
\nu_q(m)=\frac{1}{Z(q)}\frac{q^m}{\phi(m)}.
\label{eq:equilibrium-marginal-explicit-proof}
\end{equation}

Substituting \eqref{eq:equilibrium-marginal-explicit-proof} into
\eqref{eq:single-site-indicator-calculation} and grouping the two sums
over all $k\geq1$ and the two sums over $1\leq k\leq\gamma$ gives
\begin{equation}
\begin{aligned}
&\sum_{n\in\N_0}\nu_q(n)
\bigl(L_{\{i\}}^{\partial p(i)}\1_\gamma\bigr)(n)
\\
&\quad=
\frac{1}{Z(q)}
\Bigg\{
\frac{1}{\phi(\gamma)}
\sum_{k=1}^{\infty}
\left[
\left(\sum_{j\sim i}M_{i,j}(k)\right)q^{\gamma+k}
-
\left(\sum_{j\sim i}M_{j,i}(k)p_j^k\right)q^\gamma
\right]
\\
&\hspace{25mm}+
\sum_{k=1}^{\gamma}\frac{1}{\phi(\gamma-k)}
\left[
\left(\sum_{j\sim i}M_{j,i}(k)p_j^k\right)q^{\gamma-k}
-
\left(\sum_{j\sim i}M_{i,j}(k)\right)q^\gamma
\right]
\Bigg\}.
\end{aligned}
\label{eq:single-site-indicator-calculation-explicit}
\end{equation}

Integrating the left-hand side of
\eqref{eq:single-site-indicator-calculation-explicit} with respect to
$\lambda_i^{\partial p(i)}$ gives
\begin{equation}
\sum_{n\in\N_0}\mu_i^{\partial p(i)}(n)
\bigl(L_{\{i\}}^{\partial p(i)}\1_\gamma\bigr)(n).
\label{eq:single-site-mixture-indicator}
\end{equation}
The probability measure $\mu_i^{\partial p(i)}$ is stationary precisely when \eqref{eq:single-site-mixture-indicator} vanishes for every indicator $\1_\gamma$, $\gamma\in\N_0$. This proves the equivalence with \eqref{eq:single-site-stationarity-equations}.

\end{proof}

\section{Harmonic models} \label{sec:Harmonic_models}

In this section we study the (discrete) harmonic models. These are interacting particle systems that fall within the class of generalized zero-range processes introduced in Section \ref{subsec:def_generalized_zero_range}. We define them in Definition \ref{def:Harmonic_models} by specifying the mass and potential functions.

\bd[Harmonic models] \label{def:Harmonic_models}
The harmonic process with parameter $s>0$ is the generalized zero-range process, in the sense of Definition \ref{def:boundary-driven-generator}, with mass functions $M_{i,j}$, $i,j \in \bar{V}$, and potential function $\phi$ given by

\be \label{eq:mass_and_potential_functions_harmonic}
M_{i,j}(k) = a_{i,j} \frac{1}{k}, \, k\in \N \qquad \text{and} \qquad \phi(n) = \frac{\Gamma(n+1)}{\Gamma(n+2s)}, \, n \in \N_0.
\ee

Here, the coefficient $a_{i,j}\in[0,\infty)$ is called the oriented edge weight from $i$ to $j$. The oriented edge weights are called symmetric whenever $a_{i,j}=a_{j,i}$.
\ed

From \eqref{eq:mass_and_potential_functions_harmonic}, we see that the corresponding equilibrium single-site measures in Definition~\ref{def:equilibrium-single-site-measures} are discrete Gamma, i.e.,

\be
Z(q):=\sum_{m\in\N_0} q^m \frac{\Gamma(m+2s)}{\Gamma(m+1)},
\qquad
\nu_q(n):=\frac{1}{Z(q)}q^n \frac{\Gamma(n+2s)}{\Gamma(n+1)},
\quad n\in\N_0.
\ee

The admissible parameter domain is \(\mathsf D=[0,1)\).

For harmonic models it is most natural to parameterize the equilibrium single-site measures such that their mean is given by

\be
\sum_{n = 0}^\infty n \cdot \nu^{(s)}_\vartheta (n) = 2s \, \vartheta, \qquad \vartheta \geq 0.
\ee

Here, the $\vartheta$-parametrized equilibrium single-site measure $\nu^{(s)}_\vartheta$ is as defined in Definition \ref{def:theta_equilbrium_marginal}. The relation between $q$ and $\vartheta$ is

\be
q = \frac{\vartheta}{1+\vartheta}, \qquad \vartheta = \frac{q}{1-q}.
\ee

For $p_i,p_j^{\mathrm g}<1$, we define the reparametrized parameter profile $\theta$ and boundary parameters $\theta^g$ as

\be
\theta_i = \frac{p_i}{1-p_i}, \qquad \theta^g_j = \frac{p^g_j}{1-p^g_j}, \qquad i \in V, \, j \in V_g.
\ee

For ghost sites we use the convention $\theta_j:=\theta_j^{\mathrm g}$. For $S\subset V$, define
\begin{equation}
\partial\theta(S):=(\theta_j)_{j\in\partial S},
\qquad
\partial\theta(i):=\partial\theta(\{i\}).
\label{eq:exterior-scale-profile}
\end{equation}

\bd \label{def:theta_equilbrium_marginal}
Let $\vartheta \geq 0$ and $s>0$. We define
\be
\nu^{(s)}_\vartheta (n) = \frac{1}{n!} \left( \frac{\vartheta}{1 + \vartheta} \right)^{n} \frac{\Gamma(n+2s)}{\Gamma(2s)} (1+\vartheta)^{-2s}, \qquad n \in \N_0.
\ee

The measures $\nu^{(s)}_\vartheta$ are discrete Gamma with mean $2s \, \vartheta$.
\ed

We also define the discrete Gamma intertwiner associated with the discrete Gamma measures $\nu^{(s)}_\vartheta$ in Definition \ref{def:theta_equilbrium_marginal}.
 
\bd
For a function $f: \N^V_0 \rightarrow \R$ and a parameter profile $\theta \in [0, \infty)^V$, we define

\begin{equation}
(\mathcal G_s f)(\theta)
:=
\sum_{\eta\in\N_0^V}
\left(\bigotimes_{i\in V}\nu_{\theta_i}^{(s)}\right)(\eta)
f(\eta).
\label{eq:discrete-gamma-intertwiner}
\end{equation}
\ed

\subsection{Propagation of mixtures with asymmetric oriented edge weights} \label{subsec:prop_mixtures}

In \cite{giardina2025intertwining} it was shown that, when the oriented edge weights are symmetric, i.e., $a_{i,j} = a_{j,i}$ for $i,j \in V$, an inhomogeneous product measure with marginals as in Definition \ref{def:theta_equilbrium_marginal} evolves under the dynamics of the $2s$-harmonic process as if the corresponding (inhomogeneous) parameters $\theta_i$, $i \in V$, evolve according to an auxiliary continuous-time Markov process referred to as the \emph{hidden parameter process}. The proof relies heavily on duality techniques and therefore on the symmetry of the oriented edge weights. Here we revisit the result and use the decomposition in Theorem \ref{thm:decomposition} to extend it to asymmetric oriented edge weights for which the total weight directed into a given site is smaller than or equal to the total weight directed out of that site (see \eqref{eq:harmonic-outgoing-dominance}). This is the content of Theorem \ref{thm:asymmetric-harmonic-intertwining}.

\paragraph{The single-site harmonic generator}

Let $L$ denote the boundary-driven generator
for the harmonic models as in Definition \ref{def:Harmonic_models}. Denote its semigroup by $(S_t)_{t\geq0}$. For a scale profile
$\theta=(\theta_i)_{i\in V}$, set $q_i=\theta_i/(1+\theta_i)$ and use
the fixed boundary parameters for ghost neighbors. Its decomposed one-site
generator at $i$ is
\begin{equation}
\begin{aligned}
\bigl(L_i^{\partial \theta (i)}f\bigr)(n)
:= \bigl(L_{\{i\}}^{\partial p (i)}f\bigr)(n) ={}&
\sum_{k=1}^{\infty}
\frac{1}{k}
\left(\sum_{j\sim i}a_{j,i}\left(\frac{\theta_j}{1+ \theta_j}\right)^k\right)
\bigl[f(n+k)-f(n)\bigr]
\\
&+
\left(\sum_{j\sim i}a_{i,j}\right)
\sum_{k=1}^{n}\frac{1}{k}
\frac{\phi(n)}{\phi(n-k)}
\bigl[f(n-k)-f(n)\bigr].
\end{aligned}
\label{eq:asymmetric-harmonic-one-site-generator}
\end{equation}
Remark that the first line contains the incoming edge weights individually, whereas
the second depends on the outgoing coefficients only through their sum. This fact will effectively allow us to symmetrize the oriented edge weights and is at the core of the proof of Theorem \ref{thm:asymmetric-harmonic-intertwining}.

\paragraph{The sitewise symmetric comparison generators.}
In the proof of Theorem \ref{thm:asymmetric-harmonic-intertwining}, we will compare the generator $L_i^{\partial \theta (i)}$ to a symmetric single-site generator that we introduce in \eqref{eq:symmetric-one-site-comparison-generator}. Define
\begin{equation}
c_i^{\mathrm{in}}
:=
\sum_{j\sim i}a_{j,i},
\qquad
c_i^{\mathrm{out}}
:=
\sum_{j\sim i}a_{i,j},
\qquad
d_i:=c_i^{\mathrm{out}}-c_i^{\mathrm{in}},
\label{eq:harmonic-in-out-deficit}
\end{equation}
and assume that $d_i\geq0$. For $\vartheta\geq0$, let $H_i^\vartheta$ denote the
harmonic generator on one site coupled to a reservoir of scale parameter
$\vartheta$ with both oriented edge weights equal to one:
\begin{equation}
\begin{aligned}
(H_i^\vartheta f)(n)
={}&
\sum_{k=1}^{\infty}\frac{1}{k}\left(\frac{\vartheta}{1+\vartheta}\right)^k
\bigl[f(n+k)-f(n)\bigr]
\\
&+
\sum_{k=1}^{n}\frac{1}{k}
\frac{\phi(n)}{\phi(n-k)}
\bigl[f(n-k)-f(n)\bigr].
\end{aligned}
\label{eq:unit-harmonic-boundary-generator}
\end{equation}
The symmetric comparison generator associated with site $i$ is
\begin{equation}
\widetilde{L_i}^{\partial\theta(i)}
:=
\sum_{j\sim i}a_{j,i}H_i^{\theta_j}
+
d_iH_i^0.
\label{eq:symmetric-one-site-comparison-generator}
\end{equation}
We stress that this single-site generator is indeed symmetric, because for every neighboring ghost site $j$, both the incoming and the outgoing oriented edge have weight $a_{j,i}$.
Remark that this is a comparison at the level of one-site generators. The family
$\widetilde{L_i}^{\partial\theta(i)}$, $i\in V$ does not necessarily arise from a single globally symmetric particle system. The last term in
\eqref{eq:symmetric-one-site-comparison-generator} represents an
auxiliary reservoir with parameter zero and, as we will see later on, will account for the difference in total incoming weight and outgoing weight. 

\paragraph{The hidden parameter generator.}
For a function
$h:[0,\infty)^V\to\R$, define the Markov generator of the hidden process acting on site $i$ as
\begin{equation}
\begin{aligned}
(\mathscr L_i h)(\theta)
:={}&
\sum_{j\sim i}a_{j,i}
\int_0^1
\left[
h\bigl(\theta+u(\theta_j-\theta_i)\mathbf e_i\bigr)
-h(\theta)
\right]\frac{(1-u)^{2s-1}}{u}\,du
\\
&+
d_i\int_0^1
\left[
h\bigl(\theta-u\theta_i\mathbf e_i\bigr)
-h(\theta)
\right]\frac{(1-u)^{2s-1}}{u}\,du,
\end{aligned}
\label{eq:asymmetric-hidden-one-site-generator}
\end{equation}
and define the Markov generator of the boundary-driven hidden parameter process on $G$ as
\begin{equation}
\mathscr{L}
:=
\sum_{i\in V}\mathscr L_i.
\label{eq:asymmetric-hidden-generator}
\end{equation}

Denote its semigroup by $(\caS_t)_{t\geq0}$ and the hidden process by
$(\theta(t))_{t\geq0}$. The dynamics of the hidden parameter process are as follows. At rate $a_{j,i}(1-u)^{2s-1}/u$, the value of the process at site $i \in V$ is updated by replacing a portion $u$ of $\theta_i(t)$ by the value $\theta_j(t)$ at a given neighboring (ghost) site $j \sim i$, i.e $\theta_i(t) \rightarrow (1-u) \theta_i(t) + u \theta_j(t)$. The second line of
\eqref{eq:asymmetric-hidden-one-site-generator} can be interpreted as the interaction with
a ghost site whose scale parameter is zero. In the hidden parameter process only one site is updated at a time and its environment remains the same. As a consequence there is no conservation of mass.

\begin{theorem}[Propagation of mixtures for asymmetric harmonic models]
\label{thm:asymmetric-harmonic-intertwining}
Suppose that
\begin{equation}
\sum_{j\sim i}a_{j,i}
\leq
\sum_{j\sim i}a_{i,j},
\qquad i\in V.
\label{eq:harmonic-outgoing-dominance}
\end{equation}
Then, for every function $f:\N_0^V \rightarrow \R$,
\begin{equation}
\mathcal G_s\bigl(Lf\bigr)
=
\mathscr{L}\bigl(\mathcal G_s f\bigr).
\label{eq:asymmetric-harmonic-intertwining}
\end{equation}
Consequently, for all $t\geq 0$,
\begin{equation}
\mathcal G_s S_t
=
\caS_t\mathcal G_s.
\label{eq:asymmetric-harmonic-semigroup-intertwining}
\end{equation}
In particular, product measures propagate as mixtures according to
\begin{equation}
\left(\bigotimes_{i\in V}\nu_{\theta_i}^{(s)}\right)S_t
=
\E_\theta\left[
\bigotimes_{i\in V}\nu_{\theta_i(t)}^{(s)}
\right].
\label{eq:asymmetric-harmonic-propagation-of-mixtures}
\end{equation}

Here, the expectation $\E_\theta$ is with respect to the hidden process
started from $\theta$.
\end{theorem}

\begin{proof}
Apply Theorem~\ref{thm:decomposition} at the parameter profile
$(\theta_i/(1+\theta_i))_{i\in V}$. This gives
\begin{equation}
\mathcal G_s\bigl(Lf\bigr)
=
\mathcal G_s
\left(
\sum_{i\in V}L_i^{\partial \theta (i)}f
\right).
\label{eq:harmonic-application-of-decomposition}
\end{equation}
We claim that, for every $i\in V$,
\begin{equation}
L_i^{\partial \theta (i)}=\widetilde{L_i}^{\partial\theta(i)}.
\label{eq:asymmetric-comparison-generator-identity}
\end{equation}
Indeed, the insertion rate of a batch of $k$ particles under
$\widetilde{L_i}^{\partial\theta(i)}$ is
\begin{equation}
\frac{1}{k}
\left(
\sum_{j\sim i}a_{j,i}\left(\frac{\theta_j}{1+ \theta_j}\right)^k +d_i0^k
\right)
=
\frac{1}{k}\sum_{j\sim i}a_{j,i}\left(\frac{\theta_j}{1+ \theta_j}\right)^k,
\label{eq:comparison-insertion-rate}
\end{equation}
which agrees with the first line of
\eqref{eq:asymmetric-harmonic-one-site-generator}. If there are $\eta_i$ particles present at site $i$, the removal rate of $k$ particles under $\widetilde{L_i}^{\partial\theta(i)}$is

\be
\frac{1}{k} \frac{\phi(\eta_i)}{\phi(\eta_i - k)}\left( \sum_{j \sim i} a_{j,i} + d_i \right)
\ee

The term between the brackets is precisely
\begin{equation}
\sum_{j\sim i}a_{j,i}+d_i
=
c_i^{\mathrm{in}}+
\bigl(c_i^{\mathrm{out}}-c_i^{\mathrm{in}}\bigr)
=
c_i^{\mathrm{out}} = \sum_{j\sim i}a_{i,j},
\label{eq:comparison-removal-coefficient}
\end{equation}
which agrees with the second line of \eqref{eq:asymmetric-harmonic-one-site-generator}. This proves
\eqref{eq:asymmetric-comparison-generator-identity}.

The one-site harmonic intertwining from 
\cite[Theorem~3.2 and Sections~3.6--3.7]{giardina2025intertwining}
applies to each symmetric comparison generator and gives
\begin{equation}
\mathcal G_s\bigl(\widetilde{L_i}^{\partial\theta(i)} f\bigr)
=
\mathscr L_i\bigl(\mathcal G_s f\bigr).
\label{eq:symmetric-comparison-one-site-intertwining}
\end{equation}
Combining
\eqref{eq:harmonic-application-of-decomposition},
\eqref{eq:asymmetric-comparison-generator-identity}, and
\eqref{eq:symmetric-comparison-one-site-intertwining}, we find
\begin{equation}
\mathcal G_s\bigl(Lf\bigr)
=
\sum_{i\in V}\mathscr L_i(\mathcal G_s f)
=
\mathscr{L}(\mathcal G_s f),
\label{eq:asymmetric-harmonic-intertwining-proof}
\end{equation}
which proves the generator intertwining.

The semigroup identity
\eqref{eq:asymmetric-harmonic-semigroup-intertwining} follows directly from the
generator identity. Evaluating it at
the initial parameter profile $\theta$ gives
\begin{equation}
\begin{aligned}
\int
\bigl(S_tf\bigr)(\eta)
\left(\bigotimes_{i\in V}\nu_{\theta_i}^{(s)}\right)(d\eta)= \bigl(\caG_s S_t f\bigr)(\theta)=\bigl(\caS_t\mathcal G_s f\bigr)(\theta)=\E_\theta\bigl[(\mathcal G_s f)(\theta(t))\bigr],
\end{aligned}
\label{eq:asymmetric-harmonic-mixture-proof}
\end{equation}
which is exactly the propagation-of-mixtures statement
\eqref{eq:asymmetric-harmonic-propagation-of-mixtures}.
\end{proof}
\subsection{NESS for a class of harmonic models with asymmetric oriented edge weights}
\label{subsec:asymmetric-harmonic-explicit-ness}

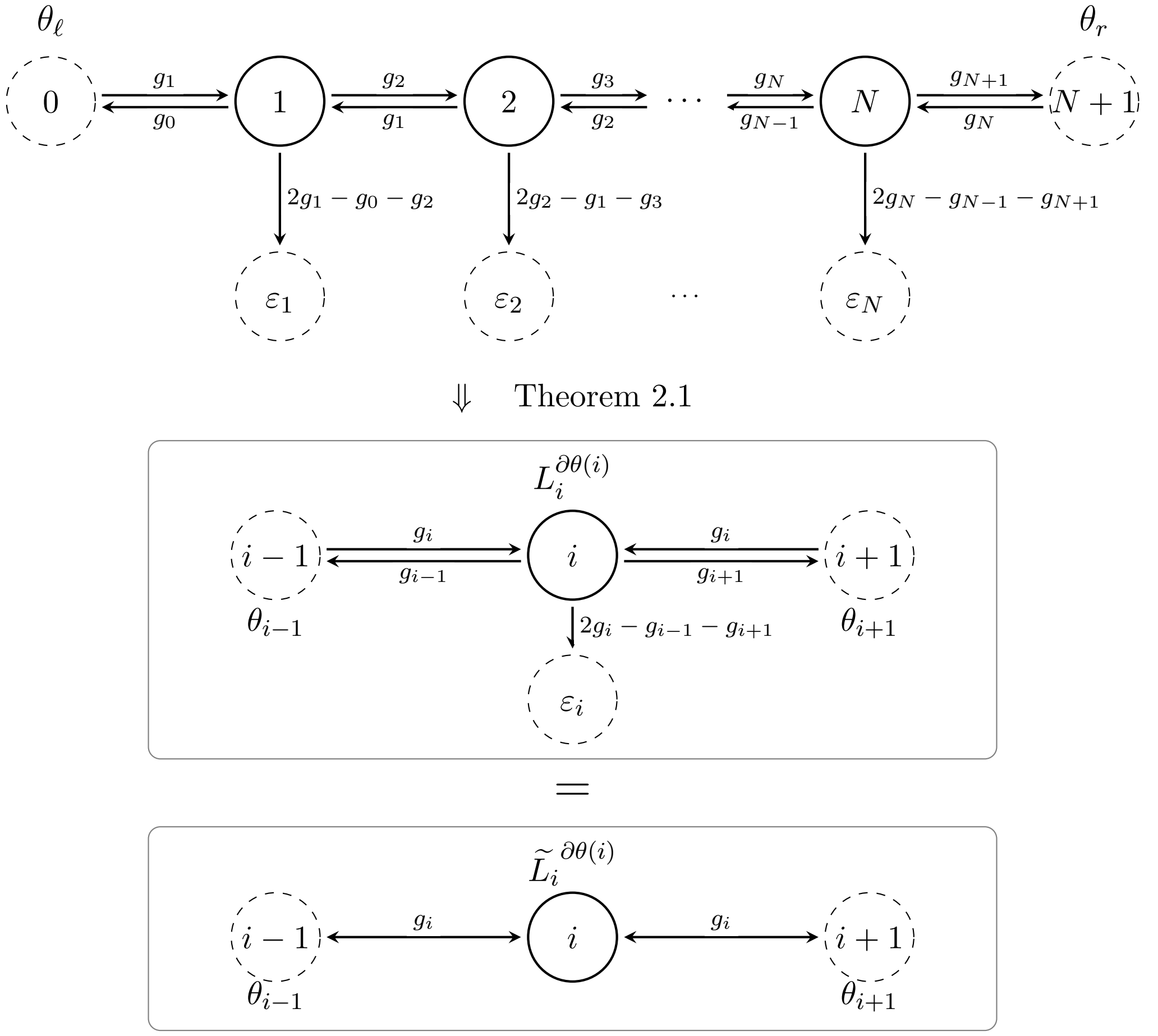
\begin{figure}[!h]
\centering
\settowidth{\dimen0}{\small $N+1$}
\edef\chainLabelWidth{\the\dimen0}
\settoheight{\dimen2}{\small $N+1$}
\edef\chainLabelHeight{\the\dimen2}
\begin{tikzpicture}[
    site/.style={circle,draw,thick,minimum size=8mm,inner sep=0pt,text width=\chainLabelWidth,
        text height=\chainLabelHeight,text depth=0pt,align=center},
    reservoir/.style={circle,draw,dashed,minimum size=8mm,inner sep=0pt,text width=\chainLabelWidth,
        text height=\chainLabelHeight,text depth=0pt,align=center},
    every node/.style={font=\small},
    rate/.style={font=\scriptsize,fill=white,inner sep=2pt},
    >=stealth
]
\node[reservoir] (gl) at (0,0) {$0$};
\node[site] (xone) at (2.7,0) {$1$};
\node[site] (xtwo) at (5.4,0) {$2$};
\node (dots) at (7.5,0) {$\cdots$};
\node[site] (xn) at (9.6,0) {$N$};
\node[reservoir] (gr) at (12.3,0) {$N+1$};

\draw[->,thick,shorten >=2pt,shorten <=2pt] ([yshift=2pt]gl.east) --
    node[rate,above] {$g_1$} ([yshift=2pt]xone.west);
\draw[->,thick,shorten >=2pt,shorten <=2pt] ([yshift=-2pt]xone.west) --
    node[rate,below] {$g_0$} ([yshift=-2pt]gl.east);
\draw[->,thick,shorten >=2pt,shorten <=2pt] ([yshift=2pt]xone.east) --
    node[rate,above] {$g_2$} ([yshift=2pt]xtwo.west);
\draw[->,thick,shorten >=2pt,shorten <=2pt] ([yshift=-2pt]xtwo.west) --
    node[rate,below] {$g_1$} ([yshift=-2pt]xone.east);
\draw[->,thick,shorten >=2pt,shorten <=2pt] ([yshift=2pt]xtwo.east) --
    node[rate,above] {$g_3$} ([yshift=2pt]dots.west);
\draw[->,thick,shorten >=2pt,shorten <=2pt] ([yshift=-2pt]dots.west) --
    node[rate,below] {$g_2$} ([yshift=-2pt]xtwo.east);
\draw[->,thick,shorten >=2pt,shorten <=2pt] ([yshift=2pt]dots.east) --
    node[rate,above] {$g_N$} ([yshift=2pt]xn.west);
\draw[->,thick,shorten >=2pt,shorten <=2pt] ([yshift=-2pt]xn.west) --
    node[rate,below] {$g_{N-1}$} ([yshift=-2pt]dots.east);
\draw[->,thick,shorten >=2pt,shorten <=2pt] ([yshift=2pt]xn.east) --
    node[rate,above] {$g_{N+1}$} ([yshift=2pt]gr.west);
\draw[->,thick,shorten >=2pt,shorten <=2pt] ([yshift=-2pt]gr.west) --
    node[rate,below] {$g_N$} ([yshift=-2pt]xn.east);

\node[font=\normalsize] at (0,0.95) {$\theta_\ell$};
\node[font=\normalsize] at (12.3,0.95) {$\theta_r$};

\node[reservoir] (envone) at (2.7,-2.3) {$\varepsilon_1$};
\node[reservoir] (envtwo) at (5.4,-2.3) {$\varepsilon_2$};
\node[reservoir] (envn) at (9.6,-2.3) {$\varepsilon_N$};
\draw[->,thick,shorten >=2pt,shorten <=2pt] (xone.south) --
    node[rate,right] {$2g_1-g_0-g_2$} (envone.north);
\draw[->,thick,shorten >=2pt,shorten <=2pt] (xtwo.south) --
    node[rate,right] {$2g_2-g_1-g_3$} (envtwo.north);
\draw[->,thick,shorten >=2pt,shorten <=2pt] (xn.south) --
    node[rate,right] {$2g_N-g_{N-1}-g_{N+1}$} (envn.north);
\node[font=\scriptsize] at (7.5,-2.3) {$\cdots$};

\node[align=center] at (6.15,-3.5) {$\Downarrow\quad\text{Theorem~\ref{thm:decomposition}}$};

\draw[gray,rounded corners] (1.15,-7.75) rectangle (11.15,-4);
\node at (6.15,-4.42) {$L_i^{\partial\theta(i)}$};
\node[reservoir] (origleft) at (2.65,-5.35) {$i-1$};
\node[site] (origsite) at (6.15,-5.35) {$i$};
\node[reservoir] (origright) at (9.65,-5.35) {$i+1$};
\node[reservoir] (origenv) at (6.15,-7.05) {$\varepsilon_i$};
\draw[->,thick,shorten >=2pt,shorten <=2pt]
    ([yshift=2pt]origleft.east) --
    node[rate,above] {$g_i$} ([yshift=2pt]origsite.west);
\draw[->,thick,shorten >=2pt,shorten <=2pt]
    ([yshift=-2pt]origsite.west) --
    node[rate,below] {$g_{i-1}$} ([yshift=-2pt]origleft.east);
\draw[->,thick,shorten >=2pt,shorten <=2pt]
    ([yshift=2pt]origright.west) --
    node[rate,above] {$g_i$} ([yshift=2pt]origsite.east);
\draw[->,thick,shorten >=2pt,shorten <=2pt]
    ([yshift=-2pt]origsite.east) --
    node[rate,below] {$g_{i+1}$} ([yshift=-2pt]origright.west);
\draw[->,thick,shorten >=2pt,shorten <=2pt]
    (origsite.south) -- node[rate,right] {$2 g_i - g_{i-1} - g_{i+1}$} (origenv.north);
\node at (2.65,-6.15) {$\theta_{i-1}$};
\node at (9.65,-6.15) {$\theta_{i+1}$};

\node[font=\Large] at (6.15,-8.15) {$=$};

\draw[gray,rounded corners] (1.15,-10.95) rectangle (11.15,-8.55);
\node at (6.15,-8.98)
    {$\widetilde{L_i}^{\partial\theta(i)}$};
\node[reservoir] (compleft) at (2.65,-9.85) {$i-1$};
\node[site] (compsite) at (6.15,-9.85) {$i$};
\node[reservoir] (compright) at (9.65,-9.85) {$i+1$};
\draw[<->,thick,shorten >=2pt,shorten <=2pt]
    (compleft.east) -- node[rate,above] {$g_i$} (compsite.west);
\draw[<->,thick,shorten >=2pt,shorten <=2pt]
    (compsite.east) -- node[rate,above] {$g_i$} (compright.west);
\node at (2.65,-10.55) {$\theta_{i-1}$};
\node at (9.65,-10.55) {$\theta_{i+1}$};
\end{tikzpicture}
\caption{The asymmetric harmonic chain for which the oriented edge weights depend only on the target site.
Solid circles are physical sites, dashed circles are ghost sites. The ghost sites different from $0$ and $N+1$ can be interpreted as reservoirs with parameter zero.
Arrow labels denote oriented edge weights rather than full
occupation-dependent transition rates. After applying $\mathcal G_s$,
the interacting generator decomposes into one-site generators in the
sense of Theorem~\ref{thm:decomposition}. The upper box displays the
original asymmetric generator at a generic site $i$. The lower box
displays its equal symmetric comparison generator.}
\label{fig:destination-weighted-harmonic-chain}
\end{figure}

In this subsection we consider the $2s$-harmonic model on the chain for which the oriented edge weights satisfy the assumptions of Theorem \ref{thm:asymmetric-harmonic-intertwining}. Moreover, they are chosen precisely so that the single-site generators obtained via the decomposition in Theorem \ref{thm:decomposition} correspond to a single-site harmonic process coupled symmetrically to two reservoirs with equal oriented edge weights. For these single-site generators, the NESS is known in closed form and has a Markovian mixture structure. Using the procedure described in Section \ref{subsec:stat_meas_Markov_mixture_structure}, we can extend the NESS to the whole graph. We find that the whole class of models has the same NESS, which is precisely the NESS of the $2s$-harmonic model on the chain found in \cite{frassek2022exact} and whose mixture structure was discovered in \cite{carinci2025large}. This is the content of Theorem \ref{thm:destination-weighted-harmonic-ness}.\\

The aforementioned class of harmonic models lives on the chain, i.e. $V=\{1,\ldots,N\}$ and $E = \{\{i,i+1\}\,:\, 1 \leq i< N \}$ with boundary reservoirs labeled
$0$ and $N+1$, whose fixed scale parameters are
$\theta_0=\theta_\ell$ and $\theta_{N+1}=\theta_r$, with
$0\leq\theta_\ell<\theta_r<\infty$. Choose positive numbers
$g_0,\ldots,g_{N+1}$ satisfying the discrete concavity condition
\begin{equation}
b_i:=2g_i-g_{i-1}-g_{i+1}\geq0,
\qquad i=1,\ldots,N.
\label{eq:concave-chain-weights}
\end{equation}
The oriented edge weights on the chain, including its boundary edges, are
\begin{equation}
a_{i,i+1}=g_{i+1},
\qquad
a_{i+1,i}=g_i,
\qquad i=0,\ldots,N.
\label{eq:destination-weighted-chain-edge weights}
\end{equation}
These oriented weights depend only on the target site. An example with $b_i = 0$ is given in Example \ref{ex:affine-destination-weights}.

In addition, particles are lost to the environment at each physical site
$i$, with harmonic rate multiplied by the oriented weight $b_i$. More precisely, if there are $n$ particles present at site, a batch of $k$ particles is removed at rate
\begin{equation}
\frac{b_i}{k}\frac{\phi(n)}{\phi(n-k)},
\qquad 1\leq k\leq n.
\label{eq:environment-harmonic-removal-rate}
\end{equation}
This loss of particles to the environment can be represented by attaching a further ghost site $\varepsilon_i$ with
parameter zero and oriented weights
$a_{i,\varepsilon_i}=b_i$ and $a_{\varepsilon_i,i}=0$.
The full particle generator is therefore
\begin{equation}
\begin{aligned}
&(Lf)(\eta)
={}\\
&\sum_{i=1}^{N}
\Biggl\{
\sum_{\substack{j\in V\\ |j-i|=1}}
\sum_{k=1}^{\eta_i}
\frac{g_j}{k}\frac{\phi(\eta_i)}{\phi(\eta_i-k)}
\bigl[f(\eta-k\mathbf e_i+k\mathbf e_j)-f(\eta)\bigr]
\\
&\hspace{25mm}
+
b_i\sum_{k=1}^{\eta_i}\frac{1}{k}
\frac{\phi(\eta_i)}{\phi(\eta_i-k)}
\bigl[f(\eta-k\mathbf e_i)-f(\eta)\bigr]
\Biggr\}
\\
&+
\sum_{k=1}^{\infty}\frac{g_1}{k}
\left(\frac{\theta_\ell}{1+\theta_\ell}\right)^k
\bigl[f(\eta+k\mathbf e_1)-f(\eta)\bigr]
+
\sum_{k=1}^{\eta_1}\frac{g_0}{k}
\frac{\phi(\eta_1)}{\phi(\eta_1-k)}
\bigl[f(\eta-k\mathbf e_1)-f(\eta)\bigr]
\\
&+
\sum_{k=1}^{\infty}\frac{g_N}{k}
\left(\frac{\theta_r}{1+\theta_r}\right)^k
\bigl[f(\eta+k\mathbf e_N)-f(\eta)\bigr]
+
\sum_{k=1}^{\eta_N}\frac{g_{N+1}}{k}
\frac{\phi(\eta_N)}{\phi(\eta_N-k)}
\bigl[f(\eta-k\mathbf e_N)-f(\eta)\bigr].
\end{aligned}
\label{eq:destination-weighted-harmonic-generator}
\end{equation}

Including the environmental ghost site, the total incoming and outgoing weights at $i$ satisfy
\begin{equation}
c_i^{\mathrm{in}}=g_i+g_i=2g_i,
\qquad
c_i^{\mathrm{out}}=g_{i-1}+g_{i+1}+b_i=2g_i.
\label{eq:destination-weighted-chain-balance}
\end{equation}
Consequently, Theorem~\ref{thm:asymmetric-harmonic-intertwining}
provides a hidden parameter process. Its generator is explicitly
\begin{equation}
\begin{aligned}
(\mathscr L^{\mathrm{hid}}h)(\theta)
={}&
\sum_{i=1}^N g_i\int_0^1
\frac{(1-u)^{2s-1}}{u}
\Bigl[
h\bigl(\theta+u(\theta_{i-1}-\theta_i)\mathbf e_i\bigr)
\\
&\hspace{38mm}+
h\bigl(\theta+u(\theta_{i+1}-\theta_i)\mathbf e_i\bigr)
-2h(\theta)
\Bigr]\,du.
\end{aligned}
\label{eq:destination-weighted-hidden-generator}
\end{equation}

Remark that there is no separate term corresponding to the loss of particles to the environmental ghost sites. This is due to the precise choice of $b_i$. For the case $N=1$, the generator in \eqref{eq:destination-weighted-hidden-generator} has a stationary measure that is Beta$(2s,\, 2s)$ (see Lemma 3.4 in \cite{giardina2025intertwining}) and can be extended to a mixing measure of the type in Proposition \ref{prop:markov_extensions_chain_continuous}. The resulting extension is a Dirichlet distribution.

\begin{example}[asymmetric oriented edge weights with $b_i = 0$]
\label{ex:affine-destination-weights}
Let $A,B>0$ such that the numbers
\begin{equation}
g_i=iA-(i-1)B,
\qquad i=0,\ldots,N+1.
\label{eq:affine-destination-weights}
\end{equation}

are greater or equal to zero. The oriented edge weights in
\eqref{eq:destination-weighted-chain-edge weights} then become
\begin{equation}
a_{i,i+1}=(i+1)A-iB,
\qquad
a_{i+1,i}=iA-(i-1)B,
\qquad i=0,\ldots,N.
\label{eq:affine-directed-edge-weights}
\end{equation}
In particular, $a_{0,1}=A$ and $a_{1,0}=B$. Since $g_i$ is
affine in $i$,
\begin{equation}
b_i=2g_i-g_{i-1}-g_{i+1}=0,
\qquad i=1,\ldots,N,
\label{eq:affine-no-environmental-loss}
\end{equation}
the environmental ghost sites have no effect and may be omitted.
These weights satisfy both
\eqref{eq:destination-weighted-chain-edge weights} and
\eqref{eq:concave-chain-weights}. Thus the difference between the rightward and leftward oriented weights is constant across every bond:
\begin{equation}
a_{i,i+1}-a_{i+1,i}=A-B.
\label{eq:affine-edge-bias}
\end{equation}
Consequently, if $A>B$, the rightward oriented edge weight is strictly larger
than the leftward one on every bond. If $B>A$, the inequality is reversed.

\begin{figure}[ht]
\centering
\begin{tikzpicture}[
    site/.style={circle,draw,minimum size=8mm,inner sep=0pt},
    reservoir/.style={circle,draw,dashed,minimum size=8mm,inner sep=0pt},
    rate/.style={font=\scriptsize,fill=white,inner sep=1.5pt}
]
\node[reservoir] (zero) at (0,0) {$0$};
\node[site] (one) at (3.1,0) {$1$};
\node[site] (two) at (6.2,0) {$2$};
\node[site] (three) at (9.3,0) {$3$};
\node[reservoir] (four) at (12.4,0) {$4$};

\draw[->,thick,shorten >=2pt,shorten <=2pt]
    ([yshift=3pt]zero.east) -- node[rate,above] {$A$}
    ([yshift=3pt]one.west);
\draw[->,thick,shorten >=2pt,shorten <=2pt]
    ([yshift=-3pt]one.west) -- node[rate,below] {$B$}
    ([yshift=-3pt]zero.east);

\draw[->,thick,shorten >=2pt,shorten <=2pt]
    ([yshift=3pt]one.east) -- node[rate,above] {$2A-B$}
    ([yshift=3pt]two.west);
\draw[->,thick,shorten >=2pt,shorten <=2pt]
    ([yshift=-3pt]two.west) -- node[rate,below] {$A$}
    ([yshift=-3pt]one.east);

\draw[->,thick,shorten >=2pt,shorten <=2pt]
    ([yshift=3pt]two.east) -- node[rate,above] {$3A-2B$}
    ([yshift=3pt]three.west);
\draw[->,thick,shorten >=2pt,shorten <=2pt]
    ([yshift=-3pt]three.west) -- node[rate,below] {$2A-B$}
    ([yshift=-3pt]two.east);

\draw[->,thick,shorten >=2pt,shorten <=2pt]
    ([yshift=3pt]three.east) -- node[rate,above] {$4A-3B$}
    ([yshift=3pt]four.west);
\draw[->,thick,shorten >=2pt,shorten <=2pt]
    ([yshift=-3pt]four.west) -- node[rate,below] {$3A-2B$}
    ([yshift=-3pt]three.east);

\node[font=\normalsize] at (0,0.78) {$\theta_\ell$};
\node[font=\normalsize] at (12.4,0.78) {$\theta_r$};
\end{tikzpicture}
\caption{The destination-weighted harmonic chain with $b_i = 0$, $i=1,2,3$, three physical sites ($N=3$). 
Each oriented edge is shown by an arrow labeled with its weight. Because $b_i=0$, the ghost sites $\epsilon_i$ are not shown.}
\label{fig:affine-destination-weighted-chain}
\end{figure}
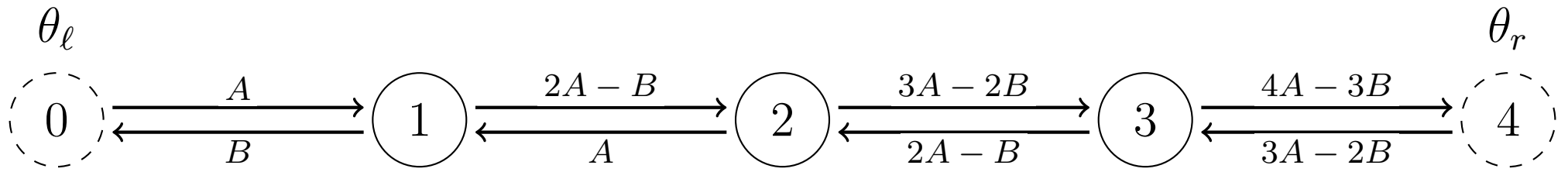
\end{example}

\begin{theorem}[Explicit stationary state of the destination-weighted chain]
\label{thm:destination-weighted-harmonic-ness}
For every $s>0$, the stationary state of
\eqref{eq:destination-weighted-harmonic-generator} is
\begin{equation}
\mu^{\theta_\ell, \theta_r}(\eta)
=
\int_{\theta_\ell\leq\theta_1\leq\cdots\leq\theta_N\leq\theta_r}
\left(\bigotimes_{i=1}^N\nu_{\theta_i}^{(s)}\right)(\eta)
\,\Lambda_{\ell,r}^{(N,s)}(d\theta),
\label{eq:destination-weighted-stationary-mixture}
\end{equation}
where, with $\theta_0=\theta_\ell$ and $\theta_{N+1}=\theta_r$,
\begin{equation}
\begin{aligned}
\Lambda_{\ell,r}^{(N,s)}(d\theta)
={}&
\frac{\Gamma(2s(N+1))}
{\Gamma(2s)^{N+1}
(\theta_r-\theta_\ell)^{2s(N+1)-1}}
\prod_{i=0}^N(\theta_{i+1}-\theta_i)^{2s-1}
\\
&\qquad\times
\1\{\theta_\ell\leq\theta_1\leq\cdots\leq\theta_N\leq\theta_r\}
\,d\theta_1\cdots d\theta_N.
\end{aligned}
\label{eq:destination-weighted-dirichlet-mixing-law}
\end{equation}
In particular, this state is independent of the weights $g_i$ and is
precisely the stationary mixture of product measures of the symmetric harmonic chain with
the same parameter $s$ and the same boundary parameters.
\end{theorem}

\begin{proof}
We use the results of Section \ref{subsubsec:Stationary_Markov_extensions}, and in particular Theorem \ref{thm:markov-mixture-stationary}. The
single-site particle generator associated with site $i\in V$ is
\begin{equation}
\begin{aligned}
\bigl(L_i^{\partial\theta(i)}f\bigr)(n)
={}&
g_i\sum_{k=1}^{\infty}\frac{1}{k}
\left[
\left(\frac{\theta_{i-1}}{1+\theta_{i-1}}\right)^k
+
\left(\frac{\theta_{i+1}}{1+\theta_{i+1}}\right)^k
\right]
\bigl[f(n+k)-f(n)\bigr]
\\
&+
2g_i\sum_{k=1}^{n}\frac{1}{k}
\frac{\phi(n)}{\phi(n-k)}
\bigl[f(n-k)-f(n)\bigr]
\\
={}&
g_i\bigl((H_i^{\theta_{i-1}}+H_i^{\theta_{i+1}})f\bigr)(n).
\end{aligned}
\label{eq:destination-weighted-single-site-identity}
\end{equation}

Indeed, both insertion weights are $g_i$, while the total removal
weight is $2g_i$ by
\eqref{eq:destination-weighted-chain-balance}. Thus this is exactly the
single-site generator of the symmetric harmonic chain, multiplied by
the positive constant $g_i$. Such multiplication changes time scales,
but not stationary measures. As a consequence, the stationary measure for each $L_i^{\partial\theta(i)}$ is of the same form, only depending on $\theta_{i-1}$ and $\theta_{i+1}$ and not on $g_i$.\\

The stationary mixing kernel for the symmetric single-site generator
$H_i^a+H_i^b$, with $a<b$, is the rescaled
$\mathrm{Beta}(2s,2s)$ law
\begin{equation}
\lambda^{a,b}(d\vartheta)
=
\frac{(\vartheta-a)^{2s-1}(b-\vartheta)^{2s-1}}
{\mathrm B(2s,2s)(b-a)^{4s-1}}
\1\{a<\vartheta<b\}\,d\vartheta.
\label{eq:destination-weighted-beta-mixing-kernel}
\end{equation}
For $a=b$, set $\lambda^{a,a}=\delta_a$.
The corresponding mixture of the marginals $\nu_\vartheta^{(s)}$ is
stationary for the single-site particle process; see
\cite[Sections~3.6--3.8]{giardina2025intertwining}.
By \eqref{eq:destination-weighted-single-site-identity}, these mixture measures are stationary for each single-site generator $L_i^{\partial\theta(i)}$.

Proposition \ref{prop:markov_extensions_chain_continuous} extends the single-site mixing measures \eqref{eq:destination-weighted-beta-mixing-kernel} to a mixing measure with a Markovian structure. This is precisely the Dirichlet measure in \eqref{eq:destination-weighted-dirichlet-mixing-law}. The stationarity of $\mu^{\theta_\ell, \theta_r}(\eta)$ in \eqref{eq:destination-weighted-stationary-mixture} follows from Theorem \ref{thm:markov-mixture-stationary}.
\end{proof}

\section{Asymmetric models with discrete-mixing NESS} \label{sec:Asymmetric models with discrete mixing structure}

In this section we introduce asymmetric generalized zero-range processes that have NESS with the mixture structure in Proposition \ref{prop:markov_extensions_chain_discrete}. We will refer to such models as having a \emph{discrete mixing structure}, because the mixing measure is supported on a discrete set of values for the parameter profile. Furthermore, we analyze the scaling behavior of these NESS and provide a large-deviation principle for the profile. The processes are introduced in Subsection \ref{subsec:asym_models}, and the large deviations of the NESS are studied in Subsection \ref{subsec:ldp}.

\subsection{Asymmetric models with discrete-mixing NESS} \label{subsec:asym_models}

The models introduced in this section are generalized zero-range processes
defined on the chain $V=\{1,2, \ldots,N\}$,
$E = \{\{i,i+1\}: 1 \leq i < N\}$, with ghost sites labeled $0$ and $N+1$
associated with sites $1$ and $N$, respectively. The ghost sites have fixed
parameters $p_0=x$ and $p_{N+1}=y$, with $x,y\in\mathsf D$. We make the
convention that $x \leq y$, i.e. the left reservoir parameter is smaller than
or equal to the right boundary parameter.

For these models, mass functions $M_{i,i \pm 1}$ are chosen such that there is an asymmetry in the rates for moving to the left and the rates for moving to the right. More precisely, let
$M_{\mathrm r},M_{\ell}:\N\to[0,\infty)$ and take, for $i \in V$ and $k \in \N$,

\be \label{eq:asym_mass_func_discrete}
M_{i,i+1}(k) := M_{\mathrm r}(k), \qquad M_{i,i-1}(k) := M_{\ell}(k).
\ee

We denote the generating functions of $M_{\mathrm r}$ and $M_{\ell}$ by
\begin{equation}
\widehat M_{\mathrm r}(q)
:=\sum_{k=1}^{\infty}M_{\mathrm r}(k)q^k,
\qquad
\widehat M_{\ell}(q)
:=\sum_{k=1}^{\infty}M_{\ell}(k)q^k,
\label{eq:asymmetric-shock-generating-functions}
\end{equation}
for $q\in\R_0^+$ such that the above series converge.

Suppose that we have a parameter profile $p \in \{x,y\}^V$ such that $p_{0} = x \leq p_1 \leq p_2 \leq \ldots \leq p_N \leq p_{N+1} = y$. Our aim is to find choices for $M_{\mathrm r},M_{\ell}$, and $\phi$ for which the single-site generator $L_i^{\partial p (i)}$, with $\partial p (i) = (p_{i-1}, p_{i+1})$ and $i \in V$, has a stationary state of the form

\be \label{eq:stat_mix_asym_one_site}
\mu^{p_{i-1},p_{i+1}}_{i}(n) = \int_{\mathsf D} \nu_{q}(n) \lambda^{p_{i-1}, p_{i+1}}_d(q) =  \frac{g_{p_{i-1}} \nu_{p_{i-1}}(n) + g_{p_{i+1}} \nu_{p_{i+1}}(n) }{g_{p_{i-1}} + g_{p_{i+1}}},
\ee

where $\lambda^{p_{i-1}, p_{i+1}}_d$ is as in \eqref{eq:discrete-single-site-kernels} and $n \in \N_0$. If this is the case, one can extend the NESS to the whole chain according to Theorem \ref{thm:markov-mixture-stationary}. In particular, the mixing measure on the chain is then the one in Proposition  \ref{prop:markov_extensions_chain_discrete}.

For $p_{i-1} = p_{i+1}$, the expression on the right-hand side of \eqref{eq:stat_mix_asym_one_site} reduces to $\nu_{p_{i-1}}$. It follows from Proposition \ref{prop:homogeneous-stationary-state} that $\nu_{p_{i-1}}$ is indeed stationary for the one-site generator. To see this, one has to realize that for the single-site model the choice in \eqref{eq:asym_mass_func_discrete} is equivalent to the choice

\be
M_{i,i+1}(k) = M_{i+1,i}(k) = M_{\ell}(k), \qquad  M_{i,i-1}(k) = M_{i-1,i}(k) = M_{\mathrm r}(k).
\ee

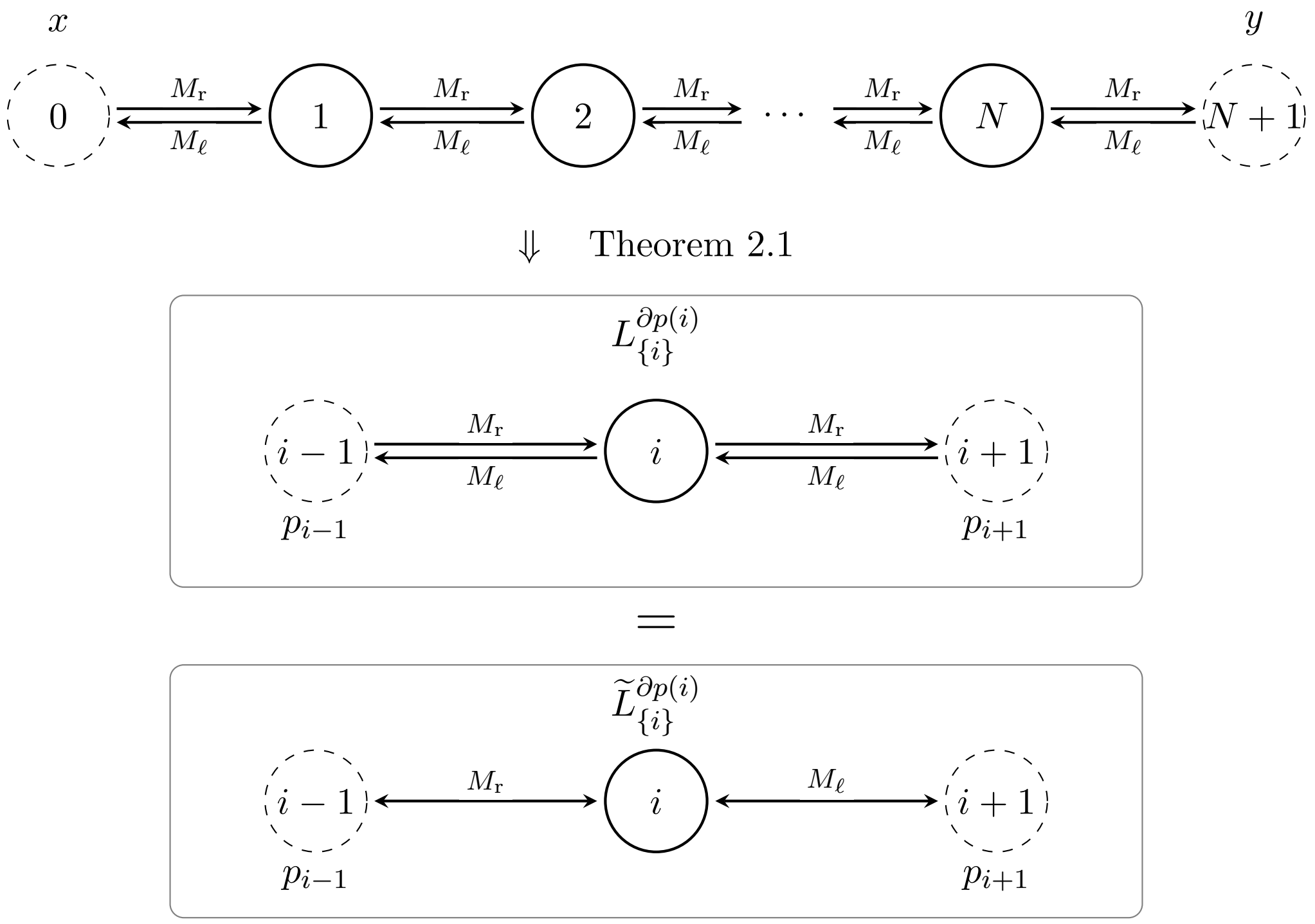
\begin{figure}[!h]
\centering
\settowidth{\dimen0}{\small $N+1$}
\edef\shockLabelWidth{\the\dimen0}
\settoheight{\dimen2}{\small $N+1$}
\edef\shockLabelHeight{\the\dimen2}
\begin{tikzpicture}[
    site/.style={circle,draw,thick,minimum size=8mm,inner sep=0pt,
        text width=\shockLabelWidth,text height=\shockLabelHeight,
        text depth=0pt,align=center},
    reservoir/.style={circle,draw,dashed,minimum size=8mm,inner sep=0pt,
        text width=\shockLabelWidth,text height=\shockLabelHeight,
        text depth=0pt,align=center},
    every node/.style={font=\small},
    rate/.style={font=\scriptsize,fill=white,inner sep=2pt},
    >=stealth
]
\node[reservoir] (shockgl) at (0,0) {$0$};
\node[site] (shockone) at (2.7,0) {$1$};
\node[site] (shocktwo) at (5.4,0) {$2$};
\node (shockdots) at (7.5,0) {$\cdots$};
\node[site] (shockn) at (9.6,0) {$N$};
\node[reservoir] (shockgr) at (12.3,0) {$N+1$};

\draw[->,thick,shorten >=2pt,shorten <=2pt]
    ([yshift=2pt]shockgl.east) -- node[rate,above] {$M_{\mathrm r}$}
    ([yshift=2pt]shockone.west);
\draw[->,thick,shorten >=2pt,shorten <=2pt]
    ([yshift=-2pt]shockone.west) -- node[rate,below] {$M_{\ell}$}
    ([yshift=-2pt]shockgl.east);
\draw[->,thick,shorten >=2pt,shorten <=2pt]
    ([yshift=2pt]shockone.east) -- node[rate,above] {$M_{\mathrm r}$}
    ([yshift=2pt]shocktwo.west);
\draw[->,thick,shorten >=2pt,shorten <=2pt]
    ([yshift=-2pt]shocktwo.west) -- node[rate,below] {$M_{\ell}$}
    ([yshift=-2pt]shockone.east);
\draw[->,thick,shorten >=2pt,shorten <=2pt]
    ([yshift=2pt]shocktwo.east) -- node[rate,above] {$M_{\mathrm r}$}
    ([yshift=2pt]shockdots.west);
\draw[->,thick,shorten >=2pt,shorten <=2pt]
    ([yshift=-2pt]shockdots.west) -- node[rate,below] {$M_{\ell}$}
    ([yshift=-2pt]shocktwo.east);
\draw[->,thick,shorten >=2pt,shorten <=2pt]
    ([yshift=2pt]shockdots.east) -- node[rate,above] {$M_{\mathrm r}$}
    ([yshift=2pt]shockn.west);
\draw[->,thick,shorten >=2pt,shorten <=2pt]
    ([yshift=-2pt]shockn.west) -- node[rate,below] {$M_{\ell}$}
    ([yshift=-2pt]shockdots.east);
\draw[->,thick,shorten >=2pt,shorten <=2pt]
    ([yshift=2pt]shockn.east) -- node[rate,above] {$M_{\mathrm r}$}
    ([yshift=2pt]shockgr.west);
\draw[->,thick,shorten >=2pt,shorten <=2pt]
    ([yshift=-2pt]shockgr.west) -- node[rate,below] {$M_{\ell}$}
    ([yshift=-2pt]shockn.east);
\node at (0,0.95) {$x$};
\node at (12.3,0.95) {$y$};

\node[align=center] at (6.15,-1.35)
    {$\Downarrow\quad\text{Theorem~\ref{thm:decomposition}}$};

\draw[gray,rounded corners] (1.15,-4.85) rectangle (11.15,-1.85);
\node at (6.15,-2.27) {$L_{\{i\}}^{\partial p(i)}$};
\node[reservoir] (shockorigleft) at (2.65,-3.45) {$i-1$};
\node[site] (shockorigsite) at (6.15,-3.45) {$i$};
\node[reservoir] (shockorigright) at (9.65,-3.45) {$i+1$};
\draw[->,thick,shorten >=2pt,shorten <=2pt]
    ([yshift=2pt]shockorigleft.east) --
    node[rate,above] {$M_{\mathrm r}$}
    ([yshift=2pt]shockorigsite.west);
\draw[->,thick,shorten >=2pt,shorten <=2pt]
    ([yshift=-2pt]shockorigsite.west) --
    node[rate,below] {$M_{\ell}$}
    ([yshift=-2pt]shockorigleft.east);
\draw[->,thick,shorten >=2pt,shorten <=2pt]
    ([yshift=2pt]shockorigsite.east) --
    node[rate,above] {$M_{\mathrm r}$}
    ([yshift=2pt]shockorigright.west);
\draw[->,thick,shorten >=2pt,shorten <=2pt]
    ([yshift=-2pt]shockorigright.west) --
    node[rate,below] {$M_{\ell}$}
    ([yshift=-2pt]shockorigsite.east);
\node at (2.65,-4.25) {$p_{i-1}$};
\node at (9.65,-4.25) {$p_{i+1}$};

\node[font=\Large] at (6.15,-5.25) {$=$};

\draw[gray,rounded corners] (1.15,-8.25) rectangle (11.15,-5.65);
\node at (6.15,-6.07) {$\widetilde L_{\{i\}}^{\partial p(i)}$};
\node[reservoir] (shockcompleft) at (2.65,-7.05) {$i-1$};
\node[site] (shockcompsite) at (6.15,-7.05) {$i$};
\node[reservoir] (shockcompright) at (9.65,-7.05) {$i+1$};
\draw[<->,thick,shorten >=2pt,shorten <=2pt]
    (shockcompleft.east) -- node[rate,above] {$M_{\mathrm r}$}
    (shockcompsite.west);
\draw[<->,thick,shorten >=2pt,shorten <=2pt]
    (shockcompsite.east) -- node[rate,above] {$M_{\ell}$}
    (shockcompright.west);
\node at (2.65,-7.85) {$p_{i-1}$};
\node at (9.65,-7.85) {$p_{i+1}$};
\end{tikzpicture}
\caption{Application of the single-site decomposition to the asymmetric
chain. Solid circles are physical sites and dashed circles are reservoirs.
Arrow labels denote mass functions rather than the full
occupation-dependent rates. The upper box is the asymmetric one-site system,
with $M_{i,j}\neq M_{j,i}$. The lower box is an equal one-site generator in
which each reservoir edge is symmetric: the left edge has mass function
$M_{\mathrm r}$ in both directions and the right edge has mass function
$M_{\ell}$ in both directions.}
\label{fig:asymmetric-shock-decomposition}
\end{figure}

The idea here is that, for the single-site generator, it does not make a difference whether particles are absorbed by the left or by the right reservoir. This idea is illustrated in Figure \ref{fig:asymmetric-shock-decomposition}.

In the case where $p_{i-1} = x \leq y = p_{i+1}$, the stationarity of $\mu^{x,y}_i$ in \eqref{eq:stat_mix_asym_one_site} for $L^{\partial p(i)}_{\{i\}}$ imposes restrictions on $M_{\ell}$, $M_{\mathrm r}$, and $\phi$. We can use the single-site stationarity conditions in Proposition \ref{prop:single-site_stationarity_equations} to define the potential function $\phi$ recursively so that $\mu^{x,y}_i$ is stationary. This recursive construction is stated in Proposition \ref{prop:asymmetric-shock-mixture}.

For the statement of Proposition \ref{prop:asymmetric-shock-mixture}, it is convenient to introduce the following notation. For $f: \N \rightarrow [0,\infty)$,

\be
\Delta_{y,x} f := f(y) - f(x).
\ee

\begin{proposition}[Asymmetric models with discrete mixing structure NESS]
\label{prop:asymmetric-shock-mixture}
Suppose that the series in
\eqref{eq:asymmetric-shock-generating-functions} are finite at $y$ and that
$x<y$. Set $\phi(0)=1$ and, for every
$\gamma\geq1$, define $\phi(\gamma)$ recursively by
\begin{equation}
(y^\gamma-x^\gamma)\frac{1}{\phi(\gamma)}
=
\sum_{k=1}^{\gamma}\frac{1}{\phi(\gamma-k)}
\left[
\frac{(y^k-x^k)M_{\mathrm r}(k)}
{\Delta_{y,x}\widehat M_{\mathrm r}}
y^{\gamma-k}
-
\frac{(y^k-x^k)M_{\ell}(k)}
{\Delta_{y,x}\widehat M_{\ell}}
x^{\gamma-k}
\right].
\label{eq:asymmetric-shock-phi-recursion}
\end{equation}
Assume that the recursion produces a potential function $\phi(\gamma)>0$ for every
$\gamma\geq1$, and that $Z(x),Z(y)<\infty$. For the generalized zero-range processes on the chain $V=\{1,\ldots,N\}$ with boundary parameters $p_0=x$
and $p_{N+1}=y$, mass functions
\begin{equation}
M_{i,i+1}(k)=M_{\mathrm r}(k),
\qquad
M_{i+1,i}(k)=M_{\ell}(k),
\qquad i=0,\ldots,N,\quad k\geq1,
\label{eq:asymmetric-shock-oriented-rates}
\end{equation}

and potential function $\phi$ as given by the recursion, we have that 

\begin{equation}
\mu^{x,y}
:=
\E\left[\bigotimes_{i=1}^N\nu_{P_i}\right]
=
\sum_{p\in\{x,y\}^N}
\Lambda_{x,y}^{(N)}(\{p\})\bigotimes_{i=1}^N\nu_{p_i}
\label{eq:asymmetric-shock-stationary-mixture}
\end{equation}
is stationary. In other words, $\mu^{x,y}$ is stationary for the boundary-driven generator
$L_V^{\partial p(V)}$. Here $P=(P_1,\ldots,P_N)$ has law $\Lambda_{x,y}^{(N)}$ as in
Proposition~\ref{prop:markov_extensions_chain_discrete}, with
\begin{equation}
g_x
:=
Z(x)\Delta_{y,x}\widehat M_{\mathrm r},
\qquad
g_y
:=
Z(y)\Delta_{y,x}\widehat M_{\ell}.
\label{eq:asymmetric-shock-mixing-weights}
\end{equation}
\end{proposition}

\begin{proof}
We verify the stationary one-site mixtures needed in
Theorem ~\ref{thm:markov-mixture-stationary}.  First suppose that
$p_{i-1}=p_{i+1}=q\in\{x,y\}$.  From
\eqref{eq:one-boundary-generator} and
\eqref{eq:asym_mass_func_discrete}, the one-site generator is
\begin{equation}
\begin{aligned}
\bigl(L_{\{i\}}^{q,q}f\bigr)(n)
={}&
\sum_{k=1}^{\infty}
\bigl(M_{\mathrm r}(k)+M_{\ell}(k)\bigr)q^k
\bigl[f(n+k)-f(n)\bigr]
\\
&+
\sum_{k=1}^{n}
\bigl(M_{\mathrm r}(k)+M_{\ell}(k)\bigr)
\frac{\phi(n)}{\phi(n-k)}
\bigl[f(n-k)-f(n)\bigr].
\end{aligned}
\label{eq:asymmetric-shock-equal-neighbor-generator}
\end{equation}
This is precisely the same generator if the two reservoir edges are replaced by the symmetric
mass functions
\begin{equation}
\begin{aligned}
\widetilde M_{i-1,i}(k)=\widetilde M_{i,i-1}(k)
&=M_{\mathrm r}(k),
\\
\widetilde M_{i,i+1}(k)=\widetilde M_{i+1,i}(k)
&=M_{\ell}(k).
\end{aligned}
\label{eq:asymmetric-shock-single-site-symmetrization}
\end{equation}
  Thus
\eqref{eq:asymmetric-shock-equal-neighbor-generator} is a symmetric
single-site generator, and $\nu_q$ is stationary by
Proposition~\ref{prop:homogeneous-stationary-state}.

It remains to consider $p_{i-1}=x$ and $p_{i+1}=y$. We substitute \eqref{eq:stat_mix_asym_one_site} for the mixture \eqref{eq:single-site-mixture-for-stationarity} into
\eqref{eq:single-site-stationarity-equations}.  At $\gamma=0$, using
$\phi(0)=1$, the stationarity equation becomes
\begin{equation}
\frac{g_x}{Z(x)}
\bigl[
\widehat M_{\ell}(x)
-\widehat M_{\ell}(y)
\bigr]
+
\frac{g_y}{Z(y)}
\bigl[
\widehat M_{\mathrm r}(y)
-\widehat M_{\mathrm r}(x)
\bigr]
=0.
\label{eq:asymmetric-shock-gamma-zero-stationarity}
\end{equation}
Consequently,
\begin{equation}
\frac{g_y}{g_x}
=
\frac{Z(y)}{Z(x)}
\frac{
\Delta_{y,x}\widehat M_{\ell}
}{
\Delta_{y,x}\widehat M_{\mathrm r}
}.
\label{eq:asymmetric-shock-weight-ratio-proof}
\end{equation}
The mixing measure $\lambda_{\mathrm d}^{x,y}$ is invariant under multiplication of
$g_x$ and $g_y$ by a common positive constant. We therefore choose $g_x$ and $g_y$ as in \eqref{eq:asymmetric-shock-mixing-weights}.

For $\gamma\geq1$, the substitution of \eqref{eq:stat_mix_asym_one_site} for the mixture measure in
\eqref{eq:single-site-stationarity-equations} yields an equation which can be rearranged into a recursion for $1 / \phi$. Combining this recursion with the choices for $g_x$ and $g_y$ in \eqref{eq:asymmetric-shock-mixing-weights} precisely gives \eqref{eq:asymmetric-shock-phi-recursion}. Through this construction, the stationarity of 
$\mu_i^{x,y}$ is immediate.

We have proved stationarity of the required one-site mixture for every increasing pair in $\{(x,x),(x,y),(y,y)\}$.  Proposition
\ref{prop:markov_extensions_chain_discrete} provides their common Markov
extension $\Lambda_{x,y}^{(N)}$, and
Theorem ~\ref{thm:markov-mixture-stationary} then gives the stationarity of
\eqref{eq:asymmetric-shock-stationary-mixture}.

\end{proof}

Proposition \ref{prop:asymmetric-shock-mixture} can be reformulated as a
statement characterizing the left and right mass functions corresponding to a
given potential function, which is given in Theorem \ref{thm:asymmetric-shock-tail-factorization}. To formulate this result, introduce random variables
$K_{\mathrm r}$ and $K_{\ell}$ with probability mass functions

\begin{equation}
m_{\mathrm r}(k)
:=
\frac{(y^k-x^k)M_{\mathrm r}(k)}{\Delta_{y,x}\widehat M_{\mathrm r}},
\qquad
m_{\ell}(k)
:=
\frac{(y^k-x^k)M_{\ell}(k)}{\Delta_{y,x}\widehat M_{\ell}},
\qquad k\geq1.
\label{eq:asymmetric-shock-normalized-mass-laws}
\end{equation}
Notice that $m_{\mathrm r}$ and $m_{\ell}$ are probability mass functions on $\N$, as they are normalized and non-negative.
We denote their probability generating functions by
\begin{equation}
\widehat m_{\mathrm r}(z)
:=
\sum_{k=1}^{\infty}m_{\mathrm r}(k)z^k,
\qquad
\widehat m_{\ell}(z)
:=
\sum_{k=1}^{\infty}m_{\ell}(k)z^k.
\label{eq:asymmetric-shock-normalized-generating-functions}
\end{equation}
The tail generating functions of $K_{\mathrm r}$ and $K_{\ell}$ are denoted by
\begin{equation}
\begin{aligned}
T_{\mathrm r}(z)
&:=
\sum_{n=0}^{\infty}\Pr(K_{\mathrm r}>n)z^n
=
\frac{1-\widehat m_{\mathrm r}(z)}{1-z},
\\
T_{\ell}(z)
&:=
\sum_{n=0}^{\infty}\Pr(K_{\ell}>n)z^n
=
\frac{1-\widehat m_{\ell}(z)}{1-z}.
\end{aligned}
\label{eq:asymmetric-shock-tail-generating-functions}
\end{equation}

\begin{theorem}[Tail-factorization characterization]
\label{thm:asymmetric-shock-tail-factorization}
Fix $0\leq x<y$ and a positive potential $\phi$ with $\phi(0)=1$ and
$Z(y)<\infty$.  Suppose that the mass functions $M_{\ell}$ and
$M_{\mathrm r}$ have generating functions $\widehat M_{\ell}$ and
$\widehat M_{\mathrm r}$ which are finite at $y$. Then the
recursion in \eqref{eq:asymmetric-shock-phi-recursion} holds if and only if
there exists a formal power series $\kappa$ with $\kappa(0)=1$ such that
\begin{equation}
T_{\mathrm r}(z)=\kappa(z)Z(xz),
\qquad
T_{\ell}(z)=\kappa(z)Z(yz).
\label{eq:asymmetric-shock-tail-factorization}
\end{equation}
If \eqref{eq:asymmetric-shock-tail-factorization} holds,
Proposition~\ref{prop:asymmetric-shock-mixture} gives that the
measure $\mu^{x,y}$ in \eqref{eq:asymmetric-shock-stationary-mixture}, with
\begin{equation}
g_x=Z(x)\Delta_{y,x}\widehat M_{\mathrm r},
\qquad
g_y=Z(y)\Delta_{y,x}\widehat M_{\ell},
\label{eq:asymmetric-shock-factorization-mixing-weights}
\end{equation}
is stationary for the boundary-driven generator
$L_V^{\partial p(V)}$
\end{theorem}

\begin{proof}
We rewrite the recursion in \eqref{eq:asymmetric-shock-phi-recursion} using generating functions. Using the notation \eqref{eq:asymmetric-shock-normalized-mass-laws}, we can rewrite \eqref{eq:asymmetric-shock-phi-recursion} as

\begin{equation}
(y^\gamma-x^\gamma)\frac{1}{\phi(\gamma)}
=
\sum_{k=1}^{\gamma}\frac{1}{\phi(\gamma-k)}
\left[
m_{\mathrm r}(k)
y^{\gamma-k}
-
m_{\ell}(k)
x^{\gamma-k}
\right].
\label{eq:asymmetric-shock-phi-recursion_prob_mass_func}
\end{equation}

Multiplying the
recursion \eqref{eq:asymmetric-shock-phi-recursion_prob_mass_func} by $z^\gamma$ and summing
over $\gamma\geq1$ gives
\begin{equation}
\bigl(1-\widehat m_{\mathrm r}(z)\bigr)Z(yz)
=
\bigl(1-\widehat m_{\ell}(z)\bigr)Z(xz).
\label{eq:asymmetric-shock-pgf-identity}
\end{equation}
After division by $1-z$, this identity says precisely that
\begin{equation}
\frac{T_{\mathrm r}(z)}{Z(xz)}
=
\frac{T_{\ell}(z)}{Z(yz)}.
\label{eq:asymmetric-shock-common-tail-factor}
\end{equation}
 Since both partition functions $Z(x)$, $Z(y)$ have constant coefficient one, their formal
power-series inverses exist.  Hence the common series in
\eqref{eq:asymmetric-shock-common-tail-factor} is the unique $\kappa$ in
\eqref{eq:asymmetric-shock-tail-factorization}. This proves both directions
of the characterization.
\end{proof}

\begin{remark}
The tail generating functions $T_\ell$ and $T_r$ determine $Z$ in the following way.  For $0<x<y$,
\eqref{eq:asymmetric-shock-common-tail-factor} can be written as
\begin{equation}
Z(yz)
=
\frac{T_{\ell}(z)}{T_{\mathrm r}(z)}Z(xz).
\label{eq:asymmetric-shock-Z-one-step}
\end{equation}
Iterating this identity $N$ times gives
\begin{equation}
Z(yz)
=
\prod_{j=0}^{N-1}
\frac{T_{\ell}\left(\left(\frac{x}{y}\right)^jz\right)}
{T_{\mathrm r}\left(\left(\frac{x}{y}\right)^jz\right)}
Z\left(y\left(\frac{x}{y}\right)^Nz\right).
\label{eq:asymmetric-shock-Z-finite-product}
\end{equation}
Since $x/y<1$ and $Z(0)=1$, letting $N\to\infty$ yields
\begin{equation}
Z(yz)
=
\prod_{j=0}^{\infty}
\frac{T_{\ell}\left(\left(\frac{x}{y}\right)^jz\right)}
{T_{\mathrm r}\left(\left(\frac{x}{y}\right)^jz\right)},
\label{eq:asymmetric-shock-Z-tail-product}
\end{equation}
whenever the product converges.  Conversely, shifting the index in the
product shows that \eqref{eq:asymmetric-shock-Z-tail-product} implies
\eqref{eq:asymmetric-shock-Z-one-step}.  Thus the two identities are
equivalent.

\end{remark}

\subsubsection{The case $\kappa(z) = 1$}
The simplest choice for $\kappa$ in Theorem~\ref{thm:asymmetric-shock-tail-factorization} is $\kappa(z)=1$. It produces a family of asymmetric generalized zero-range processes for which the total rates for moving $k$ particles to a neighboring (ghost) site are expressed in terms of Fonten\'e--Ward generalized factorials. For a sequence $(R_n)_{n\geq1}$ of positive real numbers, we define

\be
R_0!:=1,\quad R_n!:=\prod_{j=1}^{n}R_j, \quad {n\choose k}_R
:= \frac{R_n!}{R_{n-k}!R_k!} .
\ee

\begin{corollary}[Processes with Fonten\'e--Ward rates]
\label{cor:generalized-factorial-asymmetric-shock-family}
Let $(R_n)_{n\geq1}$ be a sequence satisfying
\begin{equation}
R_n\geq y,
\qquad
\sum_{n=0}^{\infty}\frac{y^n}{R_n!}<\infty,
\label{eq:generalized-factorial-asymmetric-shock-assumptions}
\end{equation}
Fix positive values of $\Delta_{y,x}\widehat M_{\mathrm r}$ and
$\Delta_{y,x}\widehat M_{\ell}$ and define, for $k\geq1$,
\begin{equation}
\begin{aligned}
M_{\mathrm r}(k)
:=
\Delta_{y,x}\widehat M_{\mathrm r} \,
\frac{x^{k-1}(R_k-x)}{R_k!(y^k-x^k)}, \quad
M_{\ell}(k)
:=
\Delta_{y,x}\widehat M_{\ell} \,
\frac{y^{k-1}(R_k-y)}{R_k!(y^k-x^k)}.
\end{aligned}
\label{eq:generalized-factorial-asymmetric-shock-masses}
\end{equation}
Then, the recursion in \eqref{eq:asymmetric-shock-phi-recursion} is solved by $\phi(n)=R_n!$, $n\in\N_0$.\\

At a site containing $n$ particles, the total rate at which a batch of
$1\leq k\leq n$ particles jumps to the right, respectively left, neighbor is given by
\begin{align}
&M_{\mathrm r}(k)\frac{\phi(n)}{\phi(n-k)}
=
\Delta_{y,x}\widehat M_{\mathrm r}
\frac{x^{k-1}(R_k-x)}{y^k-x^k}
{n\choose k}_R,\nonumber\\ 
&M_{\ell}(k)\frac{\phi(n)}{\phi(n-k)}
=
\Delta_{y,x}\widehat M_{\ell}
\frac{y^{k-1}(R_k-y)}{y^k-x^k}
{n\choose k}_R.
\label{eq:generalized-factorial-asymmetric-shock-rates}
\end{align}
The stationary measure of the boundary-driven generator $L_V^{\partial p(V)}$ is $\mu^{x,y}$ in \eqref{eq:asymmetric-shock-stationary-mixture},
with
\begin{equation}
g_x=Z(x)\Delta_{y,x}\widehat M_{\mathrm r},
\qquad
g_y=Z(y)\Delta_{y,x}\widehat M_{\ell}.
\label{eq:generalized-factorial-asymmetric-shock-weights}
\end{equation}
\end{corollary}

\begin{proof}
For fixed $t\in\{x,y\}$, define
\begin{equation}
m(t,k)
:=
\frac{t^{k-1}(R_k-t)}{R_k!}
=
\frac{t^{k-1}}{R_{k-1}!}-\frac{t^k}{R_k!}.
\label{eq:generalized-factorial-asymmetric-shock-laws}
\end{equation}
These numbers are nonnegative and telescope to
one since $t^0/R_0! = 1$.  Let $K_t$ be a random variable with probability mass function $m(t,\cdot)$.  Then
\begin{equation}
\Pr(K_t>n)=\frac{t^n}{R_n!},
\qquad
T_t(z)=\sum_{n=0}^{\infty}\frac{(tz)^n}{R_n!}=Z(tz).
\label{eq:generalized-factorial-asymmetric-shock-tails}
\end{equation}
Thus the factorization \eqref{eq:asymmetric-shock-tail-factorization} in
Theorem~\ref{thm:asymmetric-shock-tail-factorization} holds with
$\kappa(z)=1$, $m_{\mathrm r}=m(x,\cdot)$ and $m_{\ell}=m(y,\cdot)$.
Equations \eqref{eq:generalized-factorial-asymmetric-shock-masses} and
\eqref{eq:generalized-factorial-asymmetric-shock-rates} now follow
directly.
\end{proof}

\br
We remark that every potential function $\phi$ with $\phi(0):= 1$ can be written in the form $\phi(n) = R_n!$, $n \in \N_0$, as in Corollary \ref{cor:generalized-factorial-asymmetric-shock-family}. Indeed, given $\phi$, one can define $R_n = \phi(n)/\phi(n-1)$. Then

\be
\phi(n) = \frac{\phi(n)}{\phi(n-1)} \frac{\phi(n-1)}{\phi(n-2)} \ldots \frac{\phi(1)}{\phi(0)} = R_n R_{n-1} \ldots R_1 = R_n!.
\ee
\er

\begin{example}[boundary parameters $0<x<y \leq 1$]
\label{ex:four-parameter-asymmetric-shock-family}
Fix
\begin{equation}
0<x<y \leq 1,
\qquad
c>0,
\qquad
cy<1,
\qquad
A,B>0,
\label{eq:four-parameter-shock-conditions}
\end{equation}
We apply Corollary~\ref{cor:generalized-factorial-asymmetric-shock-family} with
$R_k=1/c$, $k\geq1$, and with
$\Delta_{y,x}\widehat M_{\mathrm r}=A$ and
$\Delta_{y,x}\widehat M_{\ell}=B$. 

Remark that the conditions in \eqref{eq:generalized-factorial-asymmetric-shock-assumptions} are satisfied: $R_n = 1/c \geq y$ and $\sum_{n\geq0}(cy)^n=(1-cy)^{-1}<\infty$. The corollary gives
\begin{equation}
M_{\mathrm r}(k)
:=
A\frac{(1-cx)(cx)^{k-1}}{y^k-x^k},
\qquad
M_{\ell}(k)
:=
B\frac{(1-cy)(cy)^{k-1}}{y^k-x^k}.
\label{eq:four-parameter-shock-rates}
\end{equation}
The corresponding potential is
\begin{equation}
\phi(n)=R_n!=c^{-n},
\qquad n\in\N_0.
\label{eq:four-parameter-shock-potential}
\end{equation}
Both mass functions are strictly positive for every $k\geq1$, so jumps in
both directions allow batches of every size. For a site containing $n$ particles, the total rates for a batch of
$k\leq n$ particles to jump to the right and left, respectively, are
\begin{align}
M_{\mathrm r}(k)\frac{\phi(n)}{\phi(n-k)}
&=
\frac{A(1-cx)x^{k-1}}
{c\bigl(y^k-x^k\bigr)},
\label{eq:four-parameter-shock-right-jump-rate}
\\
M_{\ell}(k)\frac{\phi(n)}{\phi(n-k)}
&=
\frac{B(1-cy)y^{k-1}}
{c\bigl(y^k-x^k\bigr)}.
\label{eq:four-parameter-shock-left-jump-rate}
\end{align}
Notice that rates are independent of $n$ apart from the constraint $k\leq n$.

For this potential,
\begin{equation}
Z(q)=\sum_{n=0}^{\infty}(cq)^n=\frac{1}{1-cq},
\label{eq:four-parameter-shock-partition-function}
\end{equation}
and the weights in \eqref{eq:asymmetric-shock-mixing-weights} are
\begin{equation}
g_x=\frac{A}{1-cx},
\qquad
g_y=\frac{B}{1-cy}.
\label{eq:four-parameter-shock-mixing-weights}
\end{equation}

\end{example}

\subsubsection{The case $x=0$ and $y=1$} \label{subsebsec:examples_x=0_y=1}

We consider the case $x=0$ and $y=1$. This naturally corresponds to a model in which the left boundary reservoir is empty and does not insert any particles into the chain, while the right boundary reservoir inserts $n$ particles with probability proportional to $1/\phi(n)$. The values $x=0$ and $y=1$ dramatically simplify the recursion for $\phi$ in Proposition \ref{prop:asymmetric-shock-mixture}. Moreover, the factorization of the tail generating functions in Theorem \ref{thm:asymmetric-shock-tail-factorization} simplifies to
$T_{\ell}(z)=T_{\mathrm r}(z)Z(z)$. In Corollary \ref{cor:asymmetric-shock-zero-one-recursion}, we provide the general form of the left and right rates for $x=0$ and $y=1$.

\begin{corollary}[General rates for $x=0$ and $y=1$]
\label{cor:asymmetric-shock-zero-one-recursion}
Fix a potential function $\phi:\N_0 \rightarrow (0,\infty)$ with
$\phi(0)=1$ and $Z(1)<\infty$, and let $m_{\mathrm r}$ be a probability
mass function on $\N$ such that
\begin{equation}
\sum_{j=1}^km_{\mathrm r}(j)\frac{1}{\phi(k-j)}
\geq \frac{1}{\phi(k)},
\qquad k\in\N.
\label{eq:requirements_phi_m}
\end{equation}

Fix $A,B>0$ and define the mass functions
\begin{equation}
M_{\mathrm r}(k):=B \, m_{\mathrm r}(k),
\qquad
M_{\ell}(k):=A \left(\sum_{j=1}^k
m_{\mathrm r}(j)\frac{1}{\phi(k-j)}-\frac{1}{\phi(k)}\right),
\qquad k\geq1.
\label{eq:zero-one-general-mass-functions}
\end{equation}
Then $\phi$ solves
the recursion in \eqref{eq:asymmetric-shock-phi-recursion}. As a consequence, the stationary measure of $L_V^{\partial p(V)}$ is the mixture in
\eqref{eq:asymmetric-shock-stationary-mixture}, with
\begin{equation}
g_0=B,
\qquad
g_1=A\sum_{n=0}^{\infty}\frac{1}{\phi(n)} = A \, Z(1).
\label{eq:zero-one-general-mixing-weights}
\end{equation}
\end{corollary}

\begin{proof}
Set
\begin{equation}
m_{\ell}(k)
:=
\sum_{j=1}^k
m_{\mathrm r}(j)\frac{1}{\phi(k-j)}
-\frac{1}{\phi(k)},
\qquad k\geq1.
\label{eq:zero-one-general-left-law}
\end{equation}
These numbers are nonnegative by assumption.  Moreover, $m_\ell$ is a probability mass function since it is normalized:
\begin{align}
\sum_{k=1}^{\infty}m_{\ell}(k)
&=
\sum_{k=1}^{\infty}\sum_{j=1}^k
m_{\mathrm r}(j)\frac{1}{\phi(k-j)}
-
\sum_{k=1}^{\infty}\frac{1}{\phi(k)}
\nonumber\\
&=
\left(\sum_{j=1}^{\infty}m_{\mathrm r}(j)\right)
\left(\sum_{l=0}^{\infty}\frac{1}{\phi(l)}\right)
-
\sum_{k=1}^{\infty}\frac{1}{\phi(k)}
=1,
\label{eq:zero-one-general-left-mass-normalization}
\end{align}
where we used $\phi(0)=1$ in the last equality. 

For $x=0$ and $y=1$, the normalized recursion
\eqref{eq:asymmetric-shock-phi-recursion_prob_mass_func} becomes
\begin{equation}
\frac{1}{\phi(\gamma)}
=
\sum_{k=1}^{\gamma}
m_{\mathrm r}(k)\frac{1}{\phi(\gamma-k)}
-m_{\ell}(\gamma),
\qquad \gamma\geq1.
\label{eq:zero-one-general-normalized-recursion}
\end{equation}
This identity is precisely the definition
\eqref{eq:zero-one-general-left-law} with $k=\gamma$, so the recursion holds.
Finally, $Z(0)=1$ and
$\Delta_{1,0}\widehat M_{\mathrm r}=B$,
$\Delta_{1,0}\widehat M_{\ell}=A$.  Therefore the weights in
\eqref{eq:asymmetric-shock-factorization-mixing-weights} reduce to
\eqref{eq:zero-one-general-mixing-weights}.
\end{proof}

\br
The condition in \eqref{eq:requirements_phi_m} is satisfied whenever

\be
m_{\mathrm r}(1)>0, \qquad \frac{\phi(k-1)}{\phi(k)}\leq m_{\mathrm r}(1),
\qquad k\geq1. 
\label{eq:zero-one-general-simple-positivity-condition}
\ee
\er

\begin{example}[Example with infinite-support mass functions]
Fix $0<c\leq1/2$ and choose
\begin{equation}
m_{\mathrm r}(k):=(1-c)c^{k-1},
\qquad
\phi(n):=c^{-n},
\qquad k\geq1,\quad n\in\N_0.
\label{eq:zero-one-infinite-support-choice}
\end{equation}
We have $m_{\mathrm r}(1)=1-c$ and
$\phi(k-1)/\phi(k)=c\leq 1-c$, so the condition
\eqref{eq:zero-one-general-simple-positivity-condition}, and as a consequence \eqref{eq:requirements_phi_m}, holds.  For
$A,B>0$, Corollary~\ref{cor:asymmetric-shock-zero-one-recursion} gives
\begin{align}
M_{\mathrm r}(k)
&=
B(1-c)c^{k-1},
\nonumber\\
M_{\ell}(k)
&=
A\left(
\sum_{j=1}^k(1-c)c^{j-1}c^{k-j}
-c^k
\right)
=
A c^{k-1}\bigl(k(1-c)-c\bigr).
\label{eq:zero-one-infinite-support-masses}
\end{align}
Since
$\phi(n)/\phi(n-k)=c^{-k}$, the rates at a site containing $n$ particles
are
\begin{equation}
M_{\mathrm r}(k)\frac{\phi(n)}{\phi(n-k)}
=B\frac{1-c}{c},
\qquad
M_{\ell}(k)\frac{\phi(n)}{\phi(n-k)}
=A\frac{k(1-c)-c}{c},
\qquad 1\leq k\leq n.
\label{eq:zero-one-infinite-support-total-rates}
\end{equation}
Notice that these rates do not depend on $n$. Moreover, the rate to the right is constant for all values of $k$, while the rate to the left depends on $k$ through the linear term $k \, A(1-c)/c$.

Also notice that $Z(1)=1/(1-c)$, and hence the mixing weights are
\begin{equation}
g_0=B,
\qquad
g_1=\frac{A}{1-c}.
\label{eq:zero-one-infinite-support-mixing-weights}
\end{equation}
Finally,
\begin{equation}
Z(z)=\frac{1}{1-cz},
\qquad
T_{\mathrm r}(z)=Z(z),
\qquad
T_{\ell}(z)=Z(z)^2.
\label{eq:zero-one-infinite-support-tail-factorization}
\end{equation}
Thus the factorization in
\eqref{eq:asymmetric-shock-tail-factorization} holds with
common factor $\kappa(z)=Z(z)$. In particular, this example is not contained
in the $\kappa(z)=1$ family of
Corollary~\ref{cor:generalized-factorial-asymmetric-shock-family}.
\end{example}

\begin{example}[Example with Fonten\'e--Ward rates]
\label{ex:fontene-ward-asymmetric-shock-rates}
In this example we recover the result of Corollary \ref{cor:generalized-factorial-asymmetric-shock-family} for $x=0$, $y=1$ by choosing $m_r(k) = \1(k=1)$ and applying Corollary \ref{cor:asymmetric-shock-zero-one-recursion}. Setting $m_r(k) = \1(k=1)$ amounts to the assumption that particles move to the right one at a time. 

Let $(R_n)_{n\geq1}$ satisfy 

\begin{equation}
R_n\geq1, \qquad \sum_{n=0}^{\infty}\frac{1}{R_n!}<\infty.
\label{eq:generalized-R-partition-finiteness}
\end{equation}

Fix $A,B>0$ and apply
Corollary~\ref{cor:asymmetric-shock-zero-one-recursion} with
\begin{equation}
\phi(n):=R_n!,
\qquad
m_{\mathrm r}(k):=\1\{k=1\}.
\label{eq:fontene-ward-potential-and-right-law}
\end{equation}
The condition \eqref{eq:requirements_phi_m} follows from
\begin{equation}
\frac{1}{R_{k-1}!}\geq\frac{1}{R_k!},
\qquad k\geq1,
\end{equation}
because $R_k\geq1$.  The corollary therefore gives
\begin{equation}
M_{\mathrm r}(k)
:=B\1\{k=1\},
\qquad
M_{\ell}(k)
:=A\left(\frac{1}{R_{k-1}!}-\frac{1}{R_k!}\right)
=A\frac{R_k-1}{R_k!},
\qquad k\geq1.
\label{eq:fontene-ward-mass-functions}
\end{equation}
For a site containing $n$ particles, the total rate at which a batch of
$1\leq k\leq n$ particles jumps to the left is
\begin{equation}
M_{\ell}(k)\frac{\phi(n)}{\phi(n-k)}
=
A(R_k-1){n\choose k}_R,
\label{eq:fontene-ward-left-rate}
\end{equation}
and the corresponding rightward rate is
\begin{equation}
M_{\mathrm r}(k)\frac{\phi(n)}{\phi(n-k)}
=
B\1\{k=1\}R_n.
\label{eq:fontene-ward-right-rate}
\end{equation}

These are precisely the rates of
Corollary~\ref{cor:generalized-factorial-asymmetric-shock-family} in \eqref{eq:generalized-factorial-asymmetric-shock-rates} with $x=0$ and $y=1$.
\end{example}

\subsubsection{Mass functions with finite support}

We consider mass transfer functions which have finite support. That is, we consider $M_\ell$ and $M_r$ for which there exist $d_\ell, \, d_r \in \N$ such that 

\be
M_r(d_r) \neq 0, \quad M_r(d_r + k) = 0 \qquad \text{for all } k \in \N,
\ee

and

\be
M_\ell(d_\ell) \neq 0, \quad M_\ell(d_\ell + k) = 0 \qquad \text{ for all } k \in \N.
\ee

We call $d_r$ (respectively $d_\ell$) the maximal point of support and write $d_r = \max \supp M_r$ (respectively $d_\ell = \max \supp M_\ell$). These mass functions are of interest because they are most similar to the mass functions of canonical zero-range models, where one particle moves at a time (i.e., $d_r, d_\ell = 1$).

\begin{proposition}[Finite-support characterization]
\label{prop:finite-support-asymmetric-shock-characterization}
The normalized mass function
$m_{\mathrm r}$ (respectively $m_{\ell}$) has finite support if and only if the corresponding tail generating function $T_{\mathrm r}$ (respectively $T_{\ell}$) is a polynomial.  More precisely, if

\begin{equation}
d_{\mathrm r}:=\max\operatorname{supp}(m_{\mathrm r}), \qquad d_{\ell}:=\max\operatorname{supp}(m_{\ell})
\end{equation}
then $T_{\mathrm r}$ has degree $d_{\mathrm r}-1$ and $T_{\ell}$ has degree $d_{\ell}-1$. The coefficients are given by 
\begin{equation}
T_{\mathrm r}(z)
=
\sum_{n=0}^{d_{\mathrm r}-1}t_{\mathrm r,n}z^n,
\quad
t_{\mathrm r,n}=\Pr(K_{\mathrm r}>n) \quad \text{and} \quad T_{\ell}(z)
=
\sum_{n=0}^{d_{\ell}-1}t_{\ell,n}z^n,
\quad
t_{\ell,n}=\Pr(K_{\ell}>n).
\label{eq:finite-support-tail-polynomial}
\end{equation}

Conversely, suppose that the coefficients satisfy
\begin{equation}
1=t_{\mathrm r,0}\geq t_{\mathrm r,1}\geq\cdots
\geq t_{\mathrm r,d_{\mathrm r}-1}>t_{\mathrm r,d_{\mathrm r}}:=0,
\qquad
1=t_{\ell,0}\geq t_{\ell,1}\geq\cdots
\geq t_{\ell,d_{\ell}-1}>t_{\ell,d_{\ell}}:=0.
\label{eq:finite-support-tail-coefficient-conditions}
\end{equation}
Then the polynomials in \eqref{eq:finite-support-tail-polynomial} are the
tail generating functions of the probability mass functions
\begin{equation}
m_{\mathrm r}(k)
=
t_{\mathrm r,k-1}-t_{\mathrm r,k},
\,\,\, 1\leq k\leq d_{\mathrm r}\quad \text{and} \quad m_{\ell}(k)
=
t_{\ell,k-1}-t_{\ell,k}, \,\,\, 1\leq k\leq d_{\ell},
\label{eq:finite-support-law-from-tail}
\end{equation}
with $m_{\mathrm r}(k)=0$ for $k>d_{\mathrm r}$ and $m_{\ell}(k)=0$ for $k>d_{\ell}$.
\end{proposition}

\begin{proof}
We give the argument for $m_r$ and $T_r$. The argument for $m_\ell$ and $T_\ell$ is identical. Suppose first that
$m_{\mathrm r}$ has finite support and set
$d_{\mathrm r}:=\max\operatorname{supp}(m_{\mathrm r})$.  Then
\begin{equation}
\Pr(K_{\mathrm r}>n)=0,
\qquad n\geq d_{\mathrm r},
\end{equation}
whereas
\begin{equation}
\Pr(K_{\mathrm r}>d_{\mathrm r}-1)
=m_{\mathrm r}(d_{\mathrm r})>0.
\end{equation}
Consequently, the series defining $T_{\mathrm r}$ terminates at degree
$d_{\mathrm r}-1$, and its coefficients are precisely
$t_{\mathrm r,n}=\Pr(K_{\mathrm r}>n)$.  This proves the forward
implication and \eqref{eq:finite-support-tail-polynomial}.

Conversely, let the coefficients satisfy
\eqref{eq:finite-support-tail-coefficient-conditions} and define
$m_{\mathrm r}$ by \eqref{eq:finite-support-law-from-tail}.  The monotonicity
of the coefficients gives $m_{\mathrm r}(k)\geq0$, while telescoping gives
\begin{equation}
\sum_{k=1}^{d_{\mathrm r}}m_{\mathrm r}(k)
=
\sum_{k=1}^{d_{\mathrm r}}
\bigl(t_{\mathrm r,k-1}-t_{\mathrm r,k}\bigr)
=t_{\mathrm r,0}-t_{\mathrm r,d_{\mathrm r}}=1.
\end{equation}
Thus $m_{\mathrm r}$ is a probability mass function.  Moreover, for
$0\leq n<d_{\mathrm r}$,
\begin{equation}
\Pr(K_{\mathrm r}>n)
=
\sum_{k=n+1}^{d_{\mathrm r}}m_{\mathrm r}(k)
=
\sum_{k=n+1}^{d_{\mathrm r}}
\bigl(t_{\mathrm r,k-1}-t_{\mathrm r,k}\bigr)
=t_{\mathrm r,n}.
\end{equation}
The tail probability vanishes for $n\geq d_{\mathrm r}$, so the associated
tail generating function is the polynomial in
\eqref{eq:finite-support-tail-polynomial}.  Finally,
$m_{\mathrm r}(d_{\mathrm r})=t_{\mathrm r,d_{\mathrm r}-1}>0$, and hence
$d_{\mathrm r}$ is indeed the maximal point of the support.
\end{proof}

\begin{example}[A finite-support model]
\label{ex:finite-support-asymmetric-shock}
Fix boundary parameters $0<x<y\leq1$ and choose the normalized mass functions
\begin{equation}
m_{\mathrm r}(1)=m_{\mathrm r}(2)=\frac12,
\qquad
m_{\ell}(1)=\frac14,
\qquad
m_{\ell}(2)=\frac58,
\qquad
m_{\ell}(3)=\frac18.
\label{eq:finite-support-example-laws}
\end{equation}
The corresponding tail generating functions are
\begin{equation}
T_{\mathrm r}(z)=1+\frac12z,
\qquad
T_{\ell}(z)=1+\frac34z+\frac18z^2
=\left(1+\frac12z\right)\left(1+\frac14z\right).
\label{eq:finite-support-example-tails}
\end{equation}
Consequently,
\begin{equation}
\frac{T_{\ell}(z)}{T_{\mathrm r}(z)}=1+\frac14z.
\end{equation}
The product representation in
\eqref{eq:asymmetric-shock-Z-tail-product} therefore gives
\begin{equation}
Z(yz)
=
\prod_{j=0}^{\infty}
\left(1+\frac14\left(\frac{x}{y}\right)^jz\right).
\label{eq:finite-support-example-partition-function}
\end{equation}
This product has strictly positive coefficients and is finite at $z=1$.
Therefore, through $Z(yz)=\sum_{n\geq0}y^nz^n/\phi(n)$, it defines an
admissible model with finite support in both directions.
\end{example}

\subsection{Large deviations for the density profile} 
\label{subsec:ldp}

In this section we discuss the large deviations of the density profiles corresponding to mixtures of product measures with discrete mixing structure. For the boundary-driven harmonic models, these large deviations were established in \cite{carinci2025large}. The proof relies heavily on the stochastic-mixture structure of the NESS and, in particular, on the availability of large deviations for the hidden parameters. The same is true for the hidden parameters of the asymmetric generalized zero-range processes in Section \ref{subsec:asym_models}. The random parameter profile satisfies a large deviation principle (LDP), and its rate function follows from a straightforward computation. In Theorem \ref{thm:asymmetric-ness-profile-ldp} we state the LDP for the empirical density profile of the asymmetric generalized zero-range models. This result is then proved using the LDP for the hidden parameter profile, which is established in the first part of the proof.\\

The stochastic mixture structure of the NESS for the asymmetric generalized zero-range models in Section \ref{subsec:asym_models}, i.e. $\mu^{x,y}$, was given in \eqref{eq:asymmetric-shock-stationary-mixture} and is repeated below,

\begin{equation}
\mu^{x,y}
=
\E\left[\bigotimes_{i=1}^N\nu_{P_i}\right]
=
\sum_{p\in\{x,y\}^N}
\Lambda_{x,y}^{(N)}(\{p\})\bigotimes_{i=1}^N\nu_{p_i}.
\label{eq:recall_asymmetric-shock-stationary-mixture}
\end{equation}

Here, $\Lambda_{x,y}^{(N)}$ is the mixing measure in \eqref{eq:discrete-chain-extension}. As a consequence, the random parameters $P_i,$ $i=1,2,\ldots,N$ take value $x$ for $i$ smaller than or equal to a random index $J$ and value $y$ for $i$ larger than $J$. In other words, the support of $\Lambda_{x,y}^{(N)}$ consists of the $N+1$ profiles, all of the form
\begin{equation}
(\underbrace{x,\ldots,x}_{J},
 \underbrace{y,\ldots,y}_{N-J}),
\qquad J=0,1,\ldots,N.
\label{eq:asymmetric-shock-profiles}
\end{equation}
We refer to this random position $J$ where the random parameters jump from $x$ to $y$ as the (microscopic) position of the \emph{shock}. There is an explicit expression for its distribution,
\begin{equation}
\P(J=j)
=
\frac{g_x^j g_y^{N-j}}
{\displaystyle\sum_{r=0}^{N}g_x^r g_y^{N-r}},
\qquad j=0,1,\ldots,N.
\label{eq:asymmetric-shock-position-law}
\end{equation}

The precise definition of the (macroscopic) shock position and the empirical density are given by 
\begin{equation}
U_N:=\frac{J}{N},
\qquad
\pi_N:=\frac1N\sum_{i=1}^N\eta_i\delta_{i/N}.
\label{eq:asymmetric-empirical-profile}
\end{equation}
We equip the space of finite positive measures on $[0,1]$ with the weak
topology. Before we state and prove the LDP for $\pi_N$, we remark that there are three regimes for the macroscopic density profile:

\begin{itemize}
    \item[1.] $g_x > g_y$: Large values of $J$ are exponentially more likely than small values of $J$. As a result the macroscopic position of the shock will concentrate near one: $U_N \rightarrow 1$  a.s. The macroscopic density profile is given by $\rho(t) = \E_{\nu_x}(\eta_1)$, $t\in (0,1)$.
    \item[2.] $g_x = g_y$: All possible configurations of the hidden parameters are equally likely. As a result the macroscopic position of the shock becomes uniformly distributed on $(0,1)$: $U_N \rightarrow \text{Unif}(0,1)$ in distribution. The macroscopic density profile is random and of the form
    \be
    \rho(t) = \E_{\nu_x}(\eta_1) \1(t \leq U) + \E_{\nu_y}(\eta_1) \1(t > U), \qquad t\in (0,1)
    \ee
    with $U \sim \text{Unif}([0,1])$. In this regime, the large deviation rate function for $U_N$, i.e., $I_{\rm sh}$ in \eqref{eq:asymmetric-joint-rate}, vanishes everywhere. In this case $U_N$ does not satisfy a large deviation principle in the canonical sense: its value does not concentrate exponentially fast on a value $u \in [0,1]$. This is because the probability that $J$ takes a given value in $\{0,1,2, \ldots,N\}$, decays as $1/(N+1)$, which is not exponentially fast in $N$.

    \item[3.] $g_x < g_y$: Small values of $J$ are exponentially more likely than large values of $J$. As a result the macroscopic position of the shock will concentrate near zero: $U_N \rightarrow 0$  a.s. The macroscopic density profile is given by $\rho(t) = \E_{\nu_y}(\eta_1)$, $t\in (0,1)$.
\end{itemize}

\begin{theorem}[Large deviations of the NESS profile]
\label{thm:asymmetric-ness-profile-ldp}
Assume that $g_x,g_y>0$ and that, for $q\in\{x,y\}$,
$Z(qe^\theta)<\infty$ for all $\theta$ in a neighborhood of zero.  Set
\begin{equation}
\psi_q(\theta):=\log Z(qe^\theta)-\log Z(q),
\qquad
h_q(a):=\sup_{\theta}\{\theta a-\psi_q(\theta)\},
\label{eq:asymmetric-local-rate}
\end{equation}
where the supremum is over the domain of $\psi_q$.  Then
$(U_N,\pi_N)$ satisfies an LDP with rate function
\begin{equation}
\mathcal I(u,\pi)=I_{\rm sh}(u)+\mathcal I_u(\pi),
\qquad
I_{\rm sh}(u):=\max\{\log g_x,\log g_y\}
-u\log g_x-(1-u)\log g_y,
\label{eq:asymmetric-joint-rate}
\end{equation}
where
\begin{equation}
\mathcal I_u(\pi):=
\sup_f\left\{
\int_0^1 f\,d\pi
-\int_0^u\psi_x(f(r))dr
-\int_u^1\psi_y(f(r))dr
\right\}.
\label{eq:asymmetric-conditional-rate-dual}
\end{equation}
The supremum is over $f\in C([0,1])$ for which the integrals are finite.
If $\pi$ is absolutely continuous, i.e. $\pi(dr)=\rho(r)dr$, 
\begin{equation}
\mathcal I_u(\rho)
=\int_0^u h_x(\rho(r))dr+
\int_u^1 h_y(\rho(r))dr.
\label{eq:asymmetric-conditional-rate-density}
\end{equation}
Consequently, the empirical particle density profile under the NESS $\mu^{x,y}$ satisfies an LDP
with rate

\begin{equation}
\mathcal I_{\mathrm{NESS}}(\pi)
=
\inf_{u\in[0,1]}
\left\{
I_{\rm sh}(u)+\mathcal I_u(\pi)
\right\}\label{eq:asymmetric-ness-profile-rate}
\end{equation}

For an absolutely continuous profile $\rho$, this becomes
\begin{equation}
\mathcal I_{\rm NESS}(\rho)=
\inf_{u\in[0,1]}
\left\{
I_{\rm sh}(u)+
\int_0^u h_x(\rho(r))dr+
\int_u^1h_y(\rho(r))dr
\right\}.
\end{equation}
\end{theorem}

\begin{remark}[Explicit form of the one-site rate]
The function $h_q$ in \eqref{eq:asymmetric-local-rate} depends on the potential function $\phi$ and can be computed as follows. For $q>0$, the substitution $z=qe^\theta$ in
\eqref{eq:asymmetric-local-rate} gives
\begin{equation}
h_q(a)=\log Z(q)+
\sup_{z:Z(z)<\infty}
\left\{a\log\frac{z}{q}-\log Z(z)\right\}.
\label{eq:asymmetric-local-rate-fugacity}
\end{equation}
The supremum is attained in $z=z(a)$, where 

\begin{equation}
a=R(z):=\frac{zZ'(z)}{Z(z)}=\E_{\nu_z}[\eta].
\label{eq:asymmetric-density-fugacity-map}
\end{equation}
Thus, whenever $a$ belongs to the range of $R$ and
$z(a):=R^{-1}(a)$,
\begin{equation}
h_q(a)=a\log\frac{z(a)}{q}-\log Z(z(a))+\log Z(q)
=\operatorname{Ent}(\nu_{z(a)}\mid\nu_q).
\label{eq:asymmetric-local-rate-explicit}
\end{equation}

Here, $\operatorname{Ent}(\nu_{z(a)}\mid\nu_q)$ denotes the relative entropy of $\nu_{z(a)}$ with respect to $\nu_q$.
\end{remark}

\begin{proof}
We first prove the LDP for the macroscopic position where the parameter profile jumps from $x$ to $y$. Then we prove the LDP for the density profile of the mixture of product measures conditioned on the parameter profile. Finally, we use the previous two steps to obtain a joint LDP for the macroscopic jump and the conditioned density profile with rate function in \eqref{eq:asymmetric-joint-rate}, which then leads to the LDP for the unconditioned density profile with rate function in \eqref{eq:asymmetric-ness-profile-rate}.\\

\noindent
\underline{LDP for $J$:} For $u\in[0,1]$,
\begin{align}
-\frac{1}{N} \log \P(J = \lfloor u N \rfloor) &= -\frac1N\log
\frac{g_x^{\lfloor uN\rfloor}g_y^{N-\lfloor uN\rfloor}}
{\displaystyle\sum_{r=0}^Ng_x^rg_y^{N-r}}
\nonumber\\
&=-\frac{\lfloor uN\rfloor}{N}\log g_x
-\left(1-\frac{\lfloor uN\rfloor}{N}\right)\log g_y
+\frac1N\log\left(\sum_{r=0}^Ng_x^rg_y^{N-r}\right).
\label{eq:asymmetric-shock-point-ldp}
\end{align}

Taking the limit $N \rightarrow \infty$ yields $I_{\rm sh}(u)$. Here we use that the sum in the second line lies between
$\max\{g_x,g_y\}^N$ and $(N+1)\max\{g_x,g_y\}^N$.\\

\noindent
\underline{LDP for $\pi_N$:} We condition on $J=j_N$, where $J_N$ is such that $j_N/N\to u$. Given $J_N$, the variables $\eta_i$, $i =1,2, \ldots,N$, are independent, with distribution $\nu_x$ for $i\leq j_N$ and
distribution $\nu_y$ for $i>j_N$. The logarithmic generating function for a single site can be computed directly from the definition of the single-site marginal in \eqref{eq:equilibrium-marginals}. For $q \in \{x,y\}$,
\begin{equation}
\log\E_{\nu_q}[e^{\theta\eta}]
=\log\left(\frac1{Z(q)}
\sum_{n=0}^{\infty}\frac{(qe^\theta)^n}{\phi(n)}\right)
=\psi_q(\theta).
\label{eq:asymmetric-single-site-log-mgf}
\end{equation}
Therefore, for every admissible $f\in C([0,1])$, we have that
\begin{align}
&\frac1N\log\E\left[
e^{N\langle f,\pi_N\rangle}\mid J=j_N
\right]\nonumber\\
&\qquad=
\frac1N\sum_{i=1}^{j_N}\psi_x(f(i/N))
+\frac1N\sum_{i=j_N+1}^N\psi_y(f(i/N))
\longrightarrow
\int_0^u\psi_x(f(r))dr+
\int_u^1\psi_y(f(r))dr.
\label{eq:asymmetric-conditional-log-mgf}
\end{align}

Here we use that, conditional on $J$, the variables $\eta_i$ are independent. Applying the G\"artner--Ellis theorem yields the conditional LDP with rate function \eqref{eq:asymmetric-conditional-rate-dual}.  The assumed exponential moment
gives exponential tightness.  For absolutely continuous profiles,
interchanging the pointwise supremum and the integral gives
\eqref{eq:asymmetric-conditional-rate-density}.\\

\noindent
\underline{Joint LDP and contraction:} We combine the LDPs for $J$ and $\pi_N$ into a joint LDP.  For an open set $\mathcal O$, fix
$(u,\pi)\in\mathcal O$, choose $J_N$ such that $j_N/N\to u$, and choose a neighborhood
$B$ of $\pi$ such that $(j_N/N,B)\subset\mathcal O$ for sufficiently large $N$.
Then
\begin{align}
\liminf_{N\to\infty}\frac1N
\log\P\bigl((U_N,\pi_N)\in\mathcal O\bigr)
&\geq\liminf_{N\to\infty}\left\{
\frac1N\log\P(J=j_N)
+\frac1N\log\P(\pi_N\in B\mid J=j_N)\right\}
\nonumber\\
&\geq-I_{\rm sh}(u)-\inf_{\widetilde\pi\in B}
\mathcal I_u(\widetilde\pi)
\nonumber\\
&\geq-I_{\rm sh}(u)-\mathcal I_u(\pi).
\label{eq:asymmetric-joint-ldp-lower-bound}
\end{align}
The last inequality uses $\pi\in B$, and hence
$\inf_{\widetilde\pi\in B}\mathcal I_u(\widetilde\pi)
\leq\mathcal I_u(\pi)$.
Taking the supremum over $(u,\pi)\in\mathcal O$ gives the lower bound.
If $\mathcal F$ is closed and
$\mathcal F_v:=\{\pi:(v,\pi)\in\mathcal F\}$, then
\begin{align}
\P\bigl((U_N,\pi_N)\in\mathcal F\bigr)
&=\sum_{j=0}^N\P(J=j)
\P(\pi_N\in\mathcal F_{j/N}\mid J=j)
\nonumber\\
&\leq(N+1)\max_{0\leq j\leq N}
\left\{\P(J=j)
\P(\pi_N\in\mathcal F_{j/N}\mid J=j)\right\},
\label{eq:asymmetric-joint-ldp-upper-bound}
\end{align}

Hence,

\be
\limsup_{N\to\infty}\frac1N
\log\P\bigl((U_N,\pi_N)\in\mathcal F\bigr)
\leq-\inf_{(u,\pi)\in\mathcal F}
\{I_{\rm sh}(u)+\mathcal I_u(\pi)\}.
\ee

This
proves the joint LDP.  Finally, the contraction principle under the
projection $(u,\pi)\mapsto\pi$ gives
\eqref{eq:asymmetric-ness-profile-rate}.
\end{proof}

\begin{example}[Total density for $g_x = g_y$]
\label{ex:equal-weight-total-density-ldp}
Suppose that $g_x=g_y$. We apply the contraction principle to the joint LDP for $(U_N,\pi_N)$ in Theorem \ref{thm:asymmetric-ness-profile-ldp} to compute the rate function for the total density

\begin{equation}
\overline\eta_N
:=
\left\langle 1,\pi_N\right\rangle
=
\frac1N\sum_{i=1}^N\eta_i.
\label{eq:equal-weight-total-density}
\end{equation}
The map $\pi\mapsto\pi([0,1])$ is continuous for the weak topology
on finite measures on $[0,1]$ and the contraction principle gives
\begin{equation}
\mathcal I_{\mathrm{tot}}(r)
=
\inf_{\substack{u\in[0,1],\,\pi\\ \pi([0,1])=r}}
\left\{
I_{\rm sh}(u)+\mathcal I_u(\pi)
\right\}.
\label{eq:equal-weight-total-density-contraction}
\end{equation}
When $g_x=g_y$, the rate function for the shock-position, i.e. $I_{\rm sh}$, vanishes for all values of $u \in [0,1]$.
For fixed $u$, contraction of $\mathcal I_u$ to the total mass gives
\begin{equation}
J_u(r)
:=
\inf_{\pi:\,\pi([0,1])=r}\mathcal I_u(\pi)
=
\sup_{\theta}
\left\{
\theta r-u\psi_x(\theta)-(1-u)\psi_y(\theta)
\right\},
\label{eq:equal-weight-fixed-shock-rate}
\end{equation}
where the supremum is over the values of $\theta$ for which $\psi_x$ and $\psi_y$
are finite. Hence
\begin{equation}
\mathcal I_{\mathrm{tot}}(r)
=
\inf_{u\in[0,1]}J_u(r).
\end{equation}

Set
\begin{equation}
\rho_x:=R(x)=\E_{\nu_x}[\eta],
\qquad
\rho_y:=R(y)=\E_{\nu_y}[\eta],
\end{equation}
where $R$ is defined in \eqref{eq:asymmetric-density-fugacity-map}. For fixed shock position $u$, the typical total density is $u \rho_x + (1-u)\rho_y$, i.e.  

\be \label{eq:J_u_vanishes}
J_u(u \rho_x + (1-u) \rho_y) = 0.
\ee

Every $r$ in the interval $[\rho_x,\rho_y]$ is of the form $r = u_r \rho_x + (1-u_r) \rho_y$, for some $u_r \in [0,1]$. Hence \eqref{eq:J_u_vanishes} gives

\be
0 \leq \mathcal I_{\mathrm{tot}}(r)
=
\inf_{u\in[0,1]}J_u(r) \leq J_{u_r}(r) = J_{u_r} (u_r \rho_x + (1-u_r) \rho_y) = 0.
\ee

We conclude that the rate function $\mathcal I_{\mathrm{tot}}(r)$ vanishes for any total density $r$ in the interval $[\rho_x, \rho_y]$.

Next we identify the rate outside of the interval $[\rho_x, \rho_y]$. We will find the following. If $r \leq \rho_x$, then the value of $J_u(r)$ will be minimal in $u=1$, for which $J_u(r) = h_x(r)$. Similarly, if $r \geq \rho_y$, then the value of $J_u(r)$ will be minimal in $u =0$, for which $J_u(r) = h_y(r)$. To prove this, we use the following observation,

\begin{align} \label{eq:psi_ineq}
&\psi_x(\theta) \leq \psi_y(\theta)  \quad \text{for } \theta \geq 0, \nonumber\\
&\psi_x(\theta) \geq \psi_y(\theta)  \quad \text{for } \theta \leq 0,
\end{align}

in combination with

\begin{align} \label{eq:diff_sup}
&h_x(r) = \sup_{\theta \leq 0} \{\theta r - \psi_x(\theta) \},
\qquad r\leq\rho_x,\nonumber\\
&h_y(r) = \sup_{\theta \geq 0} \{\theta r - \psi_y(\theta) \},
\qquad r\geq\rho_y.
\end{align}

Remark that the suprema in \eqref{eq:diff_sup} are taken over $\theta \leq 0$ and $\theta \geq 0$ respectively. We only give the proof for the expression for $h_x$. The expression for $h_y$ can be shown using a similar argument. Suppose that $\theta \geq 0$. Using the convexity of $\psi_x$ we find

\be
\psi_x(\theta) \geq \psi_x(0) + \theta \psi_x'(0) = \theta \rho_x
\ee

and hence

\be
\theta r -\psi_x(\theta) \leq \theta (r - \rho_x) \leq 0.
\ee

Because $\theta r - \psi_x(\theta)$ is precisely equal to zero for $\theta = 0$, we conclude that the supremum in \eqref{eq:asymmetric-local-rate} can be replaced by a supremum over only $\theta \leq 0$. This gives the expression for $h_x$ in \eqref{eq:diff_sup}.

Now we see that, for $r \leq \rho_x$,

\begin{align}
    J_u(r) &= \sup_{\theta}
\left\{
\theta r-u\psi_x(\theta)-(1-u)\psi_y(\theta)
\right\}\nonumber\\
&\geq \sup_{\theta \leq 0} \left\{
\theta r-u\psi_x(\theta)-(1-u)\psi_y(\theta)
\right\}\nonumber\\
&\geq \sup_{\theta \leq 0} \left\{
\theta r-\psi_x(\theta)
\right\}\nonumber = h_x(r).\\
\end{align}

For the last inequality, we use \eqref{eq:psi_ineq}. It follows that $\mathcal I_{\mathrm{tot}}(r) = J_1(r) = h_x(r)$. The proof for $\mathcal I_{\mathrm{tot}}(r) = J_0(r) = h_y(r)$ is analogous. 

We conclude that the large deviation rate function for the total density is given by

\begin{equation}
\mathcal I_{\mathrm{tot}}(r)
=
\begin{cases}
+\infty, & r<0,\\[1mm]
h_x(r), & 0\leq r<\rho_x,\\[1mm]
0, & \rho_x\leq r\leq\rho_y,\\[1mm]
h_y(r), & r>\rho_y.
\end{cases}
\label{eq:equal-weight-total-density-rate}
\end{equation}
The flat part of the rate function reflects that, when $g_x=g_y$, moving
the macroscopic shock position has no cost on the exponential scale.
\end{example}

\subsubsection*{Acknowledgments:} I want to thank Frank Redig and Cristian Giardinà for useful and inspiring discussions about the content of this paper. I also want to thank the Peter Paul Peterich Foundation for funding this research via TU Delft
University Fund. 

\bibliographystyle{plain}
\bibliography{bibliography}

@article{derrida1998exactly,
  title={An exactly soluble non-equilibrium system: the asymmetric simple exclusion process},
  author={Derrida, Bernard},
  journal={Physics Reports},
  volume={301},
  number={1-3},
  pages={65--83},
  year={1998},
  publisher={Elsevier}
}

@article{derrida2002large,
  title={Large deviation of the density profile in the steady state of the open symmetric simple exclusion process},
  author={Derrida, Bernard and Lebowitz, Joel L. and Speer, Eugene R.},
  journal={Journal of Statistical Physics},
  volume={107},
  number={3--4},
  pages={599--634},
  year={2002},
  doi={10.1023/A:1014555927320},
  publisher={Springer}
}

@article{derrida2003exact,
  title={Exact large deviation functional of a stationary open driven diffusive system: The asymmetric exclusion process},
  author={Derrida, Bernard and Lebowitz, Joel L. and Speer, Eugene R.},
  journal={Journal of Statistical Physics},
  volume={110},
  number={3--6},
  pages={775--810},
  year={2003},
  doi={10.1023/A:1022128811691},
  publisher={Springer}
}

@article{enaud2004large,
  title={Large deviation functional of the weakly asymmetric exclusion process},
  author={Enaud, Camille and Derrida, Bernard},
  journal={Journal of Statistical Physics},
  volume={114},
  number={3--4},
  pages={537--562},
  year={2004},
  doi={10.1023/B:JOSS.0000012501.43746.cf},
  publisher={Springer}
}

@article{frassek2022exact,
  title={Exact solution of an integrable non-equilibrium particle system},
  author={Frassek, Rouven and Giardin{\`a}, Cristian},
  journal={Journal of Mathematical Physics},
  volume={63},
  number={10},
  year={2022},
  publisher={AIP Publishing}
}

@article{frassek2020non,
  title={Non-compact quantum spin chains as integrable stochastic particle processes},
  author={Frassek, Rouven and Giardin{\`a}, Cristian and Kurchan, Jorge},
  journal={Journal of Statistical Physics},
  volume={180},
  number={1},
  pages={135--171},
  year={2020},
  publisher={Springer}
}

@article{carinci2025large,
  title={Large deviations and additivity principle for the open harmonic process},
  author={Carinci, Gioia and Franceschini, Chiara and Frassek, Rouven and Giardin{\`a}, Cristian and Redig, Frank},
  journal={Communications in Mathematical Physics},
  volume={406},
  number={5},
  pages={103},
  year={2025},
  publisher={Springer}
}

@article{redig2025non,
  title={Non-equilibrium Steady States with a Spatial Markov Structure},
  author={Redig, Frank and Tol, Berend van},
  journal={Journal of Statistical Physics},
  volume={192},
  number={7},
  pages={89},
  year={2025},
  publisher={Springer},
  doi={10.1007/s10955-025-03471-4}
}

@article{giardina2025intertwining,
  title={Intertwining and propagation of mixtures for generalized KMP models and harmonic models},
  author={Giardin{\`a}, Cristian and Redig, Frank and van Tol, Berend},
  journal={Journal of Statistical Physics},
  volume={192},
  number={2},
  pages={21},
  year={2025},
  publisher={Springer}
}

@article{carinci2024solvable,
  title={Solvable stationary non equilibrium states},
  author={Carinci, Gioia and Franceschini, Chiara and Gabrielli, Davide and Giardin{\`a}, Cristian and Tsagkarogiannis, Dimitrios},
  journal={Journal of Statistical Physics},
  volume={191},
  number={1},
  pages={10},
  year={2024},
  publisher={Springer}
}

@article{greenblatt2006product,
  title={Product measure steady states of generalized zero range processes},
  author={Greenblatt, RL and Lebowitz, JL},
  journal={Journal of Physics A: Mathematical and General},
  volume={39},
  number={7},
  pages={1565--1573},
  year={2006}
}

@article{evans2006factorized,
  title={Factorized steady states in mass transport models on an arbitrary graph},
  author={Evans, Martin R and Majumdar, Satya N and Zia, Royce KP},
  journal={Journal of Physics A: Mathematical and General},
  volume={39},
  number={18},
  pages={4859--4873},
  year={2006}
}

@article{fajfrova2016invariant,
  title={Invariant measures of mass migration processes},
  author={Fajfrov{\'a}, Lucie and Gobron, Thierry and Saada, Ellen},
  journal={Electronic Journal of Probability},
  volume={21},
  number={60},
  pages={1--52},
  year={2016},
  doi={10.1214/16-EJP4399}
}

@misc{redig2026ergodicpropertiesharmonicprocess,
      title={Ergodic properties of the harmonic process}, 
      author={Frank Redig and Berend van Tol},
      year={2026},
      eprint={2604.17469},
      archivePrefix={arXiv},
      primaryClass={math.PR},
      url={https://arxiv.org/abs/2604.17469}, 
}

@article{levine2005zero,
  title={Zero-range process with open boundaries},
  author={Levine, E. and Mukamel, D. and Sch\"utz, G. M.},
  journal={Journal of Statistical Physics},
  volume={120},
  number={5--6},
  pages={759--778},
  year={2005},
  doi={10.1007/s10955-005-7000-7}
}

@article{capanna2024mass,
  title={On a Class of Solvable Stationary Non Equilibrium States for Mass Exchange Models},
  author={Capanna, Monia and Gabrielli, Davide and Tsagkarogiannis, Dimitrios},
  journal={Journal of Statistical Physics},
  volume={191},
  number={2},
  pages={25},
  year={2024},
  doi={10.1007/s10955-024-03240-9}
}

@article{frassek2024madminess,
  title={The steady state of the boundary-driven multiparticle asymmetric diffusion model},
  author={Frassek, Rouven and Sz\'ecs\'enyi, Istv\'an M.},
  journal={Journal of Physics A: Mathematical and Theoretical},
  volume={57},
  number={9},
  pages={095205},
  year={2024},
  doi={10.1088/1751-8121/ad2725}
}

@article{balazs2010shock,
  title={Random Walk of Second Class Particles in Product Shock Measures},
  author={Bal\'azs, M\'arton and Farkas, Gy\"orgy and Kov\'acs, P\'eter and R\'akos, Attila},
  journal={Journal of Statistical Physics},
  volume={139},
  number={2},
  pages={252--279},
  year={2010},
  doi={10.1007/s10955-010-9933-8}
}

@article{balazs2019qzerorange,
  title={$q$-Zero Range has Random Walking Shocks},
  author={Bal\'azs, M\'arton and Duffy, Lewis and Pantelli, Dimitri},
  journal={Journal of Statistical Physics},
  volume={174},
  pages={958--971},
  year={2019},
  doi={10.1007/s10955-018-02218-8}
}

@article{gabrielli2026markovbridges,
  title={Mixtures, Markov Bridges and the Matrix Product Ansatz},
  author={Gabrielli, Davide and Iacovissi, Federica},
  journal={Annales Henri Poincar\'e},
  year={2026},
  doi={10.1007/s00023-026-01742-8}
}

@misc{gabrielli2026matrixldp,
  title={Large deviations and the matrix product ansatz},
  author={Gabrielli, Davide and Iacovissi, Federica},
  year={2026},
  eprint={2601.11820},
  archivePrefix={arXiv},
  primaryClass={math.PR},
  url={https://arxiv.org/abs/2601.11820}
}

@article{frassek2026representations,
  title={Steady state representations for the harmonic process},
  author={Frassek, Rouven},
  journal={Journal of Physics A: Mathematical and Theoretical},
  volume={59},
  pages={085202},
  year={2026},
  doi={10.1088/1751-8121/ae4799}
}

@article{bernardin2022long,
  title={Non-equilibrium Stationary Properties of the Boundary Driven Zero-Range Process with Long Jumps},
  author={Bernardin, C\'edric and Gon\c{c}alves, Patr\'icia and Jim\'enez-Oviedo, Byron and Scotta, Stefano},
  journal={Journal of Statistical Physics},
  volume={189},
  number={3},
  pages={32},
  year={2022},
  doi={10.1007/s10955-022-02987-3}
}
\end{document}